\documentclass{amsart}
\usepackage{latexsym,amsxtra,amscd,ifthen}
\usepackage[normalem]{ulem}
\usepackage{cleveref}
\usepackage{amsfonts}
\usepackage{color,graphicx}
\usepackage{verbatim}
\usepackage{amsmath}
\usepackage{amsthm}
\usepackage{amssymb}
\usepackage{url}
\usepackage{todonotes}

\usepackage{tikz}
\usepackage{multirow}
\usepackage{float}

\newtheorem{theorem}{Theorem}
\newtheorem{lemma}[theorem]{Lemma}
\newtheorem{proposition}[theorem]{Proposition}
\newtheorem{construction}[theorem]{Construction}
\newtheorem{hypothesis}[theorem]{Hypothesis}
\newtheorem{corollary}[theorem]{Corollary}

\numberwithin{theorem}{section}
\numberwithin{equation}{theorem}

\theoremstyle{definition}
\newtheorem{definition}[theorem]{Definition}
\newtheorem{example}[theorem]{Example}
\newtheorem{remark}[theorem]{Remark}
\newtheorem{question}[theorem]{Question}
\newtheorem*{question*}{Question}

\newcommand{\fm}{\mathfrak{m}}

\newcommand{\fp}{\mathfrak{p}}

\newcommand{\im}{\text{im}}

\newcommand{\bfw}{\mathbf{w}}
\newcommand{\bfd}{\mathbf{d}}

\DeclareMathOperator{\height}{ht}

\DeclareMathOperator{\rk}{rk}

\DeclareMathOperator{\Proj}{Proj}

\DeclareMathOperator{\Aut}{Aut}
\DeclareMathOperator{\gr}{gr}

\DeclareMathOperator{\hdet}{hdet}
\DeclareMathOperator{\Pdet}{Pdet}
\DeclareMathOperator{\GKdim}{GKdim}
\DeclareMathOperator{\gldim}{gldim}

\DeclareMathOperator{\MaxReg}{MaxReg}

\DeclareMathOperator{\Reg}{Reg}
\DeclareMathOperator{\Spec}{Spec}

\newcommand\kk{{\Bbbk}}

\newcommand{\jz}

\begin{document}

\title{Poisson valuations}

\author{Hongdi Huang, Xin Tang, Xingting Wang, and James J. Zhang}

\address{Huang: Department of Mathematics, Shanghai University,
Shanghai, 200444, China}

\email{hdhuang@shu.edu.cn}

\address{Tang: Department of Mathematics \& Computer Science,
Fayetteville State University, Fayetteville, NC 28301,
USA}

\email{xtang@uncfsu.edu}

\address{Wang: Department of Mathematics, 
Louisiana State University, Baton Rouge, Louisiana 70803, USA}

\email{xingtingwang@math.lsu.edu}

\address{Zhang: Department of Mathematics, Box 354350,
University of Washington, Seattle, Washington 98195, USA}
\email{zhang@math.washington.edu}

\begin{abstract}
We study Poisson valuations and provide their applications in solving problems related to rigidity, automorphisms, Dixmier property, isomorphisms, and embeddings of Poisson algebras and fields.
\end{abstract}





\maketitle


\section*{Introduction}
\label{xxsec0}
Poisson algebras are ubiquitous and an essential area of study with great significance. Various aspects of Poisson algebras have been researched in the works \cite{Ba1, Ba2, Go1, Go2, GLa, GLe, JO, LS, LuWW1, LuWW2, LvWZ, PS}, focusing on (twisted) Poincar{\' e} duality and modular derivations, Poisson Dixmier-Moeglin equivalence, Poisson enveloping algebras, Poisson spectrum and more. The role of Poisson algebras has been featured in the study of the representation theory of PI Sklyanin algebras \cite{WWY1, WWY2} and discriminants of noncommutative algebras \cite{BY1, BY2, LY1, NTY}. Additionally, the isomorphism problem and cancellation problem have been studied in the context of Poisson algebras \cite{GVW, GW, GWY}.

This paper presents the concept of Poisson valuations defined on a Poisson field. Poisson valuations can be perceived as an invariant of Poisson fields (or Poisson algebras) analogous to a prime divisor in algebraic geometry. We explore their applications by solving problems concerning rigidity, automorphisms, Dixmier property, isomorphisms, and embeddings of Poisson algebras and fields.

\subsection{Definitions}
\label{xxsec0.1}
Let $\Bbbk$ be a base field. All algebraic objects are over 
$\Bbbk$. Starting with Hypothesis \ref{xxhyp3.5}, we assume throughout that $\Bbbk$ is of characteristic zero. 

\begin{definition}
\label{xxdef0.1}
Let $K$ be a Poisson field. A map $\nu: K\to {\mathbb Z}\cup\{\infty\}$
is called a {\it Poisson valuation} on $K$ if, for all 
$a,b\in K$,
\begin{enumerate}
\item[(1)]
$\nu(a)=\infty$ if and only if $a=0$,
\item[(2)]
$\nu(a)=0$ for all $a\in \Bbbk^{\times}:=\Bbbk\setminus \{0\}$,
\item[(3)]
$\nu(ab)=\nu(a)+\nu(b)$,
\item[(4)]
$\nu(a+b)\geq \min\{\nu(a),\nu(b)\}$, 
\item[(5)]
$\nu(\{a,b\})\geq \nu(a)+\nu(b)$.
\end{enumerate}
\end{definition}

In Definition \ref{xxdef1.1}(3), we will introduce the notion of a 
$w$-valuation for any $w\in {\mathbb Z}$. The Poisson valuation described above is essentially a $0$-valuation. Given a Poisson 
valuation $\nu$, we can define a descending filtration
${\mathbb F}^{\nu}:=\{F^{\nu}_i\mid i\in {\mathbb Z}\}$ of $K$ by
\begin{equation}
\label{E0.1.1}\tag{E0.1.1}
F^{\nu}_i:=\{a\in K\mid \nu(a)\geq i\}, \quad {\text{for all $i$}}.
\end{equation}
The {\it associated graded ring} of $\nu$ is defined to be 
\begin{equation}
\label{E0.1.2}\tag{E0.1.2}
\gr_{\nu}(K):=\bigoplus_{i\in {\mathbb Z}} F^{\nu}_i/F^{\nu}_{i+1}.
\end{equation}
The degree 0 part $(\gr_{\nu}(K))_0$ is called the 
{\it residue field} of $\nu$.

 To utilize Poisson valuations more efficiently, we present several associated concepts. We use the Gelfand-Kirillov (GK) dimension instead of the conventional Krull dimension.

\begin{definition}
\label{xxdef0.2}
A Poisson valuation $\nu$ on $K$ is called a {\it faithful valuation} 
if
\begin{enumerate}
\item[(1)]
the image of $\nu$ is ${\mathbb Z}\cup\{\infty\}$,
\item[(2)]
the induced Poisson bracket on the associated graded algebra 
$\gr_{\nu}(K)$ is nonzero, and
\item[(3)]
$\GKdim (\gr_{\nu}(K))={\rm Tr.deg}(K)$.
\end{enumerate}
\end{definition}

Let $Q$ be another Poisson field. We write $K\to_{\nu} Q$ if $Q$ 
is isomorphic to the Poisson fraction field of the 
associated graded Poisson domain $\gr_{\nu}(K)$.  

\begin{definition}
\label{xxdef0.3}
Let $K$ be a Poisson field with nonzero Poisson bracket. 
\begin{enumerate}
\item[(1)] The {\it depth} of $K$ is defined to be
$$\bfd(K):=\sup\,\{n\mid K=K_0\to_{\nu_1} K_1\to_{\nu_2}
K_2 \to \cdots \to_{\nu_{n}} K_n\}$$
where each $\nu_i$ is a faithful valuation and $K_i\not\cong K_j$
for all $0\leq i\neq j\leq n$.
\item[(2)]
The {\it width} of $K$ is defined to be
$$\bfw(K):=\# [\{Q \mid K\to_{\nu} Q\}/\cong]$$
where $\nu$ runs over all faithful valuations on $K$ and 
$-/\cong$ means taking the isomorphism classes. 
\end{enumerate}
\end{definition}

Understanding all arrows $\to_{\nu}$ between Poisson fields can help us visualize the ``quiver" of all Poisson fields together. 

\subsection{Results and Secondary Invariants}
\label{xxsec0.2}
We assume that the base field $\Bbbk$ is of characteristic 0 
for the rest of this subsection. If $P$ is a Poisson domain, 
$Q(P)$ denotes its Poisson fraction field.

\begin{theorem}[Theorem \ref{xxthm6.2}]
\label{xxthm0.4}
Let $T$ be a Poisson torus $\Bbbk[x_1^{\pm 1},\cdots,x_m^{\pm 1}]$
with Poisson bracket determined by $\{x_i,x_j\}=p_{ij} x_i x_j$ 
for some $p_{ij}\in \Bbbk$ for all $1\leq i<j\leq m$. Let $K$ be 
the Poisson field $Q(T)$. Suppose $T$ is Poisson simple. Then 
$\bfd(K)=0$ and $\bfw(K)=1$. In particular, for every faithful 
valuation $\nu$, $Q(\gr_{\nu}(K))\cong K$.
\end{theorem}

A complete list of valuations on $K=Q(T)$ is provided in Theorem \ref{xxthm6.2}(1). We are particularly interested in Poisson fields of transcendence degree two (i.e., GK-dimension two). Below are two well-known examples.

\begin{example}
\label{xxexa0.5}
Let $K_{Weyl}$ denote the Weyl Poisson field that is defined to be 
$\Bbbk(x,y\mid \{x,y\}=1)$. For every $q\in \Bbbk^{\times}$, the $q$-skew Poisson field, denoted by $K_q$, is defined to be $\Bbbk(x,y\mid \{x,y\}=q xy)$.
\end{example}

In addition, our focus is on the Poisson fraction fields of some quotients of three-dimensional Poisson polynomial algebras, which have unimodular Poisson structures defined by homogeneous potentials $\Omega$.

\begin{construction}
\label{xxcon0.6}
Let $\Omega\in \Bbbk[x,y,z]$ be a non-constant homogeneous 
polynomial. We recall some well-known constructions involving 
$\Omega$ as follows.
\begin{enumerate}
\item[(1)]
First, we define a Poisson bracket on the polynomial ring 
$\Bbbk[x,y,z]$ by
\begin{equation}\label{E0.6.1}\tag{E0.6.1}
\{f,g\}_{\Omega}=~\det \begin{pmatrix} 
f_{x} & f_{y} & f_{z}\\
g_{x} & g_{y} & g_{z}\\
\Omega_{x} & \Omega_{y} & \Omega_{z}
\end{pmatrix}\quad 
\text{for all $f,g\in \kk[x, y, z]$}.
\end{equation}

This Poisson polynomial algebra 
$(\Bbbk[x,y,z],\{-,-\}_{\Omega})$ is denoted by $A_{\Omega}$.
This definition is dependent on the set of generators
$(x,y,z)$. However, if we use a new set of generators, 
the new Poisson bracket is a scalar multiple of 
$\{-,-\}_{\Omega}$, see \cite[Definition 3.1]{HTWZ1}. 
\item[(2)]
It is easy to check that $\Omega$ is in the Poisson center of 
$A_{\Omega}$. Hence the Poisson polynomial algebra $A_{\Omega}$ 
has a Poisson factor ring $P_{\Omega-\xi}:
=A_{\Omega}/(\Omega-\xi)$ where $\xi\in \Bbbk$. If $\xi=0$, 
we use $P_{\Omega}$ instead of $P_{\Omega-0}$.
Suppose the {\rm{(}}Adams{\rm{)}} degrees of $x$, $y$, and $z$ 
are 1. If $\Omega$ is homogeneous of degree 3, then 
$P_{\Omega-\xi}\cong P_{\Omega-1}$ when $\xi\neq 0$. So, in 
this case, we can assume that $\xi$ is either 0 or $1$. 
\end{enumerate}
Unless otherwise stated, $A_{\Omega}$, $P_{\Omega}$, 
$P_{\Omega-1}$ and $P_{\Omega-\xi}$ will be the Poisson 
algebras defined as above.
\end{construction}

We will consider two special cases. Case 1: $\Omega: 
=x^3+y^3+z^3+\lambda xyz$ with $\lambda^{3}\neq -3^3$. 
In this case, $\Omega$ has an isolated singularity at the origin, and $A_{\Omega}$ is called an {\it elliptic Poisson algebra}, which has been studied by several authors
\cite{FO, MTU, Pi, Po, TWZ}. Case 2: $\Omega$ is a homogeneous
element of degree $\geq 5$ with an isolated singularity at the 
origin (in this case, we say that $\Omega$ is an i.s. 
potential of degree $\geq 5$). The following result concerns 
the Poisson algebras in Case 1. 

\begin{theorem}
\label{xxthm0.7}
Let $\Omega:=x^3+y^3+z^3+\lambda xyz$ with 
$\lambda^{3}\neq -3^3$.
\begin{enumerate}
\item[(1)]{\rm{(Theorem \ref{xxthm3.8})}}
If $K$ is $Q(P_{\Omega})$, then 
$\bfd(K)=0$ and $\bfw(K)=1$.  
\item[(2)]{\rm{(Theorem \ref{xxthm3.11})}}
If $K$ is $Q(P_{\Omega-1})$, then $\bfd(K)=\bfw(K)=1$.
\end{enumerate}
\end{theorem}

Both $Q(P_{\Omega})$ and $Q(P_{\Omega-1})$ in the above 
theorem play an important role in the partial classification 
in \cite{HTWZ2}. In Case 2, we have a series of results 
[Theorems \ref{xxthm0.8} -- \ref{xxthm0.12}]. (A few results
also hold when $\deg \Omega= 4$.) Let $\Omega$ be in Case 2 for the rest of this subsection. If $P$ is a Poisson 
algebra, the group of Poisson algebra automorphisms of $P$ is 
denoted $\Aut_{Poi}(P)$.

\begin{theorem}[Theorem \ref{xxthm8.1}]
\label{xxthm0.8}
Let $\Omega$ be an i.s. potential of degree $\geq 5$.
\begin{enumerate}
\item[(1)]
$$\Aut_{Poi}(Q(P_{\Omega}))=\Aut_{Poi}(P_{\Omega}).$$
\item[(2)]
If $\xi\neq 0$, then
$$\Aut_{Poi}(Q(P_{\Omega-\xi}))=\Aut_{Poi}(P_{\Omega-\xi}).$$
\end{enumerate}
\end{theorem}

Our primary objective is to provide effective and comprehensive solutions to significant questions like automorphism and embedding problems. To accomplish this, we have introduced invariants that stem from Poisson valuations. 
\begin{enumerate}
\item[(1)]
$\alpha$-type invariant of a Poisson field $K$ is an invariant 
that is defined using the Poisson ($w$-)valuations directly on 
$K$, e.g., the number of faithful valuations of $K$ listed in 
Table 1.
\item[(2)]
$\beta$-type invariant is any invariant related to the arrow
$\longrightarrow_{\nu}$. So the depth and width of $K$ in 
Definition \ref{xxdef0.2} are $\beta$-type invariants.
\item[(3)]
$\gamma$-type invariant is any invariant related to the filtrations
${\mathbb F}^\nu$ associated to valuations $\nu$.
\item[(4)]
In \cite{HTWZ3}, we also introduce $\delta$-type invariant 
that is deduced from the associated graded ring $\gr_{\nu}(K)$
defined in \eqref{E0.1.2}.
\end{enumerate}

Secondary invariants are very useful in different Poisson algebra projects. Below is a summary of some Poisson fields and their faithful 
valuations studied in this paper. 

\begin{table}[H]
\centering
\begin{tabular}{|c|c|c|c|}
\hline
\multirow{2}{*} {Poisson fields of tr.deg 2}  & \multirow{2}{*}{fields}   & Poisson 
  &   $\#\{\text{faithful}$\\
  &   &brackets
  & $\text{Poisson valuations}\}$\\
\hline
Weyl Poisson field $K_{Weyl}$  &   \multirow{2}{*}{ $\kk(x,y)$} 
  &  \multirow{2}{*}{$\{x,y\}=1$} & $\infty$ \\
 over uncountable field $\kk$&&&uncountable\\
 
\hline
\multirow{2}{*}{$q$-Skew Poisson field $K_{q}$ ($q\neq 0$) }  &  \multirow{2}{*}{$\kk(x,y)$} 
& \multirow{2}{*}{$\{x,y\}=qxy$} &  $\infty$   \\
&&&countable\\ \hline
Graded elliptic type   & \multirow{2}{*}{$Q(P_\Omega)$}
&  \multirow{2}{*}{$\{-,-\}_\Omega$} & \multirow{2}{*}{$2$}\\
$\Omega=x^3+y^3+z^3+\lambda xyz$ ($\lambda^3\neq -3^3$) &&&\\
\hline
Elliptic type & \multirow{2}{*}{$Q(P_{\Omega-1})$} & \multirow{2}{*}{$\{-,-\}_\Omega$}  
  &  \multirow{2}{*}{$1$}\\
 $\Omega=x^3+y^3+z^3+\lambda xyz$ ($\lambda^3\neq -3^3$) &&&\\
\hline
Higher genus type  & $Q(P_{\Omega})$ &  
\multirow{2}{*}{$\{-,-\}_\Omega$}   & \multirow{2}{*}{0}\\
$\Omega$ is an i.s. potential of $\deg\ge 4$ &$Q(P_{\Omega-1})$ && \\
\hline
\end{tabular}
\vskip 1em
\caption{Some Poisson fields of tr.deg 2 and their faithful valuations}
\label{tab:my_label}
\end{table}

\subsection{More Applications}
\label{xxsec0.3}

The automorphism groups of other Poisson fields (and Nambu-Poisson fields) are computed in \cite{HTWZ1, HTWZ2, HTWZ3} 
using Poisson valuations. Motivated by the Dixmier conjecture and Poisson conjecture, see Section \ref{xxsec8.2}, we consider the Dixmier property [Definition \ref{xxdef8.5}].

\begin{theorem}[Theorem \ref{xxthm8.6}]
\label{xxthm0.9}
Let $\Omega$ be an i.s. potential of degree $\geq 5$.
Then, for any $\xi\in \Bbbk$, every injective endomorphism of 
$P_{\Omega-\xi}$ {\rm{(}}resp. $Q(P_{\Omega-\xi})${\rm{)}} is an automorphism.
\end{theorem}

Other properties, such as the uniqueness of grading/filtration,
are motivated by some work on noncommutative algebra, for 
example, \cite[Corollary 0.3]{BZ}. We refer to such a property 
as the rigidity of grading/filtration. The next two results
are in this direction. Some undefined terms can be found in 
Sections \ref{xxsec1} and \ref{xxsec2}. We say that a 
${\mathbb Z}$-graded algebra $A=\bigoplus_{i\in {\mathbb Z}} A_i$ 
is {\it connected graded} if $A_i=0$ for all $i>0$ and 
$A_0=\Bbbk$. So, $A$ lives in nonpositive degrees.  Note that 
$P_{\Omega}$ is connected graded if we set the degree of 
$x$, $y$, and $z$ to be $-1$. This is different from the 
traditional definition of connected graded algebra. The 
reason for using this non-traditional definition is to match 
up with our definition of descending filtration. A Poisson
$n$-graded algebra is defined in Definition \ref{xxdef2.1}.

\begin{theorem}[Theorem \ref{xxthm8.9}]
\label{xxthm0.10}
Let $\Omega$ be an i.s. potential of degree $\geq 4$. Then 
$P_{\Omega}$ has a unique connected grading such that it is 
Poisson $(\deg \Omega -3)$-graded. 
\end{theorem}

\begin{theorem}[Theorem \ref{xxthm8.10}]
\label{xxthm0.11}
Let $\Omega$ be an i.s. potential of degree $\geq 4$.
If $\xi\neq 0$, then $P_{\Omega-\xi}$ has a unique 
filtration ${\mathbb F}$ such that the associated graded 
ring $\gr_{\mathbb F} (P_{\Omega-\xi})$ is a 
connected graded Poisson $(\deg \Omega-3)$-graded domain. 
\end{theorem}

The results mentioned above aid in calculating the associated automorphism groups. An instance of this is 

\begin{theorem}[Theorem \ref{xxthm8.2}]
\label{xxthm0.12}
Let $\Omega$ be an i.s. potential of degree $\geq 5$.
Then $\Aut_{Poi}(Q(A_{\Omega})) =\Aut_{Poi}(A_{\Omega})=
\Aut_{Poi}(P_{\Omega})$ and every Poisson automorphism of 
$A_{\Omega}$ is graded. Consequently, 
$\Aut_{Poi}(Q(A_{\Omega}))$ is a finite subgroup of 
$GL_3(\Bbbk)$ of order bounded above by 
$42(\deg \Omega)(\deg \Omega-3)^2$.
\end{theorem}

When $d\geq 5$ and $\Omega=x^d+y^d+z^d=0$ is the Fermat curve, we can explicitly compute the Poisson automorphism groups of $P_{\Omega-\xi}$ and $A_\Omega$ according to Proposition \ref{xxpro8.4}(1,2,3). There is a close connection between the Poisson automorphism group of $A_\Omega$ and the automorphism group of the projective curve $X=\Proj (A/(\Omega))$. For instance, the order bound of $\Aut_{Poi}(Q(A_{\Omega}))$ follows from applying Hurwitz's automorphism theorem to the smooth curve $X$. We're curious if Theorems \ref{xxthm0.8}, \ref{xxthm0.9}, and \ref{xxthm0.12} hold when $\deg \Omega=4$.

\begin{remark}
\label{xxrem0.13} 
\begin{enumerate}
\item[(1)]
There is a Poisson field of transcendence degree 2 that
does not admit any nontrivial Poisson valuation, see
Lemma \ref{xxlem4.8}(3).
\item[(2)]
We will introduce a weighted version of Poisson valuations 
so that the Poisson field in part (1) admits a weighted version of a nontrivial Poisson valuation, see 
Lemma \ref{xxlem5.3}(2). 
\item[(3)]
Weighted valuations are used in the proofs of Theorems
\ref{xxthm0.8}-\ref{xxthm0.12}. 
\end{enumerate}
\end{remark}

Before we close this introduction, we want to emphasize that Poisson valuations are useful in the following topics.

\begin{enumerate}
\item[(1)]
The automorphism problem. In \cite{MTU}, 
Makar-Limanov, Turusbekova, and Umirbaev computed the Poisson automorphism groups of the elliptic Poisson algebras. The automorphism groups of other Poisson algebras are computed 
in \cite{HTWZ1, HTWZ2, HTWZ3}. Using Poisson valuations, one can 
also compute the automorphism group of a family of Poisson fields. 
Theorems \ref{xxthm0.8} and \ref{xxthm0.12} are the initial outcomes of this type. Note that the valuation method has been used to solve
the automorphism problem \cite{LY2, Ya}.
\item[(2)]
The isomorphism problem (see Section \ref{xxsec6}), the embedding problem (see Section \ref{xxsec6}), the rigidity of grading (see Theorem \ref{xxthm0.10}), the rigidity of filtration (see Theorem \ref{xxthm0.11}), and the Dixmier problem (see Theorem \ref{xxthm0.9}).
\item[(3)]
Classification of Poisson fields of transcendence degree two 
is an important project related to Artin's conjecture 
in noncommutative algebraic geometry.  
Artin's conjecture concerning division algebras of transcendence degree two  \cite{Ar} is 
(arguably) the most important open problem in noncommutative 
algebraic projective geometry today. The conjecture holds significant sway over the progress of the research field, attracting the attention and dedication of numerous renowned experts worldwide. 
We initiated this project to establish a Poisson version of Artin's conjecture. Further details about the said version can be found in \cite{HTWZ2}.
\item[(4)]
Solutions (and non-existence of solutions) to 
partial differential equations in the field $\Bbbk(x_1,\cdots,x_n)$ 
\cite{GZ, HTWZ3}.
\item[(5)]
Study of $n$-Lie Poisson (or Nambu-Poisson) algebras \cite{HTWZ3}.
\end{enumerate}

\section{Preliminaries}
\label{xxsec1}

This section reviews some introductory material related 
to valuations and Poisson valuations. We refer to the book \cite{LPV} for basic definitions 
of Poisson algebras and \cite[Chaprter VI]{Bo} (and 
\cite{Va}) for basic properties related to valuations.

\begin{definition}
\label{xxdef1.1}
Let $K$ be a commutative algebra (or a field) over $\Bbbk$. For 
parts {\rm{(3)-(5)}}, let $K$ be a Poisson algebra (or a Poisson 
field).
\begin{enumerate}
\item[(1)]
A {\it discrete valuation} or simply {\it valuation} on $K$ 
is a map
$$\nu: K \to {\mathbb Z}\cup\{\infty\}$$
which satisfies the following properties: for all $a,b\in K$, 
\begin{enumerate}
\item[(a)]
$\nu(a)=\infty$ if and only if $a=0$,
\item[(b)]
$\nu(a)=0$ for all $a\in \Bbbk^{\times}:=\Bbbk\setminus \{0\}$,
\item[(c)]
$\nu(ab)=\nu(a)+\nu(b)$ (assuming $n+\infty=\infty$ when $n\in
{\mathbb Z}\cup\{\infty\}$),
\item[(d)]
$\nu(a+b)\geq \min\{\nu(a),\nu(b)\}$ (with equality if $\nu(a)
\neq \nu(b)$),
\end{enumerate}
\item[(2)]
A valuation $\nu$ is called {\it trivial} if $\nu(a)=0$ for 
all $a\in K\setminus\{0\}$. Otherwise, it is called {\it nontrivial}.
\item[(3)]
Let $w$ be a given integer. A valuation $\nu$ on $K$ is called a
{\it $w$-valuation} if
\begin{enumerate}
\item[(e)]
$\nu(\{a,b\})\geq \nu(a)+\nu(b)-w$ for all $a,b\in K$.
\end{enumerate}
Any $0$-valuation is just a {\it Poisson valuation} given in
Definition \ref{xxdef0.1}. When $\nu$ is a $w$-valuation, we 
say $w$ is the {\it weight} of $\nu$. Let ${\mathcal V}_w(K)$ 
be the set of nontrivial $w$-valuations on $K$.
\item[(4)]
A $w$-valuation $\nu$ on $K$ is called {\it classical} if 
$\nu(\{a,b\})> \nu(a)+\nu(b)-w$ for all $a,b\in K\setminus \{0\}$. 
Otherwise, it is called {\it nonclassical}.
\item[(5)]
Let $\nu_1$ (resp. $\nu_2$) be a $w_1$-valuation
(resp. $w_2$-valuation) on $K$. We say $\nu_1$ and $\nu_2$ are
{\it equivalent} if there are positive integers $a,b$ such 
that $a\nu_1=b\nu_2$. In this case, we write $\nu_1\sim \nu_2$. 
Otherwise we say $\nu_1$ and $\nu_2$ are {\it nonequivalent}.
\end{enumerate}
\end{definition}

In this paper, we will only study discrete valuations. Thus, the term ``discrete" will be omitted throughout this paper. Note that our definition of equivalent valuations differs from the one in \cite[Proposition 1.5]{Va}. By convention, a $w$-valuation is defined on a Poisson algebra, while a valuation is defined on a commutative algebra without considering the Poisson structure. Although $0$-valuations play an important role in this paper, $w$-valuations with $w\neq 0$ also provide useful information. For instance, there 
are Poisson fields such that (a) there is no nontrivial 
$0$-valuation [Lemma \ref{xxlem4.8}(3)] and (b) there are 
nontrivial Poisson $w$-valuations for some $w>0$ [Lemmas \ref{xxlem3.10} and \ref{xxlem5.3}]. Similar to Definition 
\ref{xxdef1.1}(3), we define 
\begin{equation}
\label{E1.1.1}\tag{E1.1.1}
{\mathcal V}_{cw}(K):={\text{the set of nontrivial classical 
$w$-valuations on $K$}}
\end{equation}
and
\begin{equation}
\label{E1.1.2}\tag{E1.1.2}
{\mathcal V}_{ncw}(K):={\text{the set of nontrivial 
nonclassical $w$-valuations on $K$}}.
\end{equation}
It is clear that ${\mathcal V}_w(K)={\mathcal V}_{cw}(K)
\sqcup {\mathcal V}_{ncw}(K)$ for all $w$. Further, we define
\begin{enumerate}
\item[]
\begin{equation}
\label{E1.1.3}\tag{E1.1.3}
w(K):=\inf\{ w\mid {\mathcal V}_w(K)\neq \emptyset\}.
\end{equation}
\end{enumerate}
By Remark \ref{xxrem4.4}(2), $w(K)$ is either $-\infty$, 
$0$ or $1$. Similarly, one can define $cw(K)$ and $ncw(K)$. 

When $K$ is a field with valuation $\nu$, the {\it valuation 
ring} $D_{\nu}(K):=\{a\in K\mid \nu(a)\geq 0\}$ is local with 
maximal ideal $\fm_{\nu}(K):=\{a\in K\mid \nu(a)\geq 1\}$. The 
corresponding {\it residue field} is $R_{\nu}(K):=D_{\nu}(K)/
\fm_{\nu}(K)$. If $\nu$ is nontrivial, then the transcendence 
degree of $R_{\nu}(K)$ is strictly less than that of $K$. When 
$\nu$ is a $0$-valuation, then $R_{\nu}(K)$ is a Poisson field. 
We refer to \cite{KL} for basic definitions and properties 
related to Gelfand–Kirillov dimension (GKdim or GK-dimension). 
It's important to note that for any affine commutative algebra, its GK-dimension always equals its Krull dimension.
\begin{definition}
\label{xxdef1.2}
Let $\nu$ be a $w$-valuation on a Poisson field $K$. We say 
$\nu$ is {\it degenerate} if $\GKdim R_{\nu}(K) \leq \GKdim 
K-2$. Otherwise, we say $\nu$ is {\it nondegenerate}.
\end{definition}

The following lemma may be well-known.

\begin{lemma}
\label{xxlem1.3}
Suppose $\Bbbk$ is algebraically closed. Let $A$ be a
${\mathbb Z}$-graded domain. In parts {\rm{(1)}} and {\rm{(2)}} we suppose that $A_i\neq 0$ and $A_{-j}\neq 0$ for some $i,j>0$.
\begin{enumerate}
\item[(1)] 
Suppose $A_{k}$ is nonzero and finite-dimensional for some $k$. 
Then $A\cong \Bbbk[s^{\pm 1}]$ for some nonzero homogeneous 
element $s\in A$ of positive degree.
\item[(2)]
Suppose $\GKdim A=1$. Then $A\cong \Bbbk[s^{\pm 1}]$ for some 
nonzero homogeneous element $s\in A$ of positive degree.
\item[(3)]
Suppose $\GKdim A=1$. If $A$ is ${\mathbb N}$-graded with 
$A_i\neq 0$ for some $i>0$, then $A\subseteq \Bbbk[s]$ 
for some nonzero homogeneous element $s\in A$ of positive degree.
\end{enumerate}
\end{lemma}

\begin{proof} (1) Let $0\neq z\in A_{k}$. Since $A$ is a domain,
$A_0 z\subseteq A_k$ implies that $$\dim_{\Bbbk} A_0
=\dim_{\Bbbk} A_0 z\leq \dim_{\Bbbk} A_k<\infty.$$ 
Then $A_0$ is an algebraic field extension of the base field
$\Bbbk$. Since $\Bbbk$ is algebraically closed, $A_0=\Bbbk$.

Let $0\neq x\in A_i$ and $0\neq y\in A_{-j}$. Then
$0\neq x^j y^i\in A_0$. So $x^j y^i\in \Bbbk^{\times}$
is invertible, which implies that $x$ (and $y$) is 
invertible. Consequently, every nonzero homogeneous 
element is invertible. Thus, $A$ is a graded division
algebra. Since $A_0=\Bbbk$, $A$ must be 
$\Bbbk[s^{\pm 1}]$ for some nonzero homogeneous 
element $s\in A$ of positive degree.

(2,3) The proofs are similar to the proof of part (1) and 
omitted.
\end{proof}

\begin{lemma}
\label{xxlem1.4}
Let $K\subseteq Q$ be two Poisson fields with the same finite 
GK-dimension. Let $w$ be an integer.
\begin{enumerate}
\item[(1)]
There is a natural map induced by restriction
$${\mathcal V}_w(Q)\to {\mathcal V}_w(K).$$
As a consequence, if ${\mathcal V}_w(Q)\neq \emptyset$,
then ${\mathcal V}_w(K)\neq \emptyset$. This also implies
$w(K)\leq w(Q)$.
\item[(2)]
$${\mathcal V}_{w-1}(K)\subseteq {\mathcal V}_{cw}(K)
\subseteq {\mathcal V}_{w}(K).$$
As a consequence, if ${\mathcal V}_{w}(K)={\mathcal V}_{ncw}(K)$,
then ${\mathcal V}_{w-1}(K)=\emptyset$.
\item[(3)]
Let $d$ be a positive integer. Then the assignment
$\nu\to d \nu$ defines an injective map ${\mathcal V}_{w}(K)
\to {\mathcal V}_{dw}(K)$.
As a consequence, if ${\mathcal V}_w(K)\neq \emptyset$
for some $w<0$, then $w(K)=-\infty$.
\item[(4)]
Part {\rm{(3)}} holds for ${\mathcal V}_{cw}(K)$ and 
${\mathcal V}_{ncw}(K)$.
\end{enumerate}
\end{lemma}

\begin{proof} 
(1) Let $\nu\in {\mathcal V}_w(Q)$ and $\phi:=\nu\mid_{K}$. It 
is clear that $\phi$ is a $w$-valuation on $K$. It remains to 
show that $\phi$ is nontrivial. Suppose to the contrary that 
$\phi$ is trivial. So $\phi(f)=0$ for all $0\neq f\in K$. We 
claim that if $0\neq q\in Q$, then $\nu(q)=0$. If not, we can 
assume $\nu(q)<0$ (after replacing $q$ by $q^{-1}$ if 
necessary). Since $K$ and $Q$ have the same finite GK-dimension,
there are $a_0, \cdots, a_n\in K$ (for some $n\geq 0$) with 
$a_0\neq 0$ such that $q^{n+1} =a_0+a_1 q+\cdots+ a_{n} q^{n}$. 
Then $(n+1)\nu(q)\geq \min_{i=0}^{n}\{ \nu(a_i q^{i})\}
\geq n\nu(q)$ which contradicts the assumption 
$\nu(q)<0$. Therefore, we have proved the claim.
Since $\nu$ is nontrivial, we have that $\phi$ is nontrivial.
This is the central assertion. The consequences are clear.

(2,3,4) The proofs are easy and omitted.
%
%
%
\end{proof}

\begin{remark}
\label{xxrem1.5}
Let $K$ be a Poisson field.
\begin{enumerate}
\item[(1)]
We can use Lemma \ref{xxlem1.4}(2,3) to generate classical Poisson valuations by upgrading the weight of one $w$-valuation. Therefore, it is reasonable to remove these classical valuations and focus solely on nonclassical ones up to equivalence. To achieve this, we need to understand the ``valuation lattice" defined as 
$${\mathcal L}(K):=\frac{\bigcup_w {\mathcal V}_{ncw}(K)}{\sim}.$$
An $\alpha$-type invariant of $K$ is any invariant defined using the ($w$-)valuations on $K$. Hence, ${\mathcal L}(K)$ is an $\alpha$-type invariant of $K$.  
\item[(2)]
We can define other $\alpha$-type invariants. For example, let 
\begin{equation}
\label{E1.5.1}\tag{E1.5.1}
\alpha_{w}(K):=\# \left\{\frac{{\mathcal V}_{w}(K)}{\sim}\right\}, 
\end{equation}
which denotes the number (or cardinality) of all nontrivial 
$w$-valuations $\nu$ on $K$ where there is no $w'$-valuation 
$\nu'$ on $K$ such that $\nu'=a \nu$ for some positive rational 
number $a<1$. Computing $\alpha_{w}(K)$ for all $w$ would be interesting. One could also define $\alpha_{cw}(K)$ (resp. $\alpha_{ncw}(K)$) by replacing ${\mathcal V}_{w}(K)$ with ${\mathcal V}_{cw}(K)$-- see \eqref{E1.1.1} (resp. 
${\mathcal V}_{ncw}(K)$ -- see \eqref{E1.1.2}) in \eqref{E1.5.1}. 
Note that $w(K)$ in \eqref{E1.1.3} is also an $\alpha$-type invariant of $K$.
\item[(3)]
The depth and width defined in Definition \ref{xxdef0.3} (or any invariants related to the arrow $\to_\nu$) are referred to as $\beta$-type invariants. Additional invariants will be introduced in later sections. 
\end{enumerate}
\end{remark}

The following lemma will be used in the sequel. 

\begin{lemma}
\label{xxlem1.6}
Suppose ${\text{char}}\; \Bbbk=0$. Let $K$ be a Poisson field 
and $x,y\in K$.
\begin{enumerate}
\item[(1)]
If $x$ and $y$ are algebraically dependent, then $\{x,y\}=0$.
\item[(2)]
Suppose $K$ has GK-dimension 2 with a nontrivial Poisson bracket. 
If $\{x,y\}=0$, then $x$ and $y$ are algebraically dependent.
\end{enumerate}
\end{lemma}

\begin{proof}
(1) For every $x\in K$, let $C_x(K):=\{f\in K\mid \{x,f\}=0\}$. 
Since ${\text{char}}\; \Bbbk=0$, $C_x(K)$ is a subfield of $K$ 
that is integrally closed in $K$ (namely, $C_x(K)$ is equal to 
its integral closure in $K$). Without loss of generality, we 
may assume that $x\not\in \Bbbk$. Since $x\in C_x(K)$ and $y$ is integral over $\Bbbk(x)$, $y$ is 
integral over $C_x(K)$. Hence $y\in C_x(K)$. This means that 
$\{x,y\}=0$.

(2) Suppose to the contrary that $x$ and $y$ are algebraically
independent. Since $\{x,y\}=0$, both $x$ and $y$ are in $C_x(K)$.
Then $C_x(K)$ has a transcendence degree of at least 2. Note that $K$ 
has transcendence degree 2. Since $C_x(K)$ is integrally closed 
in $K$, we have $C_x(K)=K$. This means that $x$ is in the center 
of $K$. Similarly, $y$ is in the center of $K$. For every $f\in K$, 
then $x,y\in C_f(K)$. Since $C_f(K)$ is integrally closed in $K$, 
we obtain that $C_f(K)=K$. Thus, $K$ has a trivial Poisson bracket,
yielding a contradiction.
\end{proof}

An ideal $I$ of a Poisson algebra $P$ is called a {\it Poisson ideal} 
if $\{I, P\}\subseteq I$. A Poisson ideal $I$ is called {\it Poisson 
prime} if, for any two Poisson ideals $J_1, J_2$, $J_1 J_2\subseteq I$ 
implies that $J_1\subseteq I$ or $J_2\subseteq I$. 

\begin{lemma} \cite[Lemma 1.1(c,d)]{Go1}
\label{xxlem1.7}
If $P$ is a noetherian Poisson algebra, then an ideal
$I$ of $P$ is Poisson prime if and only if it is both
Poisson and prime. Moreover, every prime ideal minimal 
over a Poisson ideal of $P$ is Poisson prime. 
\end{lemma}

\section{Grading, degree, filtrations, and valuations}
\label{xxsec2}
In this section, we establish some foundational concepts to use valuations effectively. First, we recall the notion of a Poisson $w$-graded algebra and a slightly weaker version of a valuation called a {\it filtration} of an algebra. Throughout, let $w$ be a given integer unless otherwise stated.

\begin{definition}
\label{xxdef2.1}
A Poisson algebra $P$ is called {\it Poisson $w$-graded} if
\begin{enumerate}
\item[(1)]
$P=\oplus_{i\in {\mathbb Z}} P_i$ is a ${\mathbb Z}$-graded 
commutative algebra, and  
\item[(2)]
$\{P_i, P_j\}\subseteq P_{i+j-w}$ for all $i,j\in {\mathbb Z}$.
\end{enumerate}
\end{definition}

By Construction \ref{xxcon0.6}, the Poisson algebras $A_{\Omega}$ 
and $P_{\Omega}$ are $w$-graded where $w=-(\deg \Omega-3)$, if we
use the original grading (i.e., the standard Adams grading with $|x|=|y|=|z|=1$) of $A_{\Omega}$ 
and $P_{\Omega}$.

\begin{definition}
\label{xxdef2.2}
Let $A$ be an algebra. Let
${\mathbb F}:=\{F_{i}\mid i \in {\mathbb Z}\}$
be a set of $\Bbbk$-subspaces of $A$. 
\begin{enumerate}
\item[(1)]
We say ${\mathbb F}$ is a {\it filtration} of $A$ if it satisfies
\begin{enumerate}
\item[(a)]
 $F_{i}\supseteq F_{i+1}$ for all $i\in
{\mathbb Z}$ and $1\in F_0\setminus F_1$,
\item[(b)]
$F_i F_j\subseteq F_{i+j}$ for all $i,j\in {\mathbb Z}$,
\item[(c)]
$\bigcap_{i\in{\mathbb Z}} F_i=\{0\}$,
\item[(d)]
$\bigcup_{i\in {\mathbb Z}} F_i= A$.
\end{enumerate}
\item[(2)]
Suppose that $A$ is a Poisson algebra and ${\mathbb F}$ is a filtration of $A$. If further
\begin{enumerate}
\item[(e)]
$\{F_i,F_j\}\subseteq F_{i+j-w}$ for all $i,j\in {\mathbb Z}$,
\end{enumerate}
then ${\mathbb F}$ is called a {\it $w$-filtration} of
$A$.
\end{enumerate}
\end{definition}

For a moment, let's set aside the Poisson structure. The {\it associated graded ring} of a filtration ${\mathbb F}$ of
$A$ is defined to be
$$\gr_{\mathbb F} A:= \bigoplus_{i\in{\mathbb Z}} F_i/F_{i+1},$$
which is a ${\mathbb Z}$-graded algebra. By \cite[Lemma 6.5]{KL},
$\GKdim \gr_{\mathbb F} A\leq \GKdim A$. For any element 
$a\in F_i$, let ${\overline{a}}$ denote the element $a+F_{i+1}$ 
in the $i$th degree component $(\gr_{\mathbb F} A)_i:=F_i/F_{i+1}$. 
One canonical example is the following. Let $A$ be an algebra 
generated by a subspace $V$ containing $1$. Let
\begin{equation}
\notag 
F_i=\begin{cases} V^{-i} & i\leq 0, \\
0& i>0.\end{cases}
\end{equation}
Then ${\mathbb F}:=\{F_i\}$ is a filtration of $A$ satisfying
conditions (a,b,c,d) in Definition \ref{xxdef2.2}(1). A slightly 
more general example is given in \eqref{E2.8.3}. Next, we reintroduce the Poisson structure into our discussion.

\begin{lemma}
\label{xxlem2.3}
Suppose ${\mathbb F}$ is a $w$-filtration of a 
Poisson algebra $A$. Then $\gr_{\mathbb F} A$ is a 
Poisson $w$-graded algebra.
\end{lemma}

\begin{proof} 
It is well-known that $\gr_{\mathbb F} A$ is a graded algebra.
The addition and multiplication are defined as follows:
${\overline{f}}+{\overline{g}}=\overline{f+g}$ when
${\overline{f}},{\overline{g}}\in (\gr_{\mathbb F} A)_i$
and ${\overline{f}}{\overline{g}}=\overline{fg}$ when
${\overline{f}}\in (\gr_{\mathbb F} A)_i$ and
${\overline{g}}\in (\gr_{\mathbb F} A)_j$.

Next, we define the $w$-graded Poisson structure on $\gr_{\mathbb F} A$.
Let ${\overline{f}}$ and ${\overline{g}}$ be in
$(\gr_{\mathbb F} A)_i$ and $(\gr_{\mathbb F} A)_j$ respectively
where $f\in F_{i}$ and $g\in F_{j}$. By definition,
$\{f,g\}\in F_{i+j-w}$. Let $\overline{\{f,g\}}$ be its
class in $(\gr_{\mathbb F} A)_{i+j-w}$. We define the Poisson
bracket on $\gr_{\mathbb F} A$ by 
$$\{{\overline{f}},{\overline{g}}\}
:= \overline{\{f,g\}}\in (\gr_{\mathbb F} A)_{i+j-w}$$
for all ${\overline{f}}\in (\gr_{\mathbb F} A)_i$ and 
${\overline{g}}\in (\gr_{\mathbb F} A)_j$. It is easy to see 
that $\{-,-\}$ is well-defined (or is independent of the choices 
of preimages $f$ and $g$), bilinear, and antisymmetric.

We claim that $\{-,-\}$ is a Poisson bracket. Let $a\in
F_i, b\in F_j, c\in F_k$. Then in
$(\gr_{\mathbb F} A)_{i+j+k-w}$, we have
$$\begin{aligned}
\{\overline{a}, \overline{b}\overline{c}\}
&=\{\overline{a}, \overline{bc}\}=\overline{\{a,bc\}}\\
&=\overline{\{a,b\}c+\{a,c\}b}=\overline{\{a,b\}c}+\overline{\{a,c\}b}\\
&=\overline{\{a,b\}}\overline{c}+\overline{\{a,c\}}\overline{b}
=\{\overline{a},\overline{b}\}\overline{c}
  +\{\overline{a},\overline{c}\}\overline{b}.
\end{aligned}
$$
So $\{\overline{a}, -\}$ is a derivation. Similarly,
the Jacobi identity holds for $\overline{a}, \overline{b},
\overline{c}$. Therefore, $\gr_{\mathbb F} A$ is a Poisson
$w$-graded algebra.
\end{proof}

Given a filtration ${\mathbb F}$, we define the notion of a 
{\it degree} function, denoted by $\deg$, on elements in $A$ by
\begin{equation}
\label{E2.3.1}\tag{E2.3.1}
\deg a:=i, {\text{ if $a\in F_i\setminus F_{i+1}$}} \quad
{\text{and}} \quad \deg(0)=+\infty.
\end{equation}

\begin{lemma}
\label{xxlem2.4}
Let ${\mathbb F}$ be a filtration of an algebra $A$
such that $\gr_{\mathbb F} A$ is a domain.
\begin{enumerate}
\item[(1)]
Then $\deg$ satisfies the following conditions for $a,b\in A$,
\begin{enumerate}
\item[(a)]
$\deg(a)\in {\mathbb Z}$ for any
nonzero element $a\in A$, and $\deg a=\infty$ if and only if
$a=0$,
\item[(b)]
$\deg(c)=0$ for all $c\in \Bbbk^{\times}$,
\item[(c)]
$\deg(ab)=\deg(a)+\deg(b)$,
\item[(d)]
$\deg(a+b)\geq \min\{\deg(a),\deg(b)\}$, with
equality if $\deg(a)\neq \deg(b)$.
\end{enumerate}
Namely, $\deg$ is a valuation on $A$.
\item[(2)]
Suppose $A$ is a Poisson algebra. If ${\mathbb F}$ is a $w$-filtration of $A$,
then $\deg$ is a $w$-valuation in the sense of
Definition {\rm{\ref{xxdef1.1}(3)}}.
\item[(3)]
$F_0(A)$ is integrally closed in $A$.
\end{enumerate}
\end{lemma}

\begin{proof} (1) This follows from a routine verification.

(2) Conditions (a,b,c,d) in Definition \ref{xxdef1.1}(1) 
were proved in part (1). Condition (e) in Definition 
\ref{xxdef1.1}(3) follows from Definition \ref{xxdef2.2}(2e). 
The assertion follows.

(3) Let $f\in A$ be integral over $F_0(A)$. This means that 
there are $a_0,\cdots,a_{n-1}\in F_0(A)$ with $n>0$, and 
$\sum_{i=0}^n a_i f^i=0$ where $a_n=1$. We claim that 
$f\in F_0(A)$. Denote $\nu:=\deg$. It suffices to show that 
$\nu(f)\geq 0$. Suppose to the contrary that $\nu(f)<0$. It 
follows from the equation $\sum_{i=0}^n a_i f^i=0$ (or $f^n=-
\sum_{i=0}^{n-1}a_i f^i$) that
$$n\nu(f)=\nu(f^n) =\nu(-\sum_{i=0}^{n-1}a_i f^i)
\geq \min_{i=0}^{n-1}\{i \nu(f)+\nu(a_i)\}\geq (n-1) \nu(f),$$
which implies that $\nu(f)\geq 0$, a contradiction.
Therefore, $\nu(f)\geq 0$ is required.
\end{proof}

Conversely, if we are given a valuation denoted by $\nu$ (or 
more generally, a degree function satisfying (1a), (1b), (1c) and (1d) of Lemma \ref{xxlem2.4})
$$\nu: A\to {\mathbb Z}\cup\{\infty\},$$ 
 then we can 
define a filtration (associated with $\nu$) 
${\mathbb F}^{\nu}:=\{F_i^{\nu}\mid i\in {\mathbb Z}\}$
of $A$ by
\begin{equation}
\label{E2.4.1}\tag{E2.4.1}
F_i^{\nu}:=\{ a\in A\mid \nu(a)\geq i\}.
\end{equation}
If no confusion occurs, we will delete $^{\nu}$ from 
${\mathbb F}^{\nu}$ and $F_i^{\nu}$.

\begin{lemma}
\label{xxlem2.5}
Let $A$ be a Poisson algebra with a valuation $\nu$. 
Then ${\mathbb F}:=\{F_i\}$ defined as in \eqref{E2.4.1} is a
filtration of $A$ such that $\gr_{\mathbb F} A$ is a domain.
If $\nu$ is a $w$-valuation, then $\gr_{\mathbb F} A$
is a Poisson $w$-graded domain.
\end{lemma}

\begin{proof}
This is another routine verification. 
\end{proof}

A filtration ${\mathbb F}$ is called {\it good} if 
$\gr_{\mathbb F} A$ is a domain. One can check immediately 
that there is a one-to-one correspondence between the set 
of valuations on $A$ and the set of good filtrations 
on $A$. The following result applies in the Poisson setting.

\begin{lemma}
\label{xxlem2.6}
For a Poisson algebra $A$, there is a one-to-one
correspondence between the set of good $w$-filtrations 
of $A$ and the set of $w$-valuations on $A$.
\end{lemma}

\begin{proof}
Given a good $w$-filtration, by Lemma \ref{xxlem2.4}(2),
$\deg(=:\nu)$ defined in \eqref{E2.3.1} is a $w$-valuation on $A$.
Conversely, given a $w$-valuation $\nu$, the filtration 
${\mathbb F}$ defined in \eqref{E2.4.1} is a good filtration 
by Lemma \ref{xxlem2.5}. It is easy to see that ${\mathbb F}$ 
is a $w$-filtration. Moreover, the correspondence 
between $\deg$ and ${\mathbb F}$ is one-to-one by construction.
\end{proof}

Suppose $\gr_{\mathbb F} A$ is a domain for some filtration 
$\mathbb{F}$ of $A$, and so is $A$. The corresponding degree 
function $\deg(=:\nu)$ (see \eqref{E2.3.1}) satisfies all 
conditions in Lemma \ref{xxlem2.4}(1). Let $Q$ be the fraction field of $A$. Now, we define a {\it degree} on $Q$ by
\begin{equation}
\label{E2.6.1}\tag{E2.6.1}
\deg_Q(ab^{-1}):=\deg a-\deg b
\end{equation}
for $a\in A$ and $b\in A\setminus\{0\}$. It is easy to check 
that $\deg_Q$ is well-defined (or only depends on the 
equivalence class of $ab^{-1}$).


If $A$ is a graded domain, the graded fraction ring $Q_{gr}(A)$ of
$A$ is defined to be $AS^{-1}$ where $S$ is the set of nonzero homogeneous elements.

\begin{lemma}
\label{xxlem2.7}
Retain the above notations. Suppose $A$ is a Poisson domain and ${\mathbb F}$ is a good 
$w$-filtration of $A$.
\begin{enumerate}
\item[(1)]
$\deg_Q$ is a $w$-valuation on $Q$ and its associated filtration ${\mathbb F}(Q)$ is a good $w$-filtration of
$Q$.
\item[(2)]
The associated graded ring $\gr_{{\mathbb F}(Q)} Q$ is
canonically isomorphic to the graded fraction ring
$Q_{gr}(\gr_{\mathbb F} A)$ as Poisson $w$-graded algebras.  
\end{enumerate}
\end{lemma}

\begin{proof}
(1) By Lemma \ref{xxlem2.6}, it suffices to show that $\deg_Q$ 
is a $w$-valuation on $Q$. Recall that $\deg$ is the 
degree function on $A$ defined by \eqref{E2.3.1}. It is clear 
that $\deg_Q(a)=\deg(a)$ when $a\in A$. Conditions (a,b,c) in
Definition \ref{xxdef1.1}(1) are clear for $\deg_Q$. For
condition (d), let $a,b\in Q$. We can write $a=xy^{-1}$
and $b=zy^{-1}$. Then
$$\begin{aligned}
\deg_Q(a+b)&=\deg_Q(xy^{-1}+zy^{-1})=\deg_Q((x+z)y^{-1})\\
&=\deg(x+z)-\deg(y)\geq \min\{\deg(x),\deg(z)\}-\deg(y)\\
&=\min\{\deg(x)-\deg(y),\deg(z)-\deg(y)\}\\
&=\min\{\deg_Q(a),\deg_Q(b)\}
\end{aligned}
$$
where the $``="$ holds when $\deg(x)\neq \deg(z)$,
or equivalently, $\deg_Q(a)\neq \deg_Q(b)$. Similarly, 
for condition (e) in Definition \ref{xxdef1.1}(3), we have
$$\begin{aligned}
\deg_Q(\{a,b\})&=\deg_Q(\{xy^{-1},zy^{-1}\})\\
&=\deg_Q(\{x,z\}y^{-2}-zy^{-3}\{x,y\}-xy^{-3}\{y,z\})\\
&\geq \min\{ \deg_Q(\{x,z\}y^{-2}),
\deg_Q(zy^{-3}\{x,y\}),\deg_Q(xy^{-3}\{y,z\})\}\\
&=\min\{\deg(\{x,z\})-2\deg(y),
\deg(z)-3\deg(y)+\deg(\{x,y\}),\\
&\quad\quad \qquad
\deg(x)-3\deg(y)+\deg(\{y,z\})\}\\
&\geq  \deg(x)+\deg(z)-w-2\deg(y)\\
&=\deg_Q(a)+\deg_Q(b)-w.
\end{aligned}
$$
Therefore, the assertion holds.

(2) For any $x\in F_i(A)\setminus F_{i+1}(A)$, let 
$\overline{x}$ be $x+F_{i+1}(A)$. Similarly, for 
$y\in F_i(Q)\setminus F_{i+1}(Q)$, let $\overline{y}$ be 
$y+F_{i+1}(Q)$. Then $\overline{x} \to \overline{x}$ (where 
the second $\overline{x}$ is in $\gr_{{\mathbb F}(Q)}Q$) for 
all $x\in A$ induces a canonical inclusion 
$\gr_{\mathbb F} A\subseteq \gr_{{\mathbb F}(Q)} Q$. 

For every nonzero element $y\in Q$, the equation $y y^{-1}=1$ in $Q$ 
implies that $\overline{y} \overline{y^{-1}}=1$ in 
$\gr_{{\mathbb F}(Q)} Q$. Hence $\gr_{{\mathbb F}(Q)} Q$ is a 
graded field. As a consequence, $Q_{gr}(\gr_{{\mathbb F}} A)
\subseteq \gr_{{\mathbb F}(Q)} Q$. To prove 
$Q_{gr}(\gr_{{\mathbb F}} A)\supseteq 
\gr_{{\mathbb F}(Q)} Q$, we let $\overline{y}
\in \gr_{{\mathbb F}(Q)} Q$ where $y=ab^{-1}$
for some $a,b\in A$. Then $a=yb$ in $Q$. This
induces $\overline{a}=\overline{y}\overline{b}$ in $\gr_{{\mathbb F}(Q)} Q$. As a consequence,
$\overline{y}=\overline{a}\overline{b}^{-1}$.
Since both $\overline{a}$ and $\overline{b}$
are in $\gr_{{\mathbb F}} A$, $\overline{y}$ is in
$Q_{gr}(\gr_{{\mathbb F}} A)$ as desired.
The above natural isomorphism preserves the Poisson $w$-graded structures, and hence the assertion follows.
\end{proof}

By Lemma \ref{xxlem2.7},  the induced good filtration 
${\mathbb F}(Q)$ on $Q:=Q(A)$ is completely determined by 
the good filtration ${\mathbb F}$ on $A$, and vice versa. 
Combining with Lemma \ref{xxlem2.6}, one sees that a valuation
$\nu$ on $Q$ is completely determined by its restriction on 
$A$, and vice versa. If ${\mathbb F}$ is the filtration 
defined by a valuation $\nu$, we also use $\gr_{\nu}(A)$ for 
$\gr_{\mathbb F} (A)$.

\begin{lemma}
\label{xxlem2.8}
Let ${\mathbb F}$ and ${\mathbb G}$ be filtrations of $A$. 
Suppose ${\mathbb G}$ is a subfiltration of ${\mathbb F}$, 
namely, $G_i(A)\subseteq F_i(A)$ for all $i\in {\mathbb Z}$.
\begin{enumerate}
\item[(1)]
There is a naturally graded algebra homomorphism 
$\phi:\gr_{\mathbb G} A\to \gr_{\mathbb F} A$.
\item[(2)]
${\mathbb G}={\mathbb F}$ if and only if $\phi$ is injective.
\end{enumerate}
\end{lemma}

\begin{proof} 
(1) We define a $\Bbbk$-linear map $\phi: (\gr_{\mathbb G} A)_i
\to (\gr_{\mathbb F} A)_i$ by sending 
\begin{equation}
\label{E2.8.1}\tag{E2.8.1}
\widetilde{x}:=x + G_{i+1}(A) \quad \mapsto \quad 
\overline{x}:=x+F_{i+1}(A)
\end{equation}
for all $x\in G_{i}(A)$. Since ${\mathbb G}$ is a subfiltration 
of ${\mathbb F}$, $\phi$ is well-defined. It is clear that
$\phi$ is additive concerning homogeneous elements of the 
same degree. Let $\widetilde{x}\in (\gr_{\mathbb G} A)_i$ and 
$\widetilde{y}\in (\gr_{\mathbb G} A)_j$ for $i,j\in {\mathbb Z}$.
Then 
$$\begin{aligned}
\phi(\widetilde{x}\widetilde{y})
&=\phi((x+ G_{i+1}(A))(y+G_{j+1}(A)))=\phi(xy+G_{i+j+1}(A))\\
&=xy+F_{i+j+1}(A)=(x+F_{i+1}(A))(y+F_{j+1}(A))\\
&=\phi(\widetilde{x})\phi(\widetilde{y}).
\end{aligned}
$$
Thus, $\phi$ is a graded algebra homomorphism. Therefore, 
\eqref{E2.8.1} induces a graded algebra homomorphism as desired.

(2) One direction is clear. For the other direction, suppose
that ${\mathbb G}\neq {\mathbb F}$. Then there is some nonzero 
element $f\in F_j(A)\setminus G_j(A)$ for some $j$. Let $i$ be 
the integer such that $f\in G_i(A)\setminus G_{i+1}(A)$. 
Since $f\not\in G_j(A)$ and $f\in G_{i}(A)$, we have $i<j$. 
Consider the algebra map $\phi$ defined in part (1):
$$\phi: (\gr_{\mathbb G}(A))_i\to (\gr_{\mathbb F}(A))_i$$
which sends 
$$\widetilde{f}:=f+ G_{i+1}(A)\mapsto \overline{f}:=f+F_{i+1}(A).$$
Since $i<j$ or $i+1\leq j$, $f\in F_{j}(A)\subseteq F_{i+1}(A)$. 
This means that $\phi(\widetilde{f})=\overline{f}=0$. By the choice
of $i$, $\widetilde{f}\neq 0$. Therefore, $\phi$ is not injective 
as desired.
\end{proof}

Now, we give a general setup for the rest of this section.
Suppose that an algebra $A$ is generated by elements $\{x_k\}_k$ 
where $k$ is in a fixed index set, and we are given a degree 
assignment on the generators and $1\in A$ as follows:
\begin{equation}
\label{E2.8.2}\tag{E2.8.2}
\deg(1)=0 \quad{\text{and}} \quad 
\deg(x_k)=d_k\in {\mathbb Z} \quad
{\text{for all}} \quad k.
\end{equation}
Then, we can define a descending chain of vector spaces 
${\mathbb F}^{ind}:=\{F_i^{ind}\}_{i\in {\mathbb Z}}$ of $A$ by
\begin{equation}
\label{E2.8.3}\tag{E2.8.3}
F_i^{ind}:={\text{the span of all monomials $\prod_{k} x_k^{n_k}$ where $\sum_k n_k d_k\geq i$ }}
\end{equation}
with only finitely many 
$n_k\neq 0$. It is easy to see that ${\mathbb F}^{ind}$ satisfies 
Definition \ref{xxdef2.2}(b,d) and the first part of (a). It is 
not clear whether Definition \ref{xxdef2.2}(c) and the condition 
that $1\in F^{ind}_0 \setminus F^{ind}_1$ hold. However, this 
paper mainly considers when ${\mathbb F}^{ind}$ is a subfiltration of another filtration $\mathbb F$ of $A$. That is, the degree assignment on $x_k$ is coming from a prescribed filtration $\mathbb F$ on $A$.
Hence Definition \ref{xxdef2.2}(a) holds and 
$$\bigcap_i F^{ind}_i(A)\subseteq \bigcap_i F_i(A)=0.$$ 
Therefore Definition \ref{xxdef2.2}(c) also holds, and whence 
${\mathbb F}^{ind}$ is a filtration of $A$. In this case, 
${\mathbb F}^{ind}$ is called an {\it induced} filtration 
determined by the degree assignment given in \eqref{E2.8.2}. 
Note that the condition in Definition \ref{xxdef2.2}(c) holds 
automatically when $A$ is a domain and $\GKdim 
\gr_{{\mathbb F}^{ind}} A=\GKdim A<\infty$. It is clear that 
${\mathbb F}^{ind}$ is uniquely determined by the degree 
assignment given in \eqref{E2.8.2}. Whenever we say ${\mathbb F}^{ind}$ is an induced filtration, we always assume that it is a filtration satisfying Definition \ref{xxdef2.2}(a,b,c,d).

\begin{lemma}
\label{xxlem2.9}
Suppose that $A$ is generated by a set $\{x_k\}_k$ and that we 
are given a degree assignment $\deg$ as in \eqref{E2.8.2}. 
Write $\nu:=\deg$.
\begin{enumerate}
\item[(1)]
Let ${\mathbb F}^{ind}$ be the induced filtration defined by 
\eqref{E2.8.3}. Then $\gr_{{\mathbb F}^{ind}} A$ is generated 
by $\{\widetilde{x_k}\}_k$ where 
$\widetilde{x_k}:=x_k+F^{ind}_{\nu(x_k)+1}(A)$.
  \end{enumerate}
 In parts (2, 3), let ${\mathbb F}$ be a filtration of $A$.
 \begin{enumerate}
\item[(2)]
Suppose that $\deg$ is the degree function corresponding to 
${\mathbb F}$ via \eqref{E2.3.1} and that ${\mathbb F}^{ind}$ is defined by using \eqref{E2.8.2}-\eqref{E2.8.3} based on 
the given generating set $\{x_k\}_k$. Then, there is a natural 
graded algebra homomorphism 
\begin{equation}
\label{E2.9.1}\tag{E2.9.1}
\phi^{ind}: \gr_{{\mathbb F}^{ind}} A\to \gr_{\mathbb F} A
\end{equation}
determined by $\phi^{ind}(\widetilde{x_k})=\overline{x_k}$
for all $k$, where $\overline{x_k}:=x_k+F_{\nu(x_k)+1}(A)$.
\item[(3)]
Suppose that $A$ is a Poisson algebra and that 
${\mathbb F}^{ind}$ is the induced filtration defined as 
\eqref{E2.8.3}. Then ${\mathbb F}^{ind}$ is a $w$-filtration 
if and only if 
\begin{equation}
\label{E2.9.2}\tag{E2.9.2}
\deg(\{x_k,x_l\})\geq \deg(x_k)+\deg(x_l)-w
\end{equation} 
for all $k,l$.
\end{enumerate}
\end{lemma}

\begin{proof} (1) Let $\widetilde{f}$ be any nonzero element in 
$(\gr_{{\mathbb F}^{ind}} A)_i$. We may assume that 
$f\in F^{ind}_{i}\setminus F^{ind}_{i+1}$. Modulo $F^{ind}_{i+1}$, 
we may assume that $f$ is a linear combination of monomials 
$\prod_k x_k^{n_k}$ with $\sum_k n_k d_k=i$. Hence
$\widetilde{f}$ is a linear combination of 
$\widetilde{\prod_k x_k^{n_k}}$ with $\sum_k n_k d_k=i$ for 
finitely many $n_k\neq 0$. Since $\widetilde{\prod_k x_k^{n_k}}
=\prod_{k} \widetilde{x_k}^{n_k}$, the assertion follows.

(2) It follows from the definition that ${\mathbb F}^{ind}$ is a 
subfiltration of ${\mathbb F}$. By Lemma \ref{xxlem2.8}(1), there 
is a natural graded algebra homomorphism $\phi^{ind}$ as given in 
\eqref{E2.9.1}. By part (1), $\gr_{{\mathbb F}^{ind}}$ is 
generated by $\{\widetilde{x_k}\}_k$. Therefore, $\phi^{ind}$ is 
completely determined by $\phi^{ind}(\widetilde{x_k})
=\overline{x_k}$ for all $k$.

(3) One direction is clear. For the other direction, we assume 
that \eqref{E2.9.2} holds. Now let $f\in A$ with degree $a$ and 
$g\in A$ with degree $b$. Write $f=\sum_{\bf i} c_{\bf i} 
x_{i_1}\cdots x_{i_m}$ where $\sum_{s=1}^m \deg(x_{i_s})\geq a$
and $g=\sum_{\bf j} d_{\bf j} 
x_{j_1}\cdots x_{j_n}$ where $\sum_{s=1}^n \deg(x_{j_s})\geq b$.
By \eqref{E2.9.2} and the Poisson structure, 
$$\begin{aligned}
\deg(\{x_{i_1}\cdots x_{i_m},\,&x_{j_1}\cdots x_{j_n}\})\\
&\geq \min_{\alpha,\beta}
\{ \deg(x_{i_1}\cdots\widehat{ x_{i_\alpha}}\cdots x_{i_m}\cdot\,
x_{j_1}\cdots \widehat{x_{j_{\beta}}}\cdots x_{j_n}
\{x_{i_\alpha}, x_{j_{\beta}}\})\}\\
&\geq \deg(x_{i_1}\cdots x_{i_m})+\deg(x_{j_1}\cdots x_{j_n})
-w\\
&\geq a+b-w.
\end{aligned}
$$
Then
$$\begin{aligned}
\deg(\{f,g\})
&=\deg(\{\sum_{\bf i} c_{\bf i} x_{i_1}\cdots x_{i_m},
\sum_{\bf j} d_{\bf j} x_{j_1}\cdots x_{j_n}\})\\
&\geq \min_{{\bf i},{\bf j}}
\{\deg(\{x_{i_1}\cdots x_{i_m},x_{j_1}\cdots x_{j_n}\})\}\\
&\geq a+b-w=\deg(f)+\deg(g)-w
\end{aligned}
$$
as required.
\end{proof}

Lemma \ref{xxlem2.9}(3) offers a straightforward approach to verify whether a filtration is a $w$-filtration. The following is one of the well-known Poisson fields.

\begin{example}
\label{xxexa2.10}
Let $W$ be the Weyl Poisson polynomial ring $\Bbbk[x,y]$ whose 
Poisson structure is determined by
\begin{equation}
\label{E2.10.1}\tag{E2.10.1}
\{x,y\}=1.
\end{equation}
Let $K_{Weyl}$ denote the Weyl Poisson field which is
$Q(W)$, see Example \ref{xxexa0.5}. 

Let $\deg(x)=1$ and $\deg(y)=(w-1)$. Then, the induced filtration 
${\mathbb F}^{ind}$ is good. Since $\Bbbk[x,y]$ is generated 
by $x$ and $y$ and 
$$\deg(\{x,y\})=\deg(1)=0=1+(w-1)-w=\deg(x)+\deg(y)-w,$$
Lemma \ref{xxlem2.9}(3) implies that ${\mathbb F}^{ind}$ is a 
good $w$-filtration. One can also check that 
$\gr_{{\mathbb F}^{ind}}W \cong W$. Hence 
${\mathcal V}_{ncw}(K_{Weyl})\neq 
\emptyset$ and $\alpha_{ncw}(K_{Weyl})\neq 0$ for all 
$w\in {\mathbb Z}$. (In fact ${\mathcal V}_{fw}(K_{Weyl})\neq 
\emptyset$ for all $w$, see Definition \ref{xxdef3.1}(1)). As 
a consequence, $\alpha_{w}(K_{Weyl})\neq 0$ for all 
$w\in {\mathbb Z}$.

Next, we construct infinitely many different nontrivial nonclassical $w$-valuations on $K_{Weyl}$ when $w\leq 0$. For any $\xi\in \Bbbk$, set $f_{\xi}=x+y+\xi y^2$. Then $\{f_{\xi},y\}=1$. Define 
$\deg(f_{\xi})=1$ and $\deg(y)=w-1\leq -1$. By the above 
paragraph, the induced filtration is a $w$-filtration, and 
hence the associated $w$-valuation, denoted by $\nu_{\xi}$, 
is nontrivial and nonclassical. We claim that 
$\nu_{\xi}\neq \nu_{\xi'}$ whenever $\xi\neq \xi'$. This follows 
from the fact $\nu_{\xi}(f_{\xi})=1$ and 
$$\nu_{\xi'}(f_{\xi})=
\nu_{\xi'}(f_{\xi'}+(\xi-\xi')y^2)=
\min\{\nu_{\xi'}(f_{\xi'}),\nu_{\xi'}
(y^2)\}=\min\{1, 2(w-1)\}=2(w-1).$$

If $\Bbbk$ is uncountable, then there exist uncountably many nontrivial nonclassical $w$-valuations for $K_{Weyl}$ for each fixed $w\leq 0$ , namely $\{\nu_{\xi}\}_{\xi\in \Bbbk}$ as described in the previous paragraph. 
\end{example}

Next, we proceed to two crucial lemmas, which will have numerous applications.

\begin{lemma}
\label{xxlem2.11}
Let $P$ be a domain of finite GK-dimension, say $d$, and 
generated by $\{x_k\}_{k\in S}$ for an index set $S$. Let 
$\nu$ be a valuation on $P$ and ${\mathbb F}$ be the 
filtration associated to $\nu$. In parts {\rm{(1,2)}} let 
${\mathbb F}^{ind}$ be the induced filtration determined by 
\eqref{E2.8.3} and the degree assignment $\deg(x_k):=\nu(x_k)$ 
for all $k\in S$. Let $I$ be the image of $\phi^{ind}$, or 
the subalgebra of $\gr_{\mathbb F}(P)$ generated by 
$\{{\overline {x_k}}\}_k$ where $\overline{x_k}$ is 
$x_k+F_{\nu(x_k)+1}(P)$ for every $k$. 
\begin{enumerate}
\item[(1)]
Suppose
\begin{enumerate}
\item[(1a)]
$\GKdim I=d$,
\item[(1b)]
$\gr_{{\mathbb F}^{ind}}(P)$ is a domain.
\end{enumerate}
Then ${\mathbb F}$ agrees with the filtration ${\mathbb F}^{ind}$.
As a consequence, there are natural isomorphisms of graded algebras
$\gr_{{\mathbb F}^{ind}}(P)\cong I\cong \gr_{\mathbb F}(P)$. 
\item[(2)]
Suppose
\begin{enumerate}
\item[(2a)]
all $\{\nu(x_k)\}_k$ are zero,
\item[(2b)]
$\GKdim I=d$.
\end{enumerate}
Then, $\nu$ is trivial.
\item[(3)]
Suppose
\begin{enumerate}
\item[(3a)]
all $\{\nu(x_k)\}_k$ are nonnegative,
\item[(3b)]
$\omega$ is a nonzero element in $P$ such that the factor ring 
$P/(\omega)$ is a domain of GK-dimension $d-1$,
\item[(3c)]
$\GKdim F_0(P)/F_1(P)=d-1$ and $\omega\in F_1(P)$.
\end{enumerate}
Then ${\mathbb F}$ agrees with the new induced filtration 
${\mathbb F}^{new}$ determined by the degree assignment on
the generating set $\{x_k\}_{k\in S}\cup \{w\}$:
$$\{\deg_{new}(x_k):=\nu(x_k)\}_{k\in S} 
\cup \{\deg_{new}(\omega):=\nu(\omega)>0\}$$ 
and $\gr_{\mathbb F} P\cong (P/(\omega))[\overline{\omega}]$.
\end{enumerate}
\end{lemma}

\begin{proof}
(1) By definition, ${\mathbb F}^{ind}$ is a subfiltration of ${\mathbb F}$. By Lemma \ref{xxlem2.9}(2), there is an algebra homomorphism $\phi^{ind}:\gr_{{\mathbb F}^{ind}}P\to \gr_{{\mathbb F}} P$. By Lemma \ref{xxlem2.8}(2) it suffices to show that $\phi^{ind}$ is 
injective. By definition, $I$ is the image 
of the map $\phi^{ind}$, so we have a surjective graded
algebra homomorphism 
$$\phi^{ind}: \gr_{{\mathbb F}^{ind}} P\to I.$$
Suppose to the contrary that $\phi^{ind}$ is not injective. 
Since $\gr_{{\mathbb F}^{ind}} P$ is a domain, we obtain that
$$\GKdim I < \GKdim \gr_{{\mathbb F}^{ind}} P \leq \GKdim P=d$$
which contradicts (1a). Therefore $\phi^{ind}$ is 
injective. The consequence is clear.

(2) Since all $\nu(x_k)=0$, $F_0(P)=P$ and $I$ is a subalgebra 
of $(\gr_{\mathbb F} P)_0=F_0(P)/F_1(P)$. Since $\GKdim I=d=
\GKdim P=\GKdim F_0(P)$, we have $F_1(P)=0$. Thus, the assertion
follows.

(3) We claim that there is a positive integer $h$ such that
$$F_i(P)=\begin{cases} P & i\leq 0,\\
\omega^{\lceil i/h \rceil} P & i>0.\end{cases}
$$
Since $\nu(x_k)\geq 0$ for all $i$, $F_0(P)=P$. Consequently,
$F_i(P)=P$ for all $i\leq 0$. For $i>0$, we use induction.
By (3c), $F_0(P)/F_1(P)$ is a homomorphic image
of $P/(\omega)$. Since $\GKdim F_0(P)/F_1(P)=d-1
=\GKdim P/(\omega)$ and $P/(\omega)$ is a domain, the map
$P/(\omega)\to F_0(P)/F_1(P)$ is an isomorphism. Thus 
$F_1(P)=(\omega)=\omega P$. Hence, the claim holds for $i=1$.

Now let $h$ be the largest positive integer such that
$F_h(P)=F_1(P)$. Then, the claim holds for all $1\leq i\leq h$.
Next, we consider $i=h+1$. By definition 
$F_h(P)\supsetneq F_{h+1}(P)$ and $\omega\not\in F_{h+1}(P)$. 
Consequently, $\nu(\omega)=h$. For every element 
$f\in F_{h+j}(P)$ where $j\geq 1$, we can write it 
as $\omega g$ since $F_{h+j}(P)\subseteq F_h(P)=\omega P$. 
Then 
$$\nu(g)=\nu(\omega g)-\nu(w)\geq h+j-h=j$$
or equivalently $g\in F_j(P)$. Conversely, if $g\in F_j(P)$,
then it is clear that $\omega g\in F_{h+j}(P)$. This means 
that $F_{h+j}(P)=\omega F_j(P)$ for all $j>0$. In particular 
$F_{h+1}(P)=\omega F_1(P)=\omega^2 P$. So the claim holds for 
$i=h+1$. The claim for general $i=h+j$ with $j>0$ follows from 
the induction and the equation $F_{h+j}(P)= \omega F_j(P)$. 
The claim easily implies the main assertion in part (3).
\end{proof}

\begin{lemma}
\label{xxlem2.12}
Let $P$ be a Poisson domain generated by $\{x_k\}_{k\in S}$ 
for an index set $S$. Let $\nu$ be a $w$-valuation on 
$P$ and ${\mathbb F}$ be the filtration associated to $\nu$. 
\begin{enumerate}
\item[(1)]
Let ${\mathbb F}^{ind}$ be the induced filtration of $P$ 
determined by $\deg(x_k):=\nu(x_k)$ for all $k$. If 
${\mathbb F}^{ind}$ is a $w$-filtration,
then $\phi^{ind}: \gr_{{\mathbb F}^{ind}} P\to 
\gr_{{\mathbb F}} P$ is a Poisson algebra homomorphism.
\item[(2)]
Suppose that $P$ becomes a Poisson $w$-graded algebra 
generated by homogeneous elements $\{x_k\}_{k\in S}$ with 
the degree assignment $\deg(x_k):=\nu(x_k)$ for all $k\in S$. 
Then $\phi^{ind}: P\to \gr_{{\mathbb F}} P$ is a Poisson 
algebra homomorphism.
\end{enumerate}
\end{lemma}

\begin{proof} (1) Let $\phi$ be $\phi^{ind}$. By Lemma 
\ref{xxlem2.9}(2), $\phi$ is a graded algebra homomorphism.
Since both $\gr_{{\mathbb F}^{ind}} P$ and $\gr_{\mathbb F}P$
are Poisson $w$-graded algebras by Lemma \ref{xxlem2.3}, it 
remains to prove that $\phi$ preserves the Poisson bracket. 
Since $\gr_{{\mathbb F}^{ind}} P$ is generated by 
$\widetilde{x_k}$ where $\widetilde{x_k}=x_k+
F^{ind}_{\nu(x_k)+1}$ [Lemma \ref{xxlem2.9}(1)], it 
suffices to show that
\begin{equation}
\label{E2.12.1}\tag{E2.12.1}
\phi(\{\widetilde{x_k}, \widetilde{x_l}\})
=\{\phi(\widetilde{x_k}), \phi(\widetilde{x_l})\}
\end{equation}
for all $k,l$. Write $n_k=\nu(x_k)=\deg(x_k)$ for all $k$. 
Then $\phi(\widetilde{x_k})=\overline{x_k}=x_k+F_{n_k+1}$ 
by definition. Suppose that, for all $k,l\in S$,
$$\{x_k,x_l\}=\sum c^{k,l}_{\bf i} x_{i_1}\cdots x_{i_s}
+\sum d^{k,l}_{\bf j} x_{j_1}\cdots x_{j_t} $$
where $c^{k,l}_{\bf i}\neq 0$ only when
$\sum_{\alpha=1}^s \deg(x_{i_\alpha})=\deg(x_k)+\deg(x_l)-w$
and where 
$\sum_{\beta=1}^t \deg(x_{j_{\beta}})>\deg(x_k)+\deg(x_l)-w$. 
As a consequence, 
$$\widetilde{\{x_k,x_l\}}
=\widetilde{\sum c^{k,l}_{\bf i} x_{i_1}\cdots x_{i_s}}
=\sum c^{k,l}_{\bf i} \widetilde{x_{i_1}}\cdots \widetilde{x_{i_s}}.
$$
Now we calculate to get 
$$\begin{aligned}
\{\phi(\widetilde{x_k}),\phi(\widetilde{x_l})\}
&=\{x_k+F_{n_k+1}, x_l+F_{n_l+1}\}=\{x_k,x_l\}+F_{n_k+n_l-w+1}\\
&=\sum c^{k,l}_{\bf i} x_{i_1}\cdots x_{i_s}
   +\sum d^{k,l}_{\bf j} x_{j_1}\cdots x_{j_t} +F_{n_k+n_l-w+1}\\
&=\sum c^{k,l}_{\bf i} x_{i_1}\cdots x_{i_s}+F_{n_k+n_l-w+1}\\
&=\sum c^{k,l}_{\bf i} (x_{i_1}+F_{\nu(x_{i_1})+1})\cdots 
(x_{i_s} +F_{\nu(x_{i_s})+1})\\
&=\sum c^{k,l}_{\bf i} \phi(\widetilde{x_{i_1}})\cdots 
  \phi(\widetilde{x_{i_s}})\\
&=\phi(\sum c^{k,l}_{\bf i} \widetilde{x_{i_1}}\cdots 
  \widetilde{x_{i_s}})\\
&=\phi(\widetilde{\{x_k,x_l\}})
  =\phi(\{\widetilde{x_k}, \widetilde{x_l}\})
\end{aligned}
$$
as required.

(2) This follows from part (1) since $P$ is canonically isomorphic
to $\gr_{{\mathbb F}^{ind}} P$ under the hypotheses.
\end{proof}

\section{Faithful valuations and examples}
\label{xxsec3}
In this section, we will work out the complete set of $0$-valuations for several Poisson fields. Some of these 
$0$-valuations are unique in the following sense.

\begin{definition}
\label{xxdef3.1}
Let $K$ be a Poisson field. 
\begin{enumerate}
\item[(1)]
A $w$-valuation $\nu$ on $K$ is called a {\it faithful 
$w$-valuation} if the following hold.
\begin{enumerate}
\item[(a)]
The image of $\nu$ is ${\mathbb Z}\cup\{\infty\}$.
\item[(b)]
$\nu$ is nondegenerate.
\item[(c)]
$\nu$ is nonclassical.
\end{enumerate}
\item[(2)]
Let ${\mathcal V}_{fw}(K)$ be the set of all faithful $w$-valuations on $K$.
\item[(3)]
A $w$-valuation $\nu$ on a Poisson domain $A$ is called 
{\it faithful} if it becomes a faithful $w$-valuation when extended to the Poisson fraction field $Q(A)$.  
\end{enumerate}
\end{definition}

Note that condition (b) in Definition \ref{xxdef3.1}(1) says 
the transcendence degree of the residue field of $\nu$ is one 
less than the transcendence degree of $K$, and condition (c) means the induced Poisson bracket on the associated 
graded algebra $\gr_{\nu}(K)$ is nonzero. In particular, a {\it faithful valuation} introduced in 
Definition \ref{xxdef0.2} can be understood as a faithful $0$-valuation. 

Let $P$ be a Poisson $w$-graded domain. We fix this ${\mathbb Z}$-graded structure on $P$ and call it an 
{\it Adams grading}. We define the {\it Adams$^{Id}$ filtration} 
of $P$, denoted by ${\mathbb F}^{Id}$, by
\begin{equation}
\label{E3.1.1}\tag{E3.1.1}
F^{Id}_i(P):=\oplus_{n\geq i} P_n 
\quad {\text{for all $i\in {\mathbb Z}$}}.
\end{equation}
It is clear that ${\mathbb F}^{Id}$ is a Poisson $w$-filtration 
such that the associated graded ring $\gr_{{\mathbb F}^{Id}} P$ 
is canonically isomorphic to $P$ as Poisson $w$-graded algebras. 
The {\it Adams$^{-Id}$ filtration} of $P$, denoted by ${\mathbb F}^{-Id}$, is defined by
\begin{equation}
\label{E3.1.2}\tag{E3.1.2}
F^{-Id}_i(P):=\oplus_{n\leq -i} P_n
\quad {\text{for all $i\in {\mathbb Z}$}}.
\end{equation}
It is clear that ${\mathbb F}^{-Id}$ is a Poisson 
$(-w)$-filtration such that the associated graded ring 
$\gr_{{\mathbb F}^{-Id}} P$ is isomorphic to $P$ by flipping the 
grading ($i\leftrightarrow -i$). Let $\nu^{Id}$ (resp. 
$\nu^{-Id}$) denote the $w$-valuation (resp. $(-w)$-valuation) 
of $P$ associated to ${\mathbb F}^{Id}$ (resp. ${\mathbb F}^{-Id}$). 
We call both $\nu^{Id}$ and $\nu^{-Id}$ the {\it Adams valuations} 
of $P$. The following lemma is easy.

\begin{lemma}
\label{xxlem3.2}
Let $P$ be a ${\mathbb Z}$-graded domain with $P\neq P_0$.
\begin{enumerate}
\item[(1)]
$P$ has at least two good filtrations: ${\mathbb F}^{Id}$ and 
${\mathbb F}^{-Id}$. In both cases, the associated graded rings are isomorphic to $P$. If $P$ is a Poisson $w$-graded domain, 
then the corresponding $\nu^{Id}$ {\rm{(}}resp. $\nu^{-Id}${\rm{)}} 
is a $w$-valuation {\rm{(}}resp. $(-w)$-valuation{\rm{)}}.
\item[(2)]
If $P$ is a Poisson $0$-graded algebra, then the Poisson field 
$Q(P)$ has at least two distinct $0$-valuations $\nu^{Id}$ and 
$\nu^{-Id}$.
\item[(3)]
Suppose an algebra $A$ has a filtration ${\mathbb F}^{c}$ such 
that $\gr_{{\mathbb F}^c} A$ is canonically isomorphic to $P$ 
as ${\mathbb Z}$-graded algebras. Then $A$ has a valuation
associated to ${\mathbb F}^c$, denoted by $\nu^{c}$. If $A$ is further a Poisson algebra and ${\mathbb F}^c$ is a $w$-filtration, 
then $\nu^{c}$ is a $w$-valuation. 
\end{enumerate}
\end{lemma}

\begin{proof} The proof, which is straightforward, has been omitted.
\end{proof}

We will see that the Weyl Poisson field $K_{Weyl}$ 
[Example \ref{xxexa0.5}] plays a unique role in the study 
of 2-valuations. To achieve our objectives, we will be introducing a few temporary concepts. 

\begin{definition}
\label{xxdef3.3}
Let $K$ be a Poisson field and $\nu$ a $w$-valuation on $K$.
\begin{enumerate}
\item[(1)]
We say $\nu$ is {\it quasi-Adams} if $Q(\gr_{\nu} K)\cong K$.
\item[(2)]
$K$ is called {\it $w$-quasi-Adams} if (a) $K$ admits a 
nontrivial $w$-valuation and (b) every $w$-valuation is 
quasi-Adams.
\item[(3)]
We say $\nu$ is {\it Weyl} if $Q(\gr_{\nu} K)\cong K_{Weyl}$.
\item[(4)]
We say $K$ is {\it $w$-Weyl} if (a) there is a faithful $w$-valuation
of $K$ and (b) every faithful $w$-valuation on $K$ is Weyl.
\end{enumerate}
\end{definition}

We recall some notations introduced in Construction \ref{xxcon0.6} 
together with a review of Adams grading (internal grading). Let 
$\Bbbk[x,y,z]$ be a polynomial ring that has an internal grading, 
called {\it Adams grading}. The generators $x$, $y$, $z$ have 
Adams grading $|x|$, $|y|$ and $|z|$ being integers. If 
$|x|=|y|=|z|=1$, we say $\Bbbk[x,y,z]$ has {\it standard Adams 
grading}. Let $\Omega$ be a homogeneous element of $\Bbbk[x,y,z]$ 
with respect to the Adams grading. The Poisson structure on 
$A_{\Omega}:=\Bbbk[x,y,z]$ is determined by 
$\{x,y\}=\frac{\partial \Omega}{\partial z}$,
$\{y,z\}=\frac{\partial \Omega}{\partial x}$,
and
$\{z,x\}=\frac{\partial \Omega}{\partial y}$.
The above is equivalent to \eqref{E0.6.1}. In this case, $\Omega$ 
is called the {\it potential} of $A_{\Omega}$. Since $\Omega$ is 
in the Poisson center of $A_{\Omega}$, we have a Poisson factor 
ring $P_{\Omega-\xi}:=A_{\Omega}/(\Omega-\xi)$ for all $\xi\in 
\Bbbk$. A special case is when $|\Omega|=|x|+|y|+|z|$, where both $A_{\Omega}$ and $P_{\Omega}$ are Poisson $0$-graded algebras. 

\begin{definition}
\label{xxdef3.4}
A homogeneous element $\Omega\in A:=\Bbbk[x,y,z]$ is called to 
{\it have an isolated singularity at the origin} {\rm{(}}or 
simply {\it have an isolated singularity}{\rm{)}} if $A_{sing}:=
A/(\Omega_x,\Omega_y,\Omega_z)$ is finite-dimensional over 
$\Bbbk$. In this case, we say $\Omega$ is an {\it i.s. potential}.
\end{definition}

A generic homogeneous element $\Omega$ of degree $n>1$ has an isolated singularity. For instance, each $\Omega=x^n+y^n+z^n$ ($n>1$) is an i.s. potential. 

Note that the definition of having isolated singularities 
uses the graded generating set $\{x,y,z\}$. However,
it is well-known that $A_{sing}$ is independent of the choices 
of generating sets $\{x,y,z\}$. Consequently, the concept 
of $\Omega$ having isolated singularities is independent of the 
choices of generating sets $\{x,y,z\}$. We will focus on potentials with isolated singularities for the rest of the paper.

Henceforth, for the sake of simplicity, we will impose 
\begin{hypothesis}
\label{xxhyp3.5}
The base field $\Bbbk$ is of characteristic zero.
\end{hypothesis}

Two filtrations are {\it equivalent} if their associated 
valuations are equivalent, see Definition \ref{xxdef1.1}(5). 
In the next lemma, we use $\{x_1,x_2,x_3\}$ for any graded 
generating set of $A:=\Bbbk[x,y,z]$, which is not 
necessarily $\{x,y,z\}$. 

\begin{lemma}
\label{xxlem3.6}
Let $P$ be $P_{\Omega}$ as defined in Construction \ref{xxcon0.6}. 
Let $\{x_1,x_2,x_3\}$ be a graded generating set of $A$ 
{\rm{(}}resp. $P${\rm{)}}. Suppose that $\nu$ is a $w$-valuation 
on $P$ and that ${\mathbb F}$ is the associated $w$-filtration of $P$. 
Let $I$ be the subalgebra of $\gr_{\mathbb F} P$ generated by 
${\overline{x_1}},{\overline{x_2}},{\overline{x_3}}$.
\begin{enumerate}
\item[(1)]
Suppose that $\Omega$ is an i.s. potential of degree $|x_1|+
|x_2|+|x_3|$ and $w=0$. If $F_0(P)=P$, then either ${\mathbb F}$ is a trivial 
filtration or $F_0/F_1\cong \Bbbk$ and $\nu(x_1),\nu(x_2),
\nu(x_3)$ are all positive. As a consequence, if 
$\nu(x_1)=\nu(x_2)=\nu(x_3)=0$, then ${\mathbb F}$ is trivial. 
\item[(2)]
Suppose $(\nu(x_1),\nu(x_2),\nu(x_3))=d (|x_1|,|x_2|,|x_3|)$ for 
some positive integer $d$. If $\GKdim I=2$, then ${\mathbb F}$ 
is equivalent to the Adams$^{Id}$ filtration ${\mathbb F}^{Id}$.
\item[(3)]
Suppose $(\nu(x),\nu(y),\nu(z))=d (|x_1|,|x_2|,|x_3|)$ for some 
negative integer $d$. If $\GKdim I=2$, then ${\mathbb F}$ is 
equivalent to the Adams$^{-Id}$ filtration ${\mathbb F}^{-Id}$. 
As a consequence, $F_0(P)=\Bbbk$.
\end{enumerate}
\end{lemma}

\begin{proof}
(1) If ${\mathbb F}$ is nontrivial, then $F_1(P)\neq 0$. Since 
${\mathbb F}$ is a Poisson $0$-filtration, $F_1(P)$ is a Poisson 
prime ideal of $P$. Hence, $P/F_1(P)$ is an affine Poisson domain 
of GK-dimension $\leq 1$. So Lemma \ref{xxlem1.6}(1) implies that 
$P/F_1(P)$ has zero Poisson bracket. Then 
$\{x,y\}=\Omega_y,\{y,z\}=\Omega_x,\{z,x\}=\Omega_y\in \{P,P\}
\subseteq F_1(P)$. As a consequence, $P/F_1(P)$ is a quotient 
of $A_{sing}$. Since $\Omega$ has an isolated singularity, 
$A_{sing}$ is a finite-dimensional graded algebra. Note that 
$P/F_1(P)$ is a domain. So $P/F_1(P)=\kk$ and $F_1(P)$ is 
generated by $x,y,z$ (as well as by $x_1,x_2,x_3$). The main 
assertion follows.

Now we assume that $\nu(x_1)=\nu(x_2)=\nu(x_3)=0$. By the 
valuation axiom, $F_0(P)=P$. The consequence follows from the 
main assertion.

(2) Let ${\mathbb F}^{ind}$ be the induced filtration defined 
in \eqref{E2.8.3} determined by the degree assignment 
$(\deg(x_1),\deg(x_2),\deg(x_3))=d (|x_1|,|x_2|,|x_3|)$ where 
$d>0$. Then ${\mathbb F}^{ind}$ is equivalent to 
${\mathbb F}^{Id}$ and $\gr_{{\mathbb F}^{ind}} P=P$ since $P$ 
is graded to start with. The assertion follows from Lemma 
\ref{xxlem2.11}(1).

(3) Let ${\mathbb F}^{ind}$ be the induced filtration defined
in \eqref{E2.8.3} determined by the degree assignment
$(\deg(x_1),\deg(x_2),\deg(x_3))=d (|x_1|,|x_2|,|x_3|)$ where 
$d<0$. Then ${\mathbb F}^{ind}$ is equivalent to 
${\mathbb F}^{-Id}$ and $\gr_{{\mathbb F}^{ind}} P=P$ since 
$P$ is graded to start with. The assertion follows from Lemma 
\ref{xxlem2.11}(1). The consequence is clear.
\end{proof}

The rest of this section is to gain a better understanding of the faithful valuations on some Poisson fraction fields $Q(P_{\Omega-\xi})$ in Construction \ref{xxcon0.6} when $\Omega=x^3+y^3+z^3+\lambda xyz$. Note that such a $\Omega$ has an isolated singularity if and only if $\lambda^3\neq -3^3$. To define the fraction fields of $P_{\Omega-\xi}$, we need the polynomial $\Omega-\xi$ to be irreducible in $A$. When $\xi=0$, we need the condition $\lambda^3\ne -3^3$, which is the same as requiring $\Omega$ being smooth.

\begin{lemma}
\label{xxlem3.7}
Let $K$ be a Poisson field containing nonzero elements $x,y,z$ 
such that $\{x,y\}=3z^2+\lambda xy$, $\{y,z\}=3x^2+\lambda yz$, 
and $\{z,x\}=3y^2+\lambda xz$ for some $\lambda \in \Bbbk$. 
\begin{enumerate}
\item[(1)]
Let $\nu$ be a $0$-valuation on $K$. Then 
$\nu(x)=\nu(y)=\nu(z)$.
\item[(2)]
If $\lambda=0$, then there is no $w$-valuation for all $w<0$.
\end{enumerate}
\end{lemma}

\begin{proof} (1) Let $\nu$ be a $0$-valuation on $K$. Set 
$a=\nu(x)$, $b=\nu(y)$, and $c=\nu(z)$. Since 
$\{x,y\}=3z^2+\lambda xy$, we have
$$2c=\nu(3z^2)=\nu(\{x,y\}-\lambda xy)\geq 
\min\{\nu(\{x,y\}), \nu(xy)\}
\geq \nu(x)+\nu(y)=a+b.$$
Similarly $2a\geq b+c$ and $2b\geq a+c$. These 
inequalities together imply that $a=b=c$.

(2) Let $\nu$ be a $w$-valuation on $K$. Set $a=\nu(x)$, 
$b=\nu(y)$, and $c=\nu(z)$. Since $\{x,y\}=3z^2$ (when 
$\lambda=0$), we have
$$2c=\nu(3z^2)=\nu(\{x,y\})\geq \nu(x)+\nu(y)-w=a+b-w.$$
Similarly $2a\geq b+c-w$ and $2b\geq a+c-w$. Adding these 
three inequalities, we obtain that $2(a+b+c)\geq 2(a+b+c)-3w$ 
which implies that $w\geq 0$.
The assertion follows.
\end{proof}

A Poisson $w$-graded algebra $P$ is called {\it projectively simple}
if (a) $P$ is infinite-dimensional and (b) every nonzero graded Poisson ideal of $P$ is co-finite-dimensional. This is a Poisson analog of the `just infinite dimensional' property in noncommutative algebras \cite{RRZ}.  
When $\Omega$ has an isolated singularity, it follows from Lemma \ref{xxlem1.6} that $P_{\Omega}$ is projectively simple. We are ready to work out all faithful $0$-valuations on some Poisson fields $Q(P_\Omega)$.

\begin{theorem}
\label{xxthm3.8}
Let $P$ be $P_{\Omega}$ where $\Omega=x^3+y^3+z^3+\lambda xyz$ 
where $\lambda^3\neq -3^3$. Let $K$ be the Poisson fraction 
field $Q(P)$. 
\begin{enumerate}
\item[(1)]
$({\mathcal V}_{0}(K)/\sim)=\{\nu^{Id},\nu^{-Id}\}
={\mathcal V}_{f0}(K)$. Consequently, every nontrivial 
$0$-valuation of $K$ is equivalent to a faithful 
$0$-valuation. Further, $\alpha_0(K)=\alpha_{nc0}(K)=2$ and 
$\alpha_{c0}(K)=\alpha_{-1}(K)=0$.
\item[(2)]
Let $\nu$ be in ${\mathcal V}_0(K)$. Then $Q(\gr_{\nu}(K))
\cong K$.
\item[(3)]
$\bfd(K)=0$, $\bfw(K)=1$ and $K$ is quasi-Adams.
\end{enumerate}
\end{theorem}

\begin{proof}
Recall that the Poisson algebra $P_{\Omega}=A_{\Omega}/(\Omega)$ has Poisson bracket given
by $\{x,y\}=3z^2+\lambda xy$, $\{y,z\}=3x^2+\lambda yz$, 
and $\{z,x\}=3y^2+\lambda xz$.

(1) Let $\nu$ be a nontrivial $0$-valuation on 
$K$ and set $a=\nu(x)$, $b=\nu(y)$, and $c=\nu(z)$. By
Lemma \ref{xxlem3.7}(1), $a=b=c$. By Lemma \ref{xxlem3.6}(1),
$a\neq 0$. So we have two cases to consider.

Case 1: $a>0$. Let $I$ be the subalgebra generated by
${\overline x},{\overline y},{\overline z}$. Since
$P$ is a Poisson $0$-graded algebra with $|x|=|y|=|z|=1$,
${\mathbb F}^{ind}$ is equivalent to ${\mathbb F}^{Id}$.
Then $\phi^{ind}: P\to I\subseteq \gr_{\mathbb F} P$ is a Poisson
algebra homomorphism [Lemma \ref{xxlem2.12}(2)]. 
Since $\GKdim I\geq 1$, $\phi^{ind}$ is injective as $P$ is 
projectively simple. As a consequence, $\GKdim I=2$. By Lemma 
\ref{xxlem2.11}(1), ${\mathbb F}^{ind}={\mathbb F}$. Consequently, 
${\mathbb F}$ is equivalent to ${\mathbb F}^{Id}$ or 
$\nu$ is equivalent to $\nu^{Id}$.

Case 2: $a<0$. Let $I$ be the subalgebra generated by
${\overline x},{\overline y},{\overline z}$. Since
$P$ is a Poisson $0$-graded algebra with $|x|=|y|=|z|=1$,
${\mathbb F}^{ind}$ is equivalent to ${\mathbb F}^{-Id}$.
Then $\phi^{ind}: P\to I\subseteq \gr_{\mathbb F} P$ is a 
Poisson algebra homomorphism [Lemma \ref{xxlem2.12}(2)]. 
Since $\GKdim I\geq 1$, $\phi^{ind}$ is injective as $P$ is 
projectively simple. As a consequence, $\GKdim I=2$. By 
Lemma \ref{xxlem2.11}(1), ${\mathbb F}^{ind}={\mathbb F}$. 
Consequently, ${\mathbb F}$ is equivalent to ${\mathbb F}^{-Id}$ 
or $\nu$ is equivalent to $\nu^{-Id}$.

Therefore, up to equivalence $\sim$, there are only two 
$0$-valuations, namely $\nu^{Id}$ and $\nu^{-Id}$. 
It is easy to see that both $\nu^{Id}$ and $\nu^{-Id}$
are faithful $0$-valuations (consequently, nonclassical and 
nondegenerate since the residue field associated with 
$\nu^{\pm Id}$ is the function field of the smooth elliptic 
curve $\Omega=0$ in $\mathbb P^3$ and hence has transcendence 
degree 1). The assertions follow.

(2) Clear from the proof of part (1) and the facts that 
$\gr_{\nu^{Id}} P\cong P$ and $\gr_{\nu^{-Id}} P\cong P$, 
see discussion before Lemma \ref{xxlem3.2}.

(3) It is an easy consequence of parts (1) and (2).
\end{proof}

Next, we turn to the other case $P_{\Omega-\xi}$ in Construction \ref{xxcon0.6} for $\xi\in \Bbbk^\times$. We will use the following very nice result of Umirbaev-Zhelyabin \cite{UZ}.

\begin{lemma} 
\label{xxlem3.9}
Suppose the {\rm{(}}weighted{\rm{)}} homogeneous element 
$\Omega\in \Bbbk[x,y,z]$ has an isolated singularity. Let 
$\xi\in \Bbbk^{\times}$.
\begin{enumerate}
\item[(1)]
\cite[Theorem 1 and Corollary 2]{UZ}
$P_{\Omega-\xi}$ is a simple Poisson algebra.
\item[(2)]
$\gldim P_{\Omega-\xi}=2$.
\end{enumerate}
\end{lemma}

Part (2) of the above lemma follows quickly from the Jacobian 
criterion for hypersurfaces. 

\begin{lemma}
\label{xxlem3.10}
Suppose that the Adams degrees $|x|$, $|y|$, and $|z|$ in the 
weighted polynomial ring $\Bbbk[x,y,z]$ are all positive and 
that the potential $\Omega$ have degree $|x|+|y|+|z|+w$ for 
some integer $w$. If $\Omega$ is irreducible, then 
$P:=P_{\Omega-\xi}$ has at least one $w$-valuation $\nu^{c}$ as 
given in Lemma {\rm{\ref{xxlem3.2}(3)}}. 
\end{lemma}

\begin{proof}
Note that $P$ has a filtration determined by \eqref{E2.8.3} 
and the degree assignment $\deg(x)=-|x|$, $\deg(y)=-|y|$ and 
$\deg(z)=-|z|$. Or equivalently, we define a filtration 
${\mathbb F}:=\{F_i\}$ of $P$ by 
$$F_{-i}(P)=\{\sum_j a_j f_j\in P\mid a_j\in \Bbbk, f_j 
{\text{ are monomials of Adams degree $\leq i$}}\}$$ 
for all $i$. Then, the associated graded ring is 
$\gr_{\mathbb F} P \cong P_{\Omega}$. Since $\Omega$ is 
irreducible, both $P:=P_{\Omega-\xi}$ and $P_{\Omega}$ are 
domains. Since $|\Omega|=a+b+c+w$, we have 
$\deg(\{x,y\}) =-|\frac{\partial \Omega}{\partial z}|
=-|x|-|y|-w=\deg(x)+\deg(y)-w$. Similarly, $\deg(\{y,z\})
=\deg(y)+\deg(z)-w$ and $\deg(\{x,z\}) =\deg(x)+\deg(z)-w$. 
These equalities together with \eqref{E0.6.1} imply that 
the filtration is a $w$-filtration [Lemma \ref{xxlem2.9}(3)]. 
The assertion follows.
\end{proof}

The following result can be viewed as a twin of Theorem \ref{xxthm3.8}.

\begin{theorem}
\label{xxthm3.11}
Let $P$ be $P_{\Omega-1}$ where $\Omega=x^3+y^3+z^3+\lambda xyz$
where $\lambda^3\neq -3^3$ with standard Adams grading on 
$x,y,z$. Let $K$ be the Poisson fraction field $Q(P)$. 
\begin{enumerate}
\item[(1)]
There is a unique nontrivial $0$-valuation up to equivalence.
As a consequence, ${\mathcal V}_{f0}(K)$ is a singleton. Further
$\alpha_0(K)=\alpha_{nc0}(K)=1$ and 
$\alpha_{c0}(K)=\alpha_{-1}(K)=0$.
\item[(2)]
Let $\nu$ be a nontrivial $0$-valuation on $K$
with associated filtration ${\mathbb F}$. Then, $\nu$ is
determined by the grading of $P$ as given in Lemma 
\ref{xxlem3.10}. As a consequence,
$Q(\gr_{\nu}(K))\cong Q(P_{\Omega})$ or $Q(P)\to_{\nu}
Q(P_{\Omega})$.
\item[(3)]
$\bfd(K)=1$ and $\bfw(K)=1$.
\end{enumerate}
\end{theorem}

\begin{proof}
(1,2) By Lemma \ref{xxlem3.10}, there is a $0$-valuation $\nu$ 
such that $Q(\gr_{\nu}(K))\cong Q(P_{\Omega})$. It remains to
show every nontrivial $0$-valuation is equivalent to the one 
given in the proof of Lemma \ref{xxlem3.10}. The rest of the 
proof is similar to the one of Theorem \ref{xxthm3.8}.

Recall that $P$ is the Poisson algebra $P_{\Omega-1}:=
\Bbbk[x,y,z]/(\Omega-1)$ with $\{x,y\}=3z^2+\lambda xy$,
$\{y,z\}=3x^2+\lambda yz$, and $\{z,x\}=3y^2+\lambda xz$.

Let $\nu$ be a nontrivial $0$-valuation on $K$ and set $a=\nu(x)$, 
$b=\nu(y)$, and $c=\nu(z)$. By Lemma \ref{xxlem3.7}(1), $a=b=c$. 
So, we have three cases to consider.

Case 1: $a=0$. Then $F_0(P)=P$. Since $P$ is Poisson simple 
[Lemma \ref{xxlem3.9}(1)], $F_1(P)=0$. Thus, $\nu$ is trivial, 
yielding a contradiction.

Case 2: $a>0$. This contradicts the fact that $\Omega=1$ in $P$. 

Case 3: $a<0$. Let $I$ be the subalgebra of $\gr_{\nu}(K)$
generated by $\overline{x},\overline{y}$, and $\overline{z}$.
Let ${\mathbb F}^{ind}$ be the induced filtration determined
by $\deg(x)=\deg(y)=\deg(z)=a<0$. Then ${\mathbb F}^{ind}$ is a 
$0$-filtration with $\gr_{{\mathbb F}^{ind}} P\cong P_{\Omega}$.
By Lemma \ref{xxlem2.12}(1), 
$$\phi^{ind}: \gr_{{\mathbb F}^{ind}} P\to I\subseteq 
\gr_{\mathbb F} P$$
is a Poisson algebra morphism. Note that $\GKdim I\geq 1$. Since $\gr_{{\mathbb F}^{ind}} P
(\cong P_{\Omega})$ is projectively simple, $\phi^{ind}$ is injective. As a consequence,
$\GKdim I=2$. By Lemmas \ref{xxlem2.11}(1) and \ref{xxlem3.10},
${\mathbb F}={\mathbb F}^{ind}$. So, we obtain a unique valuation as given in Lemma \ref{xxlem3.10}. 
The assertion is proved.

(3) Let $\nu$ be the unique $0$-valuation on $K$ given in 
part (2). By parts (1,2), $Q(P_{\Omega-1})\to_{\nu}Q(P_{\Omega})$ 
is the only possible arrow. Hence $\bfw(K)=1$. By Theorem 
\ref{xxthm3.8}, there is no arrow $Q(P_{\Omega})\to_{\nu} K$ 
with $K\not\cong Q(P_{\Omega})$. Thus $\bfd(P)=1$. 
\end{proof}

\section{$\gamma$-type invariants}
\label{xxsec4}

An invariant is said to be {\it of $\gamma$-type} if it is defined by using the filtrations ${\mathbb F}^{\nu}$ \eqref{E2.4.1} associated to Poisson valuations $\nu$. This section aims to introduce several $\gamma$-type invariants of Poisson fields and utilize them to classify certain valuations of Poisson algebras/fields defined in Construction \ref{xxcon0.6}. Hypothesis \ref{xxhyp3.5} will be imposed throughout the rest of the paper, although it is not necessary for most of the definitions.

\begin{definition}
\label{xxdef4.1}
Let $K$ be a Poisson field. 
\begin{enumerate}
\item[(1)]
The {\it $^w\Gamma$-cap} of $K$ is defined to be
$$^w\Gamma(K):=\bigcap_{\nu} F^{\nu}_0(K)$$
where $\nu$ runs over all $w$-valuations $\nu$ on $K$.
Since $F^{\nu}_0(K)=K$ if $\nu$ is a trivial $w$-valuation,
it is clear that
\begin{equation}
\label{E4.1.1}\tag{E4.1.1}
^w\Gamma(K)=\bigcap_{\nu \in {\mathcal V}_w(K)} F^{\nu}_0(K).
\end{equation}
\item[(2)]
We say $K$ is {\it $^{w}\Gamma$-normal} if $C:=\,^{w}\Gamma(K)$ is a Poisson subalgebra of $K$ and $C$ is an affine normal domain with $K=Q(C)$. 
\item[(3)]
Let $C$ be a subalgebra of $K$. The {\it $C$-interval} is
defined to be
$$i(C):=
\{w\in {\mathbb Z}\mid \,^{w}\Gamma(K)=C\}.
$$
By Lemma \ref{xxlem4.2}(2) below, if $a,b\in i(C)$, then 
$c\in i(C)$ for any integer $c$ between $a$ and $b$.
\item[(4)]
The $^{\ast}\Gamma$-subalgebra collection of $K$ is defined to be
$$^{\ast}\Gamma(K):=\{^{w}\Gamma(K)\mid w\in {\mathbb Z}\}.$$
\end{enumerate} 
\end{definition}

Note that $^{w}\Gamma$-cap of $K$ will be used in later sections. It is clear that 
\begin{equation}
\label{E4.1.2}\tag{E4.1.2}
^w\Gamma(K)=\{a\in K\mid \nu(a)\geq 0,
\forall \; \nu\in {\mathcal V}_{w}(K)\}.
\end{equation}
The following lemma is straightforward, and its proof 
is omitted.

\begin{lemma}
\label{xxlem4.2}
Let $f: K\to Q$ be a Poisson algebra homomorphism between two 
Poisson fields.
\begin{enumerate}
\item[(1)]
$f$ maps $^w\Gamma(K)$ into $^w\Gamma(Q)$ via restriction for each $w$.
\item[(2)]
$^{w}\Gamma(K)\supseteq \,^{w+1}\Gamma(K)\supseteq \Bbbk$ for each $w$.
\end{enumerate}
\end{lemma}

Despite its technicality, the following theorem is a highly significant general result pertaining to $1$-valuations.

\begin{theorem}
\label{xxthm4.3}
Let $A$ be a Poisson domain of GK-dimension at least 2 
and $K$ be the Poisson fraction field $Q(A)$.
\begin{enumerate}
\item[(1)]
Let $\fp$ be a prime ideal of $A$ of height one such that
$A_{\fp}$ is regular. Then there is a unique nontrivial 
$1$-valuation on $K$, denoted by $\nu^{1,\fp}$, such that 
$F_i(A_{\fp})=\fp^i A_{\fp}$ for $i\ge 0$. In this 
case, $\im (\nu^{1, \fp})=\{\mathbb Z\}\cup\{\infty\}$.
\item[(2)]
{\rm{(}}Controlling theorem{\rm{)}}
If $A$ is a noetherian normal domain, then 
$^1\Gamma(K)\subseteq A$.
\item[(3)]
Suppose $K$ is finitely generated as a field. Then 
there are infinitely many nontrivial $1$-valuations,
namely, $\alpha_1(K)=\infty$.
\end{enumerate}
\end{theorem}

\begin{proof}
(1) Let $B$ be the localization $A_{\mathfrak p}$. Then, $B$ is a regular local ring of Krull dimension 1. Hence $B$ is a 
DVR with its maximal ideal generated by an element $\omega\in B$. 
Note that $B$ is a Poisson algebra. Since $B$ is a DVR, the 
factor ring $B/(\omega)$ is a domain. Define a filtration 
${\mathbb F}^{1,\fp}$ on $B$ as follows:
\begin{equation}
\label{E4.3.1}\tag{E4.3.1}
F^{1,\fp}_i:=\begin{cases} B & i< 0,\\
\omega^i B (=\fp^i A_{\fp}) & i\geq 0.\end{cases}
\end{equation}
It is easy to check that (a) $\gr_{{\mathbb F}^{1,\fp}} B
\cong (B/(\omega))[\overline{\omega}]$ with $\nu(\omega)=1$
and (b) for all $i,j$, $\{F^{1,\fp}_i, F^{1,\fp}_j\}
\subseteq F^{1,\fp}_{i+j-1}$. 
Hence ${\mathbb F}^{1,\fp}$ is a good $1$-filtration which 
produces a $1$-valuation $\nu^{1,\fp}$ on $B$ and whence on $K$
via localization. It is also easy to check that 
$\im (\nu^{1, \fp})=\{\mathbb Z\}\cup\{\infty\}$. The uniqueness follows 
from an idea in the proof of Lemma \ref{xxlem2.11}(3). 

(2) Continue the proof of part (1); we have 
$\nu^{1,\fp}(\omega^{-1})=-1$. Since 
$K=\cup_{i\in {\mathbb Z}} \omega^{i} A_{\fp}$, one sees that 
$F^{1,\fp}_i(K)=\omega^i A_{\fp}$ for all $i\in {\mathbb Z}$. 
Since $A$ is noetherian normal, we have
$$
A=\bigcap_{{\text{height one primes}} \; \fp} A_{\fp}
=\bigcap_{{\text{height one primes}} \; \fp} F_{0}^{1,\fp}(K)
\supseteq
\bigcap_{\nu\in {\mathcal V}_1(K)} F_0^{\nu}(K)
=\,^1\Gamma(K).
$$

(3) Let $V:=\{v_1,\cdots, v_n\}$ be a finite set of generators 
of the field $K$ and let $C$ be the subalgebra of $K$ generated 
by $V$. Then $Q(C)=K$. Write $\{v_i,v_j\}=g_{ij} z^{-1}$ where 
$z, g_{ij} \in A$ for all $i,j$. Adding $z^{-1}$ to the 
set $V$, we may assume that $C$ is a Poisson subalgebra of $K$. 
Localizing further, we may assume that $C$ is regular.

Now $C$ is an affine Poisson domain that is regular as a commutative ring. Let ${\mathfrak p}$ be a prime ideal of height 1. By part (1), there is a $1$-valuation $\nu^{1,\fp}$ on $K=Q(C)$.
By prime avoidance and principal ideal theorem, there are infinitely many such ${\mathfrak p}$. Therefore, we have infinitely many distinct nontrivial $1$-valuations on $K$.
\end{proof}

\begin{remark}
\label{xxrem4.4}
\begin{enumerate}
\item[(1)]
By the above theorem, $\alpha_1(K)$ is infinite. However, when 
$w<1$, $\alpha_w(K)$ could be finite. For example, in Theorem 
\ref{xxthm3.11}, $\alpha_{0}(K)=1$.
\item[(2)]
If ${\mathcal V}_w(K)\neq \emptyset$ for some $w<0$, then 
Lemma \ref{xxlem1.4}(3) says that $w(K)=-\infty$. Otherwise, $w(K)\geq 0$. If $w(K)> 0$, then by Theorem \ref{xxthm4.3}(3), 
$w(K)=1$. In summary, $w(K)$ can only be $-\infty$, $0$ or $1$.
\end{enumerate}
\end{remark}

Next, we will compute various $\gamma$-type invariants for 
a large class of Poisson fields. For the rest of this section, we assume

\begin{hypothesis}
\label{xxhyp4.5}
Let $A=\Bbbk[x,y,z]$ with Adams grading $|x|=|y|=|z|=1$. 
Let $\Omega\in A$ denote a homogeneous element of degree
$3+n$ where $n\geq 1$. 
\end{hypothesis}

Recall that $A_1:=\Bbbk x +\Bbbk y+\Bbbk z$ is the degree 1 part of $A$. 

\begin{lemma}
\label{xxlem4.6}
Let $\nu$ be a valuation of $A$. Then there are $x_1,x_2,x_3
\in A_1$ generating $A$ as an algebra while satisfying 
$\nu(x_1)=\nu(x_2)=\nu(x_3)$.
\end{lemma}

\begin{proof}
Let $w\in A_1$ be a nonzero element such that
$\nu(w)=\min\{\nu(x),\nu(y),\nu(z)\}$. By valuation axioms,
$\nu(w)=\min \{\nu(f)\mid f\in A_1\}$. Let $x_1=x+\lambda w$,
$x_2=y+\lambda w$ and $x_3=z+\lambda w$ for some 
$\lambda\in \Bbbk$. For a generic $\lambda$,
$\nu(x_1)=\nu(x_2)=\nu(x_3)=\nu(w)$ and $\{x_1,x_2,x_3\}$
are linearly independent. The assertion follows.
\end{proof}

The following lemma is straightforward.

\begin{lemma}
\label{xxlem4.7}
Let $\{x_1,x_2,x_3\}$ be as in Lemma \ref{xxlem4.6}.
\begin{enumerate}
\item[(1)]
The $\Bbbk$-span of $\{x_i,x_j\}$ for $1\leq i,j\leq 3$ 
is equal to the $\Bbbk$-span of $\{x,y\},
\{y,z\}$, and $\{z,x\}$.
\item[(2)]
Let $A_{\Omega}$ be defined as Construction \ref{xxcon0.6}.
Then $A_{sing}:=A/(\Omega_x,\Omega_y,\Omega_z)$ is equal to
$A/(\{x_i,x_j\},1\leq i, j\leq 3)$.
\end{enumerate}
\end{lemma}

\begin{lemma}
\label{xxlem4.8}
Let $\Omega$ be a homogeneous element of degree $3+n$ for some 
$n\geq 0$. Suppose $\Omega$ has an isolated singularity. Let 
$B$ be either $A_{\Omega}$, $P_{\Omega}$, or $P_{\Omega-\xi}$
for some $\xi\in \Bbbk^{\times}$. Let $\nu$ be a $w$-valuation
of $K:=Q(B)$ and $d=\min\{\nu(x),\nu(y),\nu(z)\}$. 
\begin{enumerate}
\item[(1)]
Suppose $w<n$. Then $d\geq 0$. As a consequence, 
$^w\Gamma(K)\supseteq B$. 
\item[(2)]
Suppose that $B=P_{\Omega-\xi}$ with $\xi\in \Bbbk^{\times}$ and $n=0$. Then there is no $(-1)$-valuation.
\item[(3)]
Suppose that $B=P_{\Omega-\xi}$ with $\xi\in \Bbbk^{\times}$ and $n>0$. Then, there is no nontrivial $0$-valuation.
\item[(4)]
Let $B$ be $P_{\Omega}$ or $P_{\Omega-\xi}$.
If $n=w>0$, then either $d\geq 0$ or $d=-1$ with ${\mathbb F}={\mathbb F}^{-Id}$ when $B=P_{\Omega}$ and
${\mathbb F}={\mathbb F}^{c}$ when $B=P_{\Omega-\xi}$. 
\item[(5)]
Let $B$ be $P_{\Omega}$ with $n>0$. Then there is no
faithful $0$-valuation on $K$.
\end{enumerate}
\end{lemma}

\begin{proof} (1) Let $\nu$ be any valuation of $K$.
By Lemma \ref{xxlem4.6}, there are linearly 
independent elements $x_1, x_2, x_3\in A_1$ such that 
$d=\nu(x_1)=\nu(x_2)=\nu(x_3)$. Let $\Omega_{k}=\{x_i, x_j\}$ 
where $(i,j,k)$ is either $(1,2,3)$, $(2,3,1)$, or $(3,1,2)$. 
Then the Adams degree of $\Omega_k$ is $2+n$. By the choice 
of $d$ and valuation axioms, we have $\nu(\Omega_k)\geq (2+n)d$ 
for $k=1,2,3$. 

Let $I$ be the subalgebra of $\gr_{\nu} K$ generated by 
$\overline{x_1}, \overline{x_2}$, and $\overline{x_3}$.

Case 1: Suppose $\nu(\Omega_k)>(2+n)d$ for all $k=1,2,3$.
Then $\overline{\Omega_k}=0$ in $I$. Or equivalently,
$\Omega_k(\overline{x_1}, \overline{x_2},\overline{x_3})=0$
in $I$. This means that $I$ is a factor ring of $A_{sing}$
by Lemma \ref{xxlem4.7}(2). This implies that $I=\Bbbk$ and
$d=0$. 

Case 2: Suppose $\nu(\Omega_k)=(2+n)d$ for some $k$. Then 
$$(2+n) d=\nu(\Omega_k)=\nu(\{x_i,x_j\})
\geq \nu(x_i)+\nu(x_j)-w=2 d- w$$
which implies that $nd \geq -w$, or equivalently, $n(-d)\leq w$. 
If $n=0$, then $w\geq 0=n$, yielding a contradiction. Otherwise 
$n>0$ and whence $-d\leq w/n<1$ where the last $<$
follows from the fact $w<n$. Since $d$ is an integer,
$-d\leq 0$, or $d\geq 0$.  

The consequence is clear.

(2) Suppose to the contrary that $\nu$ is a $(-1)$-valuation on 
$K$. By the proof of part (1), only Case 1 can happen. As a 
consequence, $I=\Bbbk$ and $d=0$. In this case $F_0^{\nu}(B)=B$
and $F_0^{\nu}(B)/F_{1}^{\nu}(B)=\Bbbk$, which contradicts Lemma \ref{xxlem3.9}(1).

(3) Let $\nu$ be a $0$-valuation on $K$. By part (1), 
$^0\Gamma(K)\supseteq B$. Consequently, $F^{\nu}_0(B)=B$. Since
$B$ is Poisson simple [Lemma \ref{xxlem3.9}(1)], the Poisson 
ideal $F^{\nu}_1(B)$ is zero. Thus, $\nu$ is trivial.

(4) Retain the notation introduced in the proof of part (1).

Case 1: Suppose $\nu(\Omega_k)>(2+n)d$ for all $k=1,2,3$.
The exact proof shows that $I=\Bbbk$ and that $d=0$. 

Case 2: Suppose $\nu(\Omega_k)=(2+n)d$ for some $k$. Then 
$$(2+n) d=\nu(\Omega_k)=\nu(\{x_i,x_j\})
\geq \nu(x_i)+\nu(x_j)-w=2 d- n$$
which implies that $nd \geq -n$. Since $n> 0$, we have $d\geq -1$.
If $d\geq 0$, we are done. Otherwise, let $d=-1$. Then 
$\overline{\Omega_{k}}\neq 0$, and consequently,
$\{\overline{x_i},\overline{x_j}\}=\overline{\Omega_{k}}
\neq 0$. This shows that $\overline{x_i}$ and $\overline{x_j}$
are algebraically independent by Lemma \ref{xxlem1.6}(1).
So $\GKdim I\geq 2$. By Lemma \ref{xxlem2.11}(1),
${\mathbb F}={\mathbb F}^{ind}$. Since $d=-1$, we further have 
${\mathbb F}={\mathbb F}^{-Id}$ when $B=P_{\Omega}$ and
${\mathbb F}={\mathbb F}^{c}$ when $B=P_{\Omega-\xi}$. 

(5) Suppose to the contrary that $\nu$ is a faithful 
$0$-valuation on $B$. By part (1), $d\geq 0$. 

Case 1: Suppose $d>0$. Let $x_1,x_2,x_3$ be as in the 
proof of part (1). Then $\deg \overline{x_i}=\nu(x_i)
=d>0$. Thus $I$ is not $\Bbbk$ and whence $\GKdim I\geq 1$.
So, Case 1 of the proof of part (1) is impossible, leaving only Case 2 as a possibility. Then $\{\overline{x_i},\overline{x_j}\}=\overline{\Omega_k}\neq 0$
for some $k$. So $\overline{x_i}$ and $\overline{x_j}$ are
algebraically independent by Lemma \ref{xxlem1.6}(1). 
Then $\GKdim I=2$. By Lemma \ref{xxlem2.11}(1),
${\mathbb F}={\mathbb F}^{ind}$ since $B$ is graded to 
start with.

Since $d>0$, for all $(i,j,k)=(1,2,3)$ or its cyclic 
permutations,
$$\nu(\{x_i,x_j\})=\nu(\Omega_k)\geq (2+n) d> 2d 
=\nu(x_i)+\nu(x_j).$$
By Lemma \ref{xxlem2.9}(3), ${\mathbb F}^{ind}$ is a
$(-1)$-valuation. So $\nu$ is not a faithful 
$0$-valuation, yielding a contradiction. 

Case 2: Suppose $d=0$. Since $\nu$ is a faithful 0-valuation, $F_0^\nu(B)/F_1^\nu(B)=B/F_1^\nu(B)$ is an affine Poisson domain of GK-dimension $\leq 1$. So Lemma \ref{xxlem1.6}(1) implies that $B/F_1^\nu(B)$ has zero Poisson bracket. Hence for all $(i, j, k) = (1, 2, 3)$ or its cyclic permutations
\[
\{\overline{x_i},\overline{x_j}\}=\overline{\Omega_k}=0,
\]
or equivalently $\Omega_k\in F_1^\nu(B)$ for all $k$. Note that $\Omega$ is an i.s. potential and $A_{sing}$ is local with residue field $\kk$. Thus, $B/F_1^\nu(B)=\kk$, which implies that  $c_i:=\overline{x_i}\in \Bbbk^{\times}$ such that 
\[
\overline{\Omega_k(x_1,x_2,x_3)}=\Omega_k(\overline{x_1},\overline{x_2}, \overline{x_3})=\Omega_k(c_1,c_2,c_3)=0
\]
for all $k$. But we must have $c_1=c_2=c_3=0$, which contradicts $d=0$.
\end{proof}

We are now able to present the primary application of $\gamma$-invariants.

\begin{theorem}
\label{xxthm4.9}
Let $\Omega$ be a homogeneous element of degree $3+n$ for some 
$n\geq 2$. Suppose $\Omega$ has an isolated singularity. Let 
$B$ be either $A_{\Omega}$, $P_{\Omega}$ or $P_{\Omega-\xi}$
for some $\xi\in \Bbbk^{\times}$. Then, for every $w$ between 
$1$ and $n-1$, $^{w}\Gamma(Q(B))=B$.
\end{theorem}

\begin{proof} Let $K:=Q(B)$. By Lemma \ref{xxlem4.8}(1), 
$^{w}\Gamma(K)\supseteq B$. We can check that $B$ is always normal 
(due to the fact that $A_\Omega$ and $P_{\Omega-\xi}$ are regular and $P_{\Omega}$ is the homogeneous coordinate ring of a smooth irreducible curve). By Theorem \ref{xxthm4.3}(2), $\,^1\Gamma(K)\subseteq B$. 
Hence
$$B\supseteq \,^1\Gamma(K)\supseteq \,^2\Gamma(K)\supseteq \cdots 
\supseteq \,^{n-1}\Gamma(K) \supseteq B.$$
The assertion follows.
\end{proof}

The above theorem will be essential in establishing the Dixmier property of the corresponding Poisson fields $K$. 

\begin{theorem}
\label{xxthm4.10}
Let $\Omega$ be a homogeneous element of degree $3+n$ for 
some $n\geq 1$. Suppose $\Omega$ has an isolated singularity. 
Let $B$ be $P_{\Omega-\xi}$ for some $\xi\in \Bbbk^{\times}$ 
and let $K=Q(B)$. Then there is no nontrivial $0$-valuation of $K$ and 
$$^w\Gamma(K)=\begin{cases} K & w\leq 0,\\
B & 1\leq w\leq n-1,\\
\Bbbk & w\geq n.
\end{cases}$$
As a consequence, $^{\ast}\Gamma(K)=\{K, B, \Bbbk\}$ 
if $n\geq 2$ and $^{\ast}\Gamma(K)=\{K, \Bbbk\}$ if $n=1$.
Further, $i(B)=\{1,\cdots, n-1\}$.
\end{theorem}

\begin{proof}
By Theorem \ref{xxlem4.8}(3), there is no nontrivial $0$-valuation on 
$K$. Hence $^0\Gamma(K)=K$. It immediately follows from Lemma \ref{xxlem4.2}(2) that $^w\Gamma(K)=K$ if $w\leq 0$. If $1\leq w\leq n-1$ (if $n=1$, then there is no such $w$), the assertion follows from Theorem \ref{xxthm4.9}. Now let $w=n$. By Lemma \ref{xxlem3.2}(3),
there is an $n$-valuation $\nu:=\nu^{c}$ such that 
$\nu(x)=\nu(y)=\nu(z)=-1$. Consequently, $F^{\nu}_0(B)=\Bbbk$. 
As a result, we have 
$$^n\Gamma(K)\subseteq \,^1\Gamma(K)\cap F^{\nu}_0(K)
\subseteq B\cap F^{\nu}_0(K)=F^{\nu}_0(B)=\Bbbk.$$
By Lemma \ref{xxlem4.2}(2), $^{w}\Gamma(K)=\Bbbk$ 
for all $w>n$. 

The consequences are clear.
\end{proof}

\section{$\beta$-type invariants and $w$-valuations for small $w$}
\label{xxsec5}
In this section, we introduce/recall $\beta$-type invariants and calculate $w$-valuations on some Poisson fields for $0\leq w\leq 2$. The previous section provided a structure theorem for certain $1$-valuations in Theorem \ref{xxthm4.3}. We will present a similar theorem for specific $2$-valuations in this section.

\begin{definition}
\label{xxdef5.1}
Let $K$ and $Q$ be two Poisson fields and $w$ an integer.
\begin{enumerate}
\item[(1)]
We say $K$ {\it $w$-controls} $Q$ if there is a $w$-valuation 
$\nu$ on $K$ such that 
$$Q(\gr_{\nu}(K))\cong Q.$$ 
We write $K \to_{\nu} Q$ in this case.
\item[(2)]
Let ${\mathcal {PF}}_{d,w}$ be the quiver with vertices 
being Poisson fields of GK-dimension $d$ and arrows 
being $\to_{\nu}$ from $K$ to $Q$ when $K$ $w$-controls $Q$.
\end{enumerate}
\end{definition}

As noted, a $\beta$-type invariant is 
defined using arrows $\to_{\nu}$. For example, 
${\mathcal {PF}}_{d,w}$ is a $\beta$-type invariant. One of 
the projects in \cite{HTWZ2} is trying to understand the quiver 
${\mathcal {PF}}_{2,0}$. 

\begin{definition}
\label{xxdef5.2}
Let $K$ be a Poisson field and $w$ an integer.
\begin{enumerate}
\item[(1)]
The {\it $w$-depth} of $K$ is defined to be
$$\bfd_w(K):=\sup\,\{n \mid 
\{K:=K_0\to_{\nu_1} K_1\to_{\nu_2} K_2\to_{\nu_3} 
\cdots \to_{\nu_{n}}K_n\}
\}$$
where $K_i\not\cong K_j$ for all $0\leq i\neq j\leq n$ and each $\nu_i$ is a faithful $w$-valuation. 
\item[(2)]
The {\it $w$-width} of $K$ is defined to be
$$\bfw_w(K):=\# \left[
\{ Q\mid K\to_{\nu} Q {\text{ where $\nu$ is a faithful 
$w$-valuation}}\}/\cong \right].$$
\end{enumerate}
\end{definition}

When $w=0$, we will omit the subscript $w$ in the above definition as in previous sections. For example, $\bfd_0(K)=\bfd(K)$ (resp. $\bfw_0(K)=\bfw(K)$) as in 
Definition \ref{xxdef5.2}(1) (resp. Definition 
\ref{xxdef5.2}(2)).

The following lemma gives some explicit examples of faithful $1$-valuations, which are special cases of Theorem \ref{xxthm4.3}(1).

\begin{lemma}
\label{xxlem5.3}
Let $\Omega=x^{3+n}+y^{3+n}+z^{3+n}$ where $n\geq 1$ and 
let $P$ be the $P_{\Omega-\xi}$ where $\xi\in \Bbbk^{\times}$. 
\begin{enumerate}
\item[(1)]
The induced filtration ${\mathbb F}^{ind}$ determined by $\deg(x)=\deg(y)=0$ and $\deg(z)=1$ is a good filtration.
\item[(2)]
${\mathbb F}^{ind}$ is a faithful $1$-filtration. As a 
consequence, ${\mathcal V}_{f1}(Q(P))\neq \emptyset$.
\end{enumerate}
\end{lemma}

\begin{proof}
(1) Note that $z$ is a nonzero divisor and $P/(z)\cong 
\Bbbk[x,y]/(x^{3+n}+y^{3+n}-\xi)$. Since 
$x^{3+n}+y^{3+n}-\xi$ is irreducible,
$P/(z)$ is a domain. By definition,
${\mathbb F}^{ind}$ is defined by
$$F_i(P)=\begin{cases} P & i\leq 0,\\
z^i P & i>0.\end{cases}$$
So $\gr_{{\mathbb F}^{ind}} P\cong (P/(z))[\overline{z}]$
which is a domain. The assertion follows.

(2) By definition, $\{x,y\}=(3+n)z^{2+n}\in F_{0+0-1}$,
$\{y,z\}=(3+n)x^{2+n}\in F_{0+1-1}$ and 
$\{z,x\}=(3+n)y^{2+n}\in F_{0+1-1}$. Since $P$ is generated 
by $x,y$ and $z$, by Lemma \ref{xxlem2.9}(3), the 
filtration in part (1) is a $1$-filtration. The faithfulness of $\nu$ is easy to check.
\end{proof}

The subsequent lemma presents an interesting observation regarding $1$-valuations.

\begin{lemma}
\label{xxlem5.4}
Let $P$ be a Poisson noetherian normal domain. Let $\nu$ be a 
valuation on $K:=Q(P)$ such that {\rm{(a)}} $F_0(P)=P$ and 
that {\rm{(b)}} $\GKdim F_0(P)/F_1(P)=\GKdim P-1$. Then, $\nu$ is equivalent to a $1$-valuation.
\end{lemma}

\begin{proof} By definition $F_1(P)$ is a prime ideal of 
$P(=F_0(P))$. Let $S=F_0(P)\setminus F_1(P)$. Then, every element 
in $S$ has $\nu$-value $0$. Then the $\nu$-value of every 
element in $B:=PS^{-1}$ is nonnegative. As a consequence, 
$F_0(B)=B$. Since $F_1(P)\subset P$ is prime of height 1 and $P$ 
is normal, $B$ is a regular local ring of Krull dimension 1.
This means that $B$ is a DVR with a generator $\omega$ of
the maximal ideal of $B$.  

By the proof of Lemma \ref{xxlem2.11}(3), we have that 
there is a positive integer $h$ such that
$$F_i(B)=\begin{cases} B & i\leq 0,\\
\omega^{\lceil i/h \rceil} B & i>0.\end{cases}
$$
Then $\nu$ is equivalent to another valuation, denoted 
by $\nu'$, with its filtration defined by
$$F'_i(B)=\begin{cases} B & i\leq 0,\\
\omega^{i} B & i>0.\end{cases}
$$
By the proof of Theorem \ref{xxthm4.3}(1), $\nu'$ is 
a $1$-valuation. The assertion follows.
\end{proof}

Recall that the {\it geometric genus} of an affine domain $A$
of GK-dimension $1$, denoted by ${\mathfrak g}(A)$, is the 
genus of the normalization of a projective closure of $\Spec A$. 
It is well-known that the geometric genus (or simply called genus) is a birational invariant of $A$. Part (1) of the following 
lemma is due to S{\' a}ndor Kov{\' a}cs. We thank him for allowing us to include his result in our paper.

\begin{lemma}
\label{xxlem5.5}
Let $B$ be an affine domain of GK-dimension 2.
\begin{enumerate}
\item[(1)]
For each positive integer $N$, there exist infinitely many 
integers $g> N$ such that there is a height one prime ideal 
$\fp$ of $B$ with ${\mathfrak g}(B/\fp)=g$.
\item[(2)]
Suppose $B$ has a nonzero Poisson algebra structure. Let 
${\mathbb F}^{ind}$ be the induced filtration of $B_{\fp}$ 
defined by
$$F_i:=\begin{cases} B_{\fp} & i\leq 0,\\
\fp^i B_{\fp} & i\geq 0.
\end{cases}
$$
Then ${\mathbb F}^{ind}$ is a faithful $1$-filtration for 
infinitely many $\fp$ given in part {\rm{(1)}}. 
\end{enumerate}
\end{lemma}

\begin{proof}
(1) We prove it for any affine domain of GK-dimension $d\geq 2$ with the condition: $\mathfrak p$ is a height $d-1$ prime 
ideal of $B$. 

First, we prove the statement for the polynomial ring 
$\Bbbk[x_1,x_2]$ in place of $B$. Indeed, choosing a generic 
irreducible polynomial $f\in \Bbbk[x_1,x_2]$ of degree $> N$ 
provides a quotient $\Bbbk[x_1,x_2]/(f)$ with genus at least $N$. 
By Bertini's theorem we may even assume that $\Bbbk[x_1,x_2]/(f)$ 
is smooth.

This easily implies the statement also for the polynomial ring
$A:=\Bbbk[x_1,\dots, x_d]$ in place of $B$. Indeed, take an $f\in 
\Bbbk[x_1,x_2]$ as above and if $d>2$, then take $\fp:= 
(f,x_3,\dots,x_d)$.

Finally, we will prove this for an arbitrary $B$.  By the Noether 
normalization theorem, there exists an embedding of the polynomial 
ring as above $A\hookrightarrow B$ such that $B$ is integral over 
$A$. Choose a prime ideal $\fp\subset A$ satisfying the desired 
condition. As the embedding is an integral extension, there exists 
a prime ideal ${\mathfrak q}\subset B$ that lies over $\fp$, i.e., 
$A\cap{\mathfrak q}=\fp$, and such that
$\height{\mathfrak q}=\height{\mathfrak p}=d-1$. It follows that 
there is an induced injective homomorphism $A/\fp\hookrightarrow 
B/{\mathfrak q}$. (Consider the composite homomorphism $A\to B\to B/
{\mathfrak q}$. The kernel of this homomorphism is $A\cap{\mathfrak q}
=\fp$.) This injective homomorphism corresponds to a morphism of curves 
$\Spec(B/{\mathfrak q})\to \Spec(A/\fp)$.  By Hurwitz's theorem, the
genus of the former is at least as large as the genus of the latter. 
This proves the statement.

(2) We claim $B$ has only finitely many Poisson prime ideals of 
height 1. Let $S$ be the set of all Poisson prime ideals of 
height 1. Consider the Poisson ideal 
$I=\cap_{\mathfrak p\in S}\mathfrak p$. If $I\neq 0$, there are 
only finitely many primes minimal over $I$, each of which is 
Poisson prime of height one by Lemma \ref{xxlem1.7}. This implies 
that $S$ is finite. So it remains to show that $I\neq 0$. Suppose 
$I=0$. Let $\mathfrak p\in S$. Since $B/\mathfrak p$ is an affine 
Poisson domain of GK-dimension 1, Lemma \ref{xxlem1.6}(1) implies 
that $B/\mathfrak p$ has zero Poisson bracket. So 
$\{B,B\}\subseteq \mathfrak p$ and $\{B,B\}\subseteq I=0$. This 
contradicts the fact that $B$ has a nonzero Poisson bracket. So, our 
claim is proven. 

Therefore, after removing those finitely many $\fp$, we can assume 
that every $\fp$ is not a Poisson ideal. By construction, 
$F_1(B)=B\cap F_1(B_{\fp})=\fp$. Since $\fp$ is not a Poisson 
ideal, $\{B, \fp\} \not\subseteq \fp$. Therefore 
$\{F_0(B), F_1(B)\}\not\subseteq F_1(B)$. This means that 
${\mathbb F}^{ind}$ is not classical. Since $\GKdim (B/\fp)=1$,
the residue field of $\nu$ has GK-dimension 1, which implies
that $\nu$ is nondegenerate. It is clear that the image of $\nu$ 
is ${\mathbb Z}\cup\{\infty\}$. So $\nu$ is a faithful $1$-valuation.
\end{proof}

As a consequence, we have

\begin{proposition}
\label{xxpro5.6}
Let $B$ be an affine Poisson domain of GK-dimension 2 with 
nonzero Poisson bracket and let $K=Q(B)$. Then 
$\bfw_1(K)=\bfd_1(K)=\infty$.
\end{proposition}

\begin{proof}
Without loss of generality, we may assume $B$ is normal after localization and $K=Q(A[t^{\pm 1}])$ for some affine domain $A$ of GK-dimension 1. Let $N$ be 
an integer $>{\mathfrak g}(A)$. Applying Lemma \ref{xxlem5.5}(2) to $B$, there are infinitely many $(g,\fp)$ such that (i) $g\geq N$, 
(ii) ${\mathfrak g}(B/\fp)=g$, and (iii) ${\mathbb F}^{ind}$ 
in Lemma \ref{xxlem5.5}(2) is a faithful $1$-filtration. As a 
consequence, $\gr_{{\mathbb F}^{ind}} B_\fp \cong 
(B_\fp/\fp B_\fp)[t]\cong Q(B/\fp)[t]$. So, $Q(\gr_{{\mathbb F}^{ind}}(B))=Q((B/\fp)[t^{\pm 1}])\not\cong Q(A[t^{\pm 1}])$ 
by the cancellation theorem \cite[Theorem 2]{De}. Since there 
are infinitely many $g$, there are infinitely many pairwise
non-isomorphic Poisson fields $Q((B/\fp)
[t^{\pm 1}])$ with nonzero Poisson bracket. This implies that
there are infinitely many faithful $1$-valuations $\nu$ on $K$
such that $Q(\gr_{\nu}(K))=Q(\gr_{\nu}(B))=Q((B/\fp)
[t^{\pm 1}])$ are non-isomorphic. 
This means that $\bfw_1(K)=\infty$. 

Next we prove that $\bfd_1(K)=\infty$. Let $A_0:=A$ with genus 
$N_0$ and let $B_0=B$. By Lemma \ref{xxlem5.5}(2), there is a 
height one prime, but non-Poisson, ideal $\fp_0$ of $B_0$ such 
that the genus of $A_1:= B_0/\fp_0$ is at least one larger 
than $N_0$. Let $N_1$ be the genus of $A_1$. So $N_1>N_0$. By 
induction, for $i\geq 1$, let $B_i=\gr_{\nu_{i-1}} B_{i-1}$ where $\nu_{i-1}$ is the faithful $1$-valuation on $Q(B_{i-1})=Q(A_{i-1}[t^{\pm 1}])$ associated with $\fp_{i-1}$. We may still assume $B_i$ is normal after localization. We can define $A_{i+1}:=B_i/\fp_i$ 
such that the genus of $A_{i+1}$ is at least one larger than 
the genus of $N_{i}$. Let $N_{i+1}$ be the genus of $A_{i+1}$. 
Then $N_{i+1}>N_{i}$ by definition. By construction, we obtain 
a sequence of arrows 
$Q(A_i[t^{\pm 1}]) \to_{\nu_i} Q(A_{i+1}[t^{\pm 1}])$ with 
$A_i\not\cong A_j$ for all $i\neq j$. By \cite[Theorem 2]{De}, 
$Q(A_i[t^{\pm 1}])\not\cong Q(A_j[t^{\pm 1}])$ for all $i\neq j$. 
Thus $\bfd_1(K)=\infty$. 
\end{proof}

Next, we provide a structure theorem for certain $2$-valuations of some Poisson fields. Let $\Reg A$ denote the regular locus 
of $\Spec A$. If $A$ is affine, then $\Reg A$ is an open 
variety of $\Spec A$. If $A$ is a Poisson algebra, let $A_{cl}$ 
be the factor ring $A/(\{A,A\})$.

\begin{theorem}
\label{xxthm5.7}
Let $P$ be an affine Poisson domain of GK-dimension 2 and let $K$ be the 
Poisson fraction field $Q(P)$. Let $\fp$ be a maximal ideal
in $\Reg P$.
\begin{enumerate}
\item[(1)] 
There is a unique $2$-valuation, denoted by $\nu^{2,\fp}$, 
of $K$ such that $F_i(P)=\fp^i$ for all $i\geq 0$. 
\item[(2)]
$\nu^{2,\fp}$ is classical if and only of $\fp\in 
\Spec P_{cl}$.
\item[(3)] 
$\nu^{2,\fp}$ is either classical or Weyl, see 
Definition \ref{xxdef3.3}(3).
\end{enumerate}
\end{theorem}

In fact, in parts (1) and (2) of the above theorem, we
don't need to assume that $\GKdim P=2$. 

\begin{proof}[Proof of Theorem \ref{xxthm5.7}] 
(1) Let $\fp$ be a maximal ideal in $\Reg P$. Then the 
localization $B:=P_{\fp}$ is a regular local ring with 
maximal ideal, denoted by $\fm$. Define a filtration 
${\mathbb F}^{2,\fp}:=\{F_i\mid i\in {\mathbb Z}\}$ of
$B$ by
$$F^{2,\fp}_i=\begin{cases} \fm^i & i\geq 0.\\
B& i<0.\end{cases}$$
Since $B$ is regular local, $\gr_{{\mathbb F}^{2,\fp}}B$ 
is isomorphic to a polynomial $B_0[t_1,\cdots,t_d]$
where $B_0$ is the residue field and $d=\gldim B$
(when $\GKdim P=2$, then $d=2$). Hence 
${\mathbb F}^{2,\fp}$ is a good filtration. It is well-known in commutative
algebra that $\rk_{P/\fp} (\fp/\fp^2)=
\rk_{B/\fm} (\fm/\fm^2)=\gldim B$. When restricted to $P$, we have $F^{2,\fp}_0(P)=P$ and, for $i>0$, $F^{2,\fp}_i(P)
= P\cap \fm^i=\fp^i$ as required.

Next, we show that the valuation associated with 
${\mathbb F}^{2,\fp}$ is a $2$-valuation. Let 
$S:=\{s_1,\cdots,s_d\}$ be a regular system of parameters for the regular local ring $B=P_{\fp}$. Then 
$\fm=\sum_{\alpha=1}^d s_{\alpha} B$. Consequently, 
for all $i\geq 0$, $F^{2,\fp}_i=\sum_{\alpha_1,\cdots,\alpha_i}
s_{\alpha_1}\cdots s_{\alpha_i} B$. Using the fact that
$\{s_{\alpha}, s_{\beta}\}\in B$ and $\{B,B\}\subseteq B$, we 
obtain that $\{F^{2,\fp}_i,F^{2,\fp}_j\}\subseteq F^{2,\fp}_{i+j-2}$
for all $i,j$. Therefore ${\mathbb F}^{2,\fp}$ is a good 
$2$-filtration. We use $\nu^{2,\fp}$ for the associated 
$2$-valuation.

(2) Since $\gr_{\mathbb F} P$ is generated in degree 1, we have 
$$\begin{aligned}
{\text{$\nu^{2,\fp}$ is classical}}
&\Leftrightarrow
{\text{ $\{F_1(P),F_1(P)\}\subseteq F_{1+1-2+1}(P)$}}\\
&\Leftrightarrow
{\text{ $\{\fp,\fp\}\subseteq \fp$}}
\Leftrightarrow
{\text{ $\{\Bbbk+\fp,\Bbbk+\fp\}\subseteq \fp$}}\\
&\Leftrightarrow
{\text{ $\{P,P\}\subseteq \fp$}}
\Leftrightarrow
{\text{ $\fp\in \Spec P_{cl}$.}}
\end{aligned}
$$

(3) Since $P$ has GK-dimension two, $\gr_{\mathbb F} 
P$ is isomorphic to $\Bbbk[t_1,t_2]$ where $t_1$ and
$t_2$ are in degree 1. If $\nu^{2,\fp}$ is not classical, 
then $\{t_1,t_2\}\neq 0$. Since $\gr_{\mathbb F} P$
is Poisson $2$-graded, $\{t_1,t_2\}\in \Bbbk^{\times}$.
In this case $\nu^{2,\fp}$ is Weyl.
\end{proof}

By the proof of Theorem \ref{xxthm5.7}(1), there is a 
one-to-one correspondence between $\MaxReg P$ and the set of 
Poisson $2$-valuations $\{\nu^{2,\fp}\mid \fp\in \MaxReg P\}$.
If we can add an appropriate algebraic structure to ${\mathcal V}_{2}(K)$, then it makes sense to state that $\MaxReg P$ parameterizes the set of Poisson $2$-valuations $\{\nu^{2,\fp}\mid \fp\in \MaxReg P\}$.

By Lemma \ref{xxlem5.4}, if $\fp$ has height one, then 
$\nu^{2,\fp}$ is equivalent to a $1$-valuation. Therefore, it is a classical $2$-valuation. When $\GKdim P=2$, then 
$\nu^{2,\fp}$ is a faithful 2-valuation only when $\fp$ is a 
maximal ideal not in $\Spec P_{cl}$. In this case, Theorem \ref{xxthm5.7}(3) implies that $\nu^{2,\fp}$ is Weyl. Therefore, the Weyl Poisson field $K_{Weyl}$ is crucial in studying $2$-valuations.

\section{Embedding and isomorphism Problems}
\label{xxsec6}
In this section, we will discuss the embedding problem of determining whether one Poisson field can be embedded into another and the isomorphism problem of deciding whether two given Poisson fields are isomorphic for certain families of Poisson fields.

Firstly, let's make a simple observation about the embedding problem.

\begin{lemma}
\label{xxlem6.1}
Let $\Omega$ and $\Omega'$ be two i.s. potentials of 
degrees at least 4. Let $A_{\Omega}$, $P_{\Omega'}$, $P_{\Omega-\xi}$ and $P_{\Omega'-\xi}$, for $\xi\in \Bbbk^{\times}$ be defined
as in Construction \ref{xxcon0.6}.
\begin{enumerate}
\item[(1)]
$Q(A_{\Omega})$ cannot be embedded into $Q(P_{\Omega'})$
and $Q(P_{\Omega'-\xi})$.
\item[(2)]
$Q(P_{\Omega-\xi})$ cannot be embedded into $Q(P_{\Omega'})$.
\end{enumerate}
\end{lemma}

\begin{proof}
(1) This is clear since 
$$\GKdim Q(A_{\Omega})=3>2=\GKdim Q(P_{\Omega'})=
\GKdim Q(P_{\Omega'-\xi}).$$

(2) Suppose to the contrary that $Q(P_{\Omega-\xi})$ is a 
Poisson subfield of $Q(P_{\Omega'})$. Since $P_{\Omega'}$ 
is $w'$-graded Poisson domain where $w'=\deg \Omega'-3$, by 
Lemma \ref{xxlem3.2}(1), $Q(P_{\Omega'})$ has a nontrivial 
$(-w')$-valuation $\nu$. By Lemma \ref{xxlem1.4}(2), $\nu$ is 
also a nontrivial $0$-valuation on $Q(P_{\Omega'})$. By Lemma 
\ref{xxlem1.4}(1), $\nu$ restricts to a nontrivial 
$0$-valuation on $Q(P_{\Omega-\xi})$. This contradicts Lemma 
\ref{xxlem4.8}(3). 
\end{proof}

It is not immediately obvious to us whether $Q(P_{\Omega})$ 
(resp. $Q(P_{\Omega-\xi})$) can be embedded into $Q(A_{\Omega'})$.

We will make more statements about embedding and isomorphism 
problems that help us to distinguish Poisson fields.
Let $m$ be an integer at least $2$. Let $T(m,p_{ij})$ denote a 
Poisson torus $\Bbbk[x_1^{\pm 1},\cdots,x_m^{\pm 1}]$ with 
Poisson bracket determined by $\{x_i,x_j\}=p_{ij} x_i x_j$ 
for $p_{ij}\in \Bbbk$ for all $1\leq i<j\leq m$.

\begin{theorem}
\label{xxthm6.2}
Let $T$ be the Poisson torus $T(m,p_{ij})$ defined above and 
assume that $T$ is Poisson simple. Let $K$ be the Poisson 
fraction field of $T$.
\begin{enumerate}
\item[(1)]
There is a one-to-one correspondence between the set 
of $0$-valuations and the set ${\mathbb Z}^m$. Further, for every $0$-valuation $\nu$, $Q(\gr_{\nu}(K))\cong K$. 
As a consequence, $\alpha_0(K)=\infty$, $\bfd_0(K)=0$, 
$\bfw_0(K)=1$, and $K$ is $0$-quasi-Adams.
\item[(2)]
There is a one-to-one correspondence between the set 
${\mathcal V}_{f0}(K)$ and the set 
$\{(v_1,\cdots,v_m)\in {\mathbb Z}^m \mid \gcd(v_1,\cdots,v_n)=1\}$. 
\item[(3)]
There is no $w$-valuation for all $w<0$.
\end{enumerate}
\end{theorem}

\begin{proof}
(1,2) First we prove that there is a one-to-one correspondence between the set of $0$-valuations and the set ${\mathbb Z}^m$. 
Let $\nu$ be a $0$-valuation on $K$. Let $v_i=\nu(x_i)$. Then 
$(v_1,\cdots,v_m)\in {\mathbb Z}^m$. Conversely, given 
$(v_1,\cdots,v_m)\in {\mathbb Z}^m$ we claim that there is a 
unique $0$-valuation $\nu$ such that $\nu(x_i)=v_i$ for all $i$. 
Let ${\mathbb F}^{ind}$ be the induced filtration determined by \eqref{E2.8.3} with the degree assignment $\deg(x_i)=v_i$ 
and $\deg(x_i^{-1})=-v_i$ for all $i$. By the definition of 
$T$, it is easy to see that $T$ is in fact ${\mathbb Z}$-graded 
with degree assignment $\deg(x_i^{\pm 1})=\pm v_i$ for all $i$. 
By Lemma \ref{xxlem3.2}(1), ${\mathbb F}^{ind}$ agrees with the 
Adams$^{Id}$ filtration ${\mathbb F}^{Id}$. As a consequence,
${\mathbb F}^{ind}$ is a good filtration which provides a 
$0$-valuation $\nu$ such that $\nu(x_i)=v_i$ for all $i$. It 
remains to show that such $\nu$ is unique. 

Now let $\nu$ be (another) $0$-valuation on $K$ such that 
$\nu(x_i)=v_i$ for all $i$ and let ${\mathbb F}$ be the 
associated filtration. By the valuation axioms, 
$\nu(x_i^{-1})=-v_i$ for all $i$. By the argument in the 
previous paragraph and Lemma \ref{xxlem2.9}(2), we have a 
sequence of algebra homomorphisms
$$T\cong \gr_{{\mathbb F}^{Id}} T\cong
\gr_{{\mathbb F}^{ind}} T\to I\subseteq \gr_{\mathbb F} T.$$
It is clear that $T$ is Poisson $0$-graded. By Lemma 
\ref{xxlem2.12}, $\phi^{ind}: \gr_{{\mathbb F}^{ind}} T \to 
\gr_{\mathbb F} T$ is a Poisson $0$-graded algebra homomorphism.
Since $T$ is Poisson simple, $\phi^{ind}$ is injective. Thus 
$\GKdim I\geq \GKdim T=n$. By Lemma \ref{xxlem2.11}(1),
${\mathbb F}$ agrees with ${\mathbb F}^{ind}(={\mathbb F}^{Id})$. 
Therefore $\nu$ agrees with $\nu^{Id}$. This proves the 
uniqueness. 

Note that ${\nu}^{Id}$ is a faithful $0$-valuation if and only 
if $\gcd\{\nu(x_1),\cdots,\nu(x_n)\}=1$. Assertion (2) 
follows. It is clear that there are infinitely many 
$\{c_1,\cdots,c_n\}\in {\mathbb Z}^n$ such that
$\gcd(c_1,\cdots,c_n)=1$. Hence $\alpha_0(K)=\infty$.
By the above proof, $\gr_{\mathbb F} T\cong T$ for any 
$0$-valuation $\nu$. By Lemma \ref{xxlem2.7}(2),
$Q(\gr_{\mathbb F} K)\cong K$. Or equivalently, $K\to_{\nu} K$
for all $0$-valuations $\nu$. By definition, we have 
$\bfd_0(K)=0$, $\bfw_0(K)=1$, and $K$ is $0$-quasi-Adams.

(3) By the proof of part (1,2), every $0$-valuation is 
nonclassical. The assertion follows from Lemma \ref{xxlem1.4}(2).
\end{proof}

\begin{example}
\label{xxexa6.3}
For every $q\in \Bbbk^{\times}$, the $q$-skew Poisson field,
denoted by $K_q$, is defined to be the $Q(T(2, p_{12}=q))$. 
Similar to Example \ref{xxexa2.10}, we can write
$$K_{q}=\Bbbk(x,y\mid \{x,y\}=q xy).$$ 
\end{example}
It is worth noting that Poisson valuations cannot separate non-isomorphic $q$-skew Poisson fields such as $K_{1}$ and $K_{2}$. One of the key results regarding the isomorphism problem is presented below.

\begin{corollary}
\label{xxcor6.4} 
The following Poisson fields are nonisomorphic.
\begin{enumerate}
\item[(1)]
The $q$-skew Poisson field $K_q$;
\item[(2)]
The Weyl Poisson field $K_{Weyl}$;
\item[(3)]
$Q_1:=Q(P_{\Omega})$ where $\Omega=x^3+y^3+z^3+\lambda xyz$
where $\lambda^3\neq -3^3$;
\item[(4)]
$Q_2:=Q(P_{\Omega-1})$ where $\Omega=x^3+y^3+z^3+\lambda xyz$
where $\lambda^3\neq -3^3$;
\item[(5)]
$Q(P_{\Omega})$ where $\Omega$ is an i.s. potential 
of degree $\geq 4$;
\item[(6)]
$Q(P_{\Omega-\xi})$ where $\xi\in \Bbbk^{\times}$ and 
$\Omega$ is an i.s. potential of degree $\geq 4$.
\end{enumerate}
\end{corollary}

\begin{proof} We have 
\begin{enumerate}
\item[(1)]
$\alpha_{-1}(K_q)=0$ [Theorem \ref{xxthm6.2}(3)]
and $\alpha_0(K_q)=\infty$ [Theorem \ref{xxthm6.2}(1)].
\item[(2)]
$\alpha_{-1}(K_{Weyl})\neq 0$ [Example \ref{xxexa2.10}].
\item[(3)]
$\alpha_{-1}(Q_1)=0$ and $\alpha_0(Q_1)=2$ [Theorem \ref{xxthm3.8}(1)].
\item[(4)]
$\alpha_{-1}(Q_2)=0$ and $\alpha_0(Q_2)=1$ [Theorem \ref{xxthm3.11}(1)].
\end{enumerate}
Hence, these four Poisson fields are pair-wisely nonisomorphic.

In cases (1-4), there is a nontrivial $0$-valuation. By Lemma 
\ref{xxlem4.8}(3), (6) is not isomorphic to (1-4). By 
Lemma \ref{xxlem6.1}(2), (5) and (6) are nonisomorphic. 
Finally, (5) is not isomorphic to (1-4) since the Poisson 
fields in (5) have no faithful $0$-valuation by \Cref{xxlem4.8}(5), but each Poisson
field in (1-4) has at least one faithful $0$-valuation.
\end{proof}

We remark that these fields in (1), (2) and (4) of the above corollary are isomorphic to $\kk(x,y)$, but the fields in other items are not isomophic to each other. The fields in (3) and (5) are not isomorphic due to the Zariski cancellation. The field in (6) is not isomorphic to that of items (1)-(5) due to the fact that the arithmetic genus of a smooth hypersurface  is a birational invariant.

In the next section, we will continue working on the isomorphism problem within classes (5) and (6).
Note that Corollary \ref{xxcor6.4} and its proof are closely 
related to \cite[Proposition 5.3]{Ar}. 

For the rest of this section, we provide further detailed 
information about Poisson fields listed in Corollary 
\ref{xxcor6.4}(1-4). Hopefully, these results are also helpful in understanding other aspects of Poisson fields in Corollary 
\ref{xxcor6.4}(1-4).

\begin{lemma}
\label{xxlem6.5}
Let $K$ be a Poisson field containing $K_{q}$ where 
$q\in \Bbbk^{\times}$.
\begin{enumerate}
\item[(1)]
Every $0$-valuation on $K$ is nonclassical. As a consequence, 
${\mathcal V}_{-1}(K)=\emptyset$ and $ncw(K)\geq w(K)\geq 0$.
\item[(2)]
$K_q$ cannot be embedded into $Q(P_{\Omega})$ for $\deg \Omega
\neq |x|+|y|+|z|$.
\item[(3)]
For every $0$-valuation $\nu$ on $K$, $Q(\gr_{\nu}(K))$ 
contains $K_{q}$ as a Poisson subfield.
\end{enumerate}
\end{lemma}

\begin{proof}
(1) Let $\nu$ be a $0$-valuation on $K$ and $\mu$ be the
restriction of $\nu$ on $K_{q}$. Then $\gr_{\mu} K_{q}
\subseteq \gr_{\nu} K$ as Poisson algebras. By Theorem 
\ref{xxthm6.2}(1), $Q(\gr_{\mu} K_{q})\cong K_{q}$. So $\mu$ 
is nonclassical. As a consequence, $\nu$ is nonclassical. 
Therefore ${\mathcal V}_0(K)={\mathcal V}_{nc0}(K)$. By Lemma 
\ref{xxlem1.4}(2), ${\mathcal V}_{-1}(K)=\emptyset$. Hence 
$ncw(K)\geq w(K)\geq 0$.

(2) This follows from part (1) and the fact that 
${\mathcal V}_{-1}(Q(P_{\Omega}))\neq \emptyset$ when 
$\deg \Omega \neq |x|+|y|+|z|$.

(3) By the proof of part (1), $K_{q}\cong Q(\gr_{\mu} K_{q})
\subseteq Q(\gr_{\nu} K)$ as Poisson algebras. 
\end{proof}

Recall from Definition \ref{xxdef3.3}(3) that a valuation 
$\nu$ is called {\it Weyl} if $Q(\gr_{\nu}(K))$ is 
isomorphic to the Weyl Poisson field $K_{Weyl}$. Next we 
describe $1$-valuations and $2$-valuations on $K_q$. The 
following lemma is clear.

\begin{lemma}
\label{xxthm6.6}
Let $T$ be $T(2, p_{12}=q)$ where $q\in \Bbbk^\times$ 
generated by $x:=x_1,y:=x_2$. Let $\nu$ be any valuation 
on $T$. Then, after a base change 
$$x\mapsto x^a y^b, \quad y\mapsto x^cy^d$$
for some $\begin{pmatrix} a& b\\ c& d\end{pmatrix}
\in SL_2({\mathbb Z})$, we may assume that $\nu(x)=0$.
\end{lemma}

By definition, $K_q=Q(T(2, p_{12}=q))$. 

\begin{proposition}
\label{xxpro6.7}
Assume that the base field $\Bbbk$ is algebraically closed.
Let $K_q$ be the $q$-skew Poisson field where $q\neq 0$. 
\begin{enumerate}
\item[(1)]
Every nontrivial nonclassical $1$-valuation on $K_q$ is either 
Weyl or $\nu^{1,\fp}$ for some height one prime $\fp\in 
\Spec T$ with $T=T(2, p_{12}=q)$ as in Theorem {\rm{\ref{xxthm4.3}(1)}}.
\item[(2)]
Every faithful $2$-valuation on $K_q$ is Weyl.
\end{enumerate}
\end{proposition}

\begin{proof} Let $T$ be $T(2, p_{12}=q)$ and $P$ be the Poisson 
polynomial ring $\Bbbk[x,y]$ with $\{x,y\}=qxy$. Then 
$P\subseteq T\subseteq K_q$ and $K_q=Q(P)=Q(T)$. 

(1) Let $\nu$ be a nontrivial nonclassical $1$-valuation on 
$K_q$ with associated filtration ${\mathbb F}$. We consider the 
following cases.

Case 1: $F_0(T)=T$. If $\GKdim F_0(T)/F_1(T)=0$, 
then $T/F_1(T)=\Bbbk$. So $u:=x-c_1$ and $v:=y-c_2$ have positive $\nu$-values for
some $c_1,c_2\in \Bbbk^{\times}$.
Then
$$0=\nu(qxy)=\nu(\{x,y\})=\nu(\{u,v\})\geq \nu(u)+\nu(v)-1,$$
yielding a contradiction. If $\GKdim F_0(T)/F_1(T)=1$, then 
$\fp:=F_1(T)$ is a prime ideal of $T$ of height one. Let 
$B$ be the localization $T_{\fp}$ which is a DVR. By 
Theorem \ref{xxthm4.3}(1), $\nu=\nu^{1,\fp}$. If 
$\GKdim F_0(T)/F_1(T)=2$, then $F_1(T)=0$ and whence 
$\nu$ is trivial.  In this case, we only obtain valuations
of the form $\nu^{1,\fp}$.

Case 2: $F_0(T)\neq T$. By Lemma \ref{xxthm6.6}, we may 
assume that $\nu(x)=0$, and consequently, $a:=\nu(y)\neq 0$
(if $a= 0$, then $F_0(T)=T$). We may assume $a<0$ (after 
replacing $y$ by $y^{-1}$ if necessary).
Let $I$ be the subalgebra of $\gr_{\mathbb F} P$ generated 
by $\overline{x}$ and $\overline{y}$. If $\GKdim I=2$, then 
Lemma \ref{xxlem2.11}(1), ${\mathbb F}={\mathbb F}^{ind}$.
In this case, $\nu$ is classical by Lemma \ref{xxlem2.9}(3). 
If $\GKdim I=1$, then $I\subseteq \Bbbk[s]$ with $\deg s<0$
[Lemma \ref{xxlem1.3}(3)], and whence, $u:=x-c$ has 
positive $\nu$-value for some $c\in \Bbbk^{\times}$. Then 
$$a=\nu(qxy)=\nu(\{x,y\})=\nu(\{u,y\})\geq \nu(u)+\nu(y)-1=
\nu(u)+a-1.$$
This implies that $\nu(u)=1$. In $(\gr_{{\mathbb F}} T)_{a}$, 
we have
\begin{equation}
\label{E6.7.1}\tag{E6.7.1}
\{\overline{u},\overline{y}\}=\overline{\{u,y\}}
=\overline{\{x,y\}} =\overline{qxy}=cq\overline{y} \neq 0.
\end{equation}
Let $\widetilde{I}$ be the subalgebra of $\gr_{{\mathbb F}} P$ 
generated by $\overline{u},\overline{y}$ (note that 
$P$ is the subalgebra of $K_q$ generated by $u$ and $y$). Then, by Lemma \ref{xxlem1.6}(2), $\GKdim \widetilde{I}=2$.  Then, by 
Lemma \ref{xxlem2.11}(1), ${\mathbb F}(P)={\mathbb F}^{ind}(P)$. 
By \eqref{E6.7.1}, $Q(\gr_{{\mathbb F}^{ind}}(P))$ is 
$K_{Weyl}$. Hence, $\nu$ is Weyl. In this case, we only obtain 
Weyl valuations. 

The assertion follows by combining the above two cases.

(2) Let $\nu$ be a faithful $2$-valuation on $K_q$ with associated 
filtration ${\mathbb F}$. We consider the following cases.

Case 1: $F_0(T)=T$. If $\GKdim F_0(T)/F_1(T)=0$, then 
$F_0(T)/F_1(T)=\Bbbk$. Thus $u:=x-c_1$ and $v:=y-c_2$ possess positive $\nu$-values for
some $c_1,c_2\in \Bbbk^{\times}$.
Then
$$0=\nu(qxy)=\nu(\{x,y\})=\nu(\{u,v\})\geq \nu(u)+\nu(v)-2.$$
This implies that $\nu(u)=\nu(v)=1$. Further, $\{\overline{u},
\overline{v}\}=qc_1c_2\neq 0$ in $(\gr_{\nu} T)_{0}$. Let $I$ 
be the subalgebra of $\gr_{\mathbb F} P$ generated by 
$\overline{u}$ and $\overline{v}$. By Lemma \ref{xxlem1.6}(1), 
$\GKdim I=2$. It follows from Lemma \ref{xxlem2.11}(1), 
${\mathbb F}={\mathbb F}^{ind}$ where ${\mathbb F}^{ind}$
is the induced filtration of $P$ determined by $\deg(u)=\deg(v)=1$
(which is a good filtration). It is clear that 
$Q(\gr_{{\mathbb F}^{ind}} P)\cong K_{Weyl}$. So $\nu$ is Weyl.
If $\GKdim F_0(T)/F_1(T)=1$, then $\fp:=F_1(T)$ is a prime
ideal of $T$ of height one. Let $B$ be the localization $T_{\fp}$. Then $B$ is a DVR. By Theorem \ref{xxthm4.3}(1), 
$\nu=\nu^{1,\fp}$. This is classical considered as a $2$-valuation. 
If $\GKdim F_0(T)/F_1(T)=2$, then $F_1(P)=0$ and whence $\nu$ is 
trivial. 

Case 2: $F_0(T)\neq T$. By Lemma \ref{xxthm6.6}, we may assume 
that $\nu(x)=0$, and by the proof of part (1) Case 1, $a:=\nu(y)<0$. 
Let $I$ be the subalgebra of $\gr_{\mathbb F} P$ generated by 
$\overline{x}$ and $\overline{y}$. If $\GKdim I=2$, then Lemma 
\ref{xxlem2.11}(1), ${\mathbb F}={\mathbb F}^{ind}$. In this case, 
$\nu$ is classical by Lemma \ref{xxlem2.9}(3), yielding a 
contradiction. 

If $\GKdim I=1$, then $u:=x-c$ has positive $\nu$-value for 
some $c\in \Bbbk^{\times}$. Then 
$$a=\nu(qxy)=\nu(\{x,y\})=\nu(\{u,y\})\geq \nu(u)+\nu(y)-2=
\nu(u)+a-2.$$
This implies that either $\nu(u)=1$ or $\nu(u)=2$. If $\nu(u)=2$, 
then $\{\overline{u},\overline{y}\}=qc \overline{y} \neq 0$ in 
$(\gr_{\mathbb F} P)_a$. An argument similar to Case 2 of the 
proof of part (1) shows that $\nu$ is Weyl. Next, we consider the 
case when $\nu(u)=1$. Let $\widetilde{I}$ be the subalgebra 
generated by $\overline{u}$ and $\overline{y}$. If 
$\GKdim \widetilde{I}=2$, similar to the argument at the 
beginning of the proof of Case 2, we have that $\nu$ is 
classical. If $\GKdim \widetilde{I}=1$, then 
$\overline{y}=c' \overline{u}^a$ for some $c'\in \Bbbk^{\times}$
[Lemma \ref{xxlem1.3}(2)]. Let $v:=y-c'u^a$. Then $\nu(v)\geq \nu(y)$. Now 
$$a=\nu(qxy)=\nu(\{x,y\})=\nu(\{u,y\})=\nu(\{u, v\})
\geq \nu(u)+\nu(v)-2=\nu(v)-1.$$
This implies that $\nu(v)=a+1$ and $\{\overline{u},\overline{v}\}
=qcc' \overline{u}^a\neq 0$ in $(\gr_{\mathbb F} K_q)_a$.
By Lemma \ref{xxlem1.6}(1), the subalgebra generated by
$\overline{u}$ and $\overline{v}$, denoted by $I'$, has GK-dimension
$2$. Let $A$ be the subalgebra of $K_q$ generated by $u,u^{-1},v$.
It is clear that $Q(A)=K_q$. Let ${\mathbb F}^{ind}$ be the 
induced filtration determined by $\deg(u^{\pm 1})=\pm 1$
and $\deg(v)=a+1$. By Lemma \ref{xxlem2.11}(1), 
${\mathbb F}={\mathbb F}^{ind}$. In this case, $\nu$ is Weyl since
$\{\overline{u},\overline{v}/\overline{u}^a\}=qcc'$.

The assertion follows by combining the above two cases.
\end{proof}

Proposition \ref{xxpro6.7} provides useful information about 
the $1$- and $2$-valuations on $K_q$. Such a
detailed description can rarely be obtained for any Poisson field. 
In general, it is always challenging to work out the 
complete set of $1$- and $2$-valuations. 

The Sklyanin Poisson field $Skly_3$ is defined to be 
$Q(P_{\Omega-1})$ where $\Omega=x^3+y^3+z^3+\lambda xyz$
where $\lambda^3\neq -3^3$ as in Corollary \ref{xxcor6.4}. 

\begin{corollary}
\label{xxcor6.8}
Let $Skly_3$ be the Sklyanin Poisson field and $K$ be a
Poisson field of GK-dimension two. 
\begin{enumerate}
\item[(1)]
If ${\mathcal V}_{-1}(K)\neq \emptyset$, then there is no 
Poisson algebra homomorphism from $Skly_3$ to $K$.
\item[(2)]
If $K$ is either the Weyl Poisson field or $Q(P_{\Omega})$
where $|\Omega|\neq |x|+|y|+|z|$, then there is  no 
Poisson algebra homomorphism from $Skly_3$ to $K$.
\end{enumerate}
\end{corollary}

\begin{proof}
(1) Suppose there is an embedding $Skly_3\hookrightarrow K$. Let $\nu$ be an 
element in ${\mathcal V}_{-1}(K)$ and let $\mu$ be the 
restriction of $\nu$ on $Skly_3$. Then $\nu\in 
{\mathcal V}_{-1}(Skly_3)$, which contradicts Theorem \ref{xxthm3.8}(1).

(2) The assertion follows because ${\mathcal V}_{-1}(K)\neq \emptyset$
for all Poisson fields $K$ in part (2).
\end{proof}

\begin{corollary}
\label{xxcor6.9}
Suppose $\Omega$ is an i.s. potential of degree $\geq 4$.
Let $Q$ be $Q(P_{\Omega-\xi})$ where $\xi\in \Bbbk^{\times}$ 
and $K$ be a Poisson field of GK-dimension two. 
\begin{enumerate}
\item[(1)]
If ${\mathcal V}_{0}(K)\neq \emptyset$, then there is no 
Poisson algebra homomorphism from $Q$ to $K$.
\item[(2)]
If $K$ is either $K_{Weyl}$, or $K_q$, or $Q(P_{\Omega'})$ where 
$\Omega'$ is any homogeneous potential, so there is no Poisson algebra homomorphism from $Q$ to $K$.
\end{enumerate}
\end{corollary}

\begin{proof}
(1) Suppose there is an embedding $Q\hookrightarrow K$. Let $\nu$ be an 
element in ${\mathcal V}_{0}(K)$ and let $\mu$ be the 
restriction of $\nu$ on $Q$. By Lemma \ref{xxlem1.4}(1),
$\mu$ is nontrivial. This yields a contradiction 
to Theorem \ref{xxlem4.8}(3).

(2) In all cases, one has that ${\mathcal V}_{0}(K)\neq 
\emptyset$. The assertion follows from part (1).
\end{proof}

Without using valuations, it probably isn't easy to demonstrate any of the non-embedding results above. 

\section{$\epsilon$-morphisms}
\label{xxsec7}
To solve embedding and isomorphism problems within 
the class $Q(P_{\Omega-\xi})$ (as well as the automorphism problem in the next section), it is convenient to consider 
$\epsilon$-morphisms to be introduced in this section.

Let $A$ be a Poisson algebra with a Poisson bracket 
$\{-,-\}$ and let $e$ be in $\Bbbk^{\times}$. Then
we can define a new Poisson structure $\{-,-\}_e$ on $A$ by
$$\{f,g\}_e:= e\{f,g\} \quad {\text{for all $f,g\in A$.}}$$
This new Poisson algebra is denoted by $A(e)$. 

\begin{definition}
\label{xxdef7.1}
Let $A$ and $B$ be two Poisson algebras. A $\Bbbk$-linear 
map $\phi: A\to B$ is called an {\it $e$-morphism}
if it is a Poisson algebra homomorphism from $A\to B(e)$ for
some $e\in \Bbbk^{\times}$. If $e$ is not specified, we say $\phi$ is an {\it $\epsilon$-morphism}.

Similarly, one can define $e$-($\epsilon$)-versions of Poisson
endomorphism, Poisson isomorphism, Poisson automorphism, etc.
\end{definition}

The following lemma is easy.

\begin{lemma}
\label{xxlem7.2}
Let $A$ be a Poisson algebra and $e\in \Bbbk^{\times}$.
The following hold.
\begin{enumerate}
\item[(1)]
$\nu$ is a $w$-valuation on $A$ if and only if 
it is a $w$-valuation on $A(e)$.
\item[(2)]
$\Aut_{Poi}(A)=\Aut_{Poi}(A(e))$.
\end{enumerate}
\end{lemma}

By the above lemma, we can replace Poisson homomorphisms 
by $\epsilon$-morphisms when $w$-valuations are 
concerned.

\begin{definition}
\label{xxdef7.3}
Let $A$ be a Poisson algebra with a nontrivial Poisson 
bracket. If $\phi$ is an $\epsilon$-automorphism of $A$,
its {\it Poisson determinant} of $\phi$ is defined to be the
scalar $e$ such that $\phi$ is an $e$-automorphism of $A$. 
We write $\Pdet(\phi)=e$ in this case.
\end{definition}

Homological determinant of a graded algebra automorphism 
$\sigma$ of an Artin-Schelter regular algebra $A$ was introduced
in \cite{JZ}. In the study of Artin-Schelter regular algebras and
superpotentials, Smith and Mori gave a nice interpretation of the 
homological determinant in terms of the following 
formula \cite[Theorem 1.2]{MS}
$$\sigma^{\otimes l} (w)=\hdet(\sigma) w$$
where $w\in (A_1)^{\otimes l}$ is the superpotential associated with an $m$-Koszul Artin-Schelter regular algebra $A=\bigoplus_{i\ge 0}A_i$ of AS-index $l$. Motivated by this, we 
propose the following definition.

\begin{definition}
\label{xxdef7.4}
Let $A_{\Omega}=\kk[x,y,z]$ be the Poisson polynomial algebra given in Construction \ref{xxcon0.6}.
Suppose that $\Omega$ is an i.s. potential. 
Let $\phi$ be an $\epsilon$-automorphism of $A_{\Omega}$.
The {\it homological determinant} of $\phi$ is defined to 
be the scalar $c\in \kk^\times$ such that $\phi(\Omega)=c\,\Omega$.
\end{definition}

For three polynomials $f,g,h$ in $\Bbbk[x,y,z]$ (with a 
given graded generators $(x,y,z)$), the Jacobian determinant
of $(f,g,h)$ is defined to be
\begin{equation}
\label{E7.4.1}\tag{E7.4.1}
J(f,g,h):=\det \begin{pmatrix}
f_x &f_y & f_z\\ g_x& g_y& g_z\\h_x&h_y&h_z
\end{pmatrix}
=:{\rm det}\left(\frac{\partial(f,g,h)}{\partial(x,y,z)}\right)
\end{equation}
which was used in Construction \ref{xxcon0.6}.

As noted at the beginning of this section, we are trying to 
understand embedding and isomorphism problems within the classes in (5) and (6) in Corollary \ref{xxcor6.4}.
For the rest of this section, we assume that $\Omega$ and
$\Omega'$ are two i.s. potentials of degrees at least five unless stated otherwise. 
Let $\xi,\xi'\in \Bbbk$, which could be zero.

\begin{lemma}
\label{xxlem7.5}
Let $\phi$ be an $\epsilon$-morphism from $Q(P_{\Omega-\xi})$ 
to $Q(P_{\Omega'-\xi'})$. 
\begin{enumerate}
\item[(1)]
Then $\phi$ restricts to an 
injective $\epsilon$-morphism $P_{\Omega-\xi}\to 
P_{\Omega'-\xi'}$.
\item[(2)]
$\deg \Omega\leq \deg \Omega'$.
\end{enumerate}
\end{lemma}

\begin{proof} (1) By Theorem \ref{xxthm4.9}, 
$\,^1\Gamma(Q(P_{\Omega-\xi}))=P_{\Omega-\xi}$ and 
$\,^1\Gamma(Q(P_{\Omega'-\xi'}))=P_{\Omega'-\xi'}$. 
The assertion follows from Lemma \ref{xxlem4.2}(1). 

(2) Suppose to the contrary that 
$n+3:=\deg \Omega> \deg \Omega'=: n'+3$. Let $w=n'<n$.
By Lemmas \ref{xxlem4.2}(1), Theorem \ref{xxthm4.9} and Theorem \ref{xxthm4.10},
$$P_{\Omega-\xi}=^{w}\Gamma(Q(P_{\Omega-\xi}))
\subseteq ^{w}\Gamma(Q(P_{\Omega'-\xi'}))=\Bbbk,$$
yielding a contradiction.
\end{proof}

Next, we consider the case when $\deg \Omega=\deg \Omega'$. 
Let $\phi: P_{\Omega-\xi}\to P_{\Omega'-\xi'}$ be an 
algebra homomorphism. We say $\phi$ is {\it linear} if
$\phi(x), \phi(y), \phi(z)\in \Bbbk+\Bbbk x+\Bbbk y+ \Bbbk z$.
When $\xi=\xi'=0$, we say $\phi$ is {\it graded} if 
$\phi(x), \phi(y), \phi(z)\in \Bbbk x+\Bbbk y+ \Bbbk z$.

\begin{lemma}
\label{xxlem7.6}
Suppose $\deg \Omega=\deg \Omega'\geq 4$. Let $\phi$ be an 
injective $\epsilon$-morphism from $P_{\Omega-\xi}\to 
P_{\Omega'-\xi'}$.
\begin{enumerate}
\item[(1)]
$\phi$ is linear. Furthermore,
$\phi(x),\phi(y),\phi(z)\in \Bbbk x+\Bbbk y+\Bbbk z$. 
\item[(2)]
$\phi$ is bijective.
\item[(3)]
$\xi=0$ if and only if $\xi'=0$.
\item[(4)]
If $\xi=\xi'=0$, then $\phi$ is graded. 
\end{enumerate}
\end{lemma}

\begin{proof}
 We use $x,y,z$ for the generators of both $P_{\Omega-\xi}$ and 
$P_{\Omega'-\xi'}$. Let $w=\deg \Omega-3=\deg \Omega'-3\geq 1$. We employ the notations
introduced in Lemma \ref{xxlem4.8}. For example, $d$ and $d'$ 
 denote the minimum value of a valuation on the generators $x,y,z$ of $(P_{\Omega-\xi})_1$ and $(P_{\Omega'-\xi'})_1$, respectively.

(1,2) By Lemma \ref{xxlem4.8}(4), there is a unique $w$-valuation 
on $P_{\Omega'-\xi'}$ with $d'=-1$, denoted by $\nu'$. In this 
case, $\nu'(x)=\nu'(y)=\nu'(z)=-1$ when applied to $x,y,z\in 
Q(P_{\Omega'-\xi'})$ and $F^{\nu'}_0(P_{\Omega'-\xi'})=\Bbbk$. 

Let $\nu$ be the pullback valuation on $P_{\Omega-\xi}$ via 
$\phi$. Since $P_{\Omega-\xi}$ is a subalgebra of 
$P_{\Omega'-\xi'}$ via $\phi$, we have that 
$F^{\nu}_0(P_{\Omega-\xi})=\Bbbk$. So $d< 0$. As a consequence
of Lemma \ref{xxlem4.8}(4), $d=-1$, whence $\nu$ is the unique 
valuation on $P_{\Omega-\xi}$ whose filtration is either 
${\mathbb F}^{-Id}$ or ${\mathbb F}^c$. As a consequence, 
$\nu(x)=\nu(y)=\nu(z)=-1$. This means that 
$\phi$ maps $\Bbbk+\Bbbk x+\Bbbk y+\Bbbk z$ in $P_{\Omega-\xi}$
to $\Bbbk+\Bbbk x+\Bbbk y+\Bbbk z$ in $P_{\Omega'-\xi'}$.
Thus, $\phi$ is surjective. So $\phi$ is bijective and linear.

Now assume that
$$\begin{aligned}
\phi(x)&=x_1+a,\\
\phi(y)&=y_1+b,\\
\phi(z)&=z_1+c,\\
\end{aligned}
$$
where $x_1,y_1,z_1\in \Bbbk x+\Bbbk y+\Bbbk z$ and 
$a,b,c\in \Bbbk$. We compute
$$\begin{aligned}
0&= \phi(\Omega-\xi)=\Omega(\phi(x),\phi(y),\phi(z))-\xi\\
&= \Omega(x_1+a,y_1+b,z_0+c)-\xi\\
&=\Omega(a,b,c)-\xi+\Omega_x(a,b,c) x_1+\Omega_{y}(a,b,c)
y_1+\Omega_z(a,b,c) z_1+hdt \in P_{\Omega'-\xi'}
\end{aligned}
$$
where $hdt$ is a linear combination of higher Adams degree terms. Hence $\Omega_x(a,b,c)=\Omega_y(a,b,c)
=\Omega_z(a,b,c)=0$. As $\Omega$ has an isolated singularity 
at the origin, $a=b=c=0$ as desired.

(3) Note that $P_{\Omega-\xi}$ has finite global
dimension if and only if $\xi\neq 0$. Since $\phi$ is an algebra isomorphism from $P_{\Omega-\xi}$ to $P_{\Omega'-\xi'}$, the assertion follows.

(4) When $\xi=\xi'=0$, then $P_{\Omega-\xi}$ and $P_{\Omega'-\xi'}$ 
are connected Adams graded. The assertion follows from part (1).
\end{proof}

Combining Lemma \ref{xxlem7.5}(1) with Lemma \ref{xxlem7.6},
we deduce 

\begin{corollary}
\label{xxcor7.7}
Suppose $\deg \Omega=\deg \Omega'\geq 5$.
Let $\phi$ be a $\epsilon$-morphism from $Q(P_{\Omega-\xi})$ 
to $Q(P_{\Omega'-\xi'})$. Let $\phi$ also denote the restriction
$P_{\Omega-\xi}\to P_{\Omega'-\xi'}$ as in Lemma \ref{xxlem7.5}(1).
Then (1)-(4) of Lemma \ref{xxlem7.6} hold.
\end{corollary}

If $\phi$ is an automorphism of the commutative algebra
$\Bbbk[x,y,z]$. We write $J(\phi):=J(\phi(x),\phi(y),\phi(z))$
which is $\det\left( \frac{\partial(\phi(x),\phi(y),\phi(z))}{\partial(x,y,z)}\right)$ by definition. 

\begin{lemma}
\label{xxlem7.8}
Let $\phi: A_{\Omega}\to A_{\Omega'}$ be an algebra 
isomorphism without considering their Poisson structures.
\begin{enumerate}
\item[(1)]
Then $\phi$ is an $e$-isomorphism if and only if $\phi(\Omega)=
a\Omega'$ where $a=eJ(\phi)\in \Bbbk^{\times}$.
\item[(2)]
Then $\phi$ is a Poisson algebra isomorphism if and only if 
$\phi(\Omega)=J(\phi)\,\Omega'$.
\item[(3)]
Then $\phi$ is a Poisson algebra automorphism of $A_{\Omega}$ if and only if 
$\phi(\Omega)=J(\phi)\,\Omega$.
\item[(4)] If $\phi$ is a Poisson $\epsilon$-automorphism of $A_{\Omega}$, then $\Pdet(\phi)=aJ(\phi)^{-1}$ where $\phi(\Omega)=a\Omega$ for some $a\in \kk^\times$.
\item[(5)] If $\phi$ is a Poisson $\epsilon$-automorphism of $A_{\Omega}$, then the homological determinant of $\phi$ is $\Pdet(\phi)J(\phi)$.
\end{enumerate}
\end{lemma}

\begin{proof} (1) ``$\Rightarrow$''
Since the Poisson center of $A_{\Omega}$ is $\Bbbk[\Omega]$,
$\phi(\Omega)=a \Omega' +b$ for some $a\in \Bbbk^{\times}$ 
and $b\in \Bbbk$. If $b\neq 0$, then we obtain an $e$-isomorphism
from $P_{\Omega}\to P_{a\Omega'+b}$, which contradicts 
Lemma \ref{xxlem7.6}(3). So $b=0$. 

Let $\{-,-\}_{\Omega}$ (resp. $\{-,-\}_{\Omega'}$) be the 
Poisson bracket of $A_{\Omega}$ (resp. $A_{\Omega'}$).
Now let $f,g\in A_{\Omega}$. We compute by 
\eqref{E0.6.1}
$$\begin{aligned}
e\{\phi(f)&,\phi(y)\}_{\Omega'}
=
\phi(\{f,g\}_{\Omega})=\phi\left(\det
\frac{\partial(f,g,\Omega)}{\partial(x,y,z)}
\right)\\
&=\det
\frac{\partial(\phi(f),\phi(g),\phi(\Omega))}{\partial(\phi(x),\phi(y),\phi(z))}
\\
&=\det
\frac{\partial(\phi(f),\phi(g),\phi(\Omega))}{\partial(\phi(x),\phi(y),\phi(z))}
\det \frac{\partial(\phi(x),\phi(y),\phi(z))}{\partial(x,y,z)} 
J(\phi(x),\phi(y),\phi(z))^{-1}\\
&= \det
\frac{\partial(\phi(f),\phi(g),\phi(\Omega))}{\partial(x,y,z)} 
J(\phi)^{-1}\\
&=\det
\frac{\partial(\phi(f),\phi(g),a\Omega')}{\partial(x,y,z)} 
J(\phi)^{-1}\\
&= a J(\phi)^{-1} \{\phi(f),\phi(g)\}_{\Omega'},
\end{aligned}
$$
which implies that $e=a J(\phi)^{-1}$. The assertion follows.

``$\Leftarrow$'' Partially reverse the above proof. 

(2,3,4,5) Easy consequences of the part (1).
\end{proof}

Note that the isomorphism problem for the class $\{A_{\Omega}\mid {\text{$\Omega$ i.s. potentials}}\}$ 
is solved by Lemma \ref{xxlem7.8}(2).

\begin{lemma}
\label{xxlem7.9}
Let $\Omega$ and $\Omega'$ be i.s. potentials of degree at least $5$.
Let $\phi: P_{\Omega-\xi}\to P_{\Omega'-\xi'}$ be an $e$-isomorphism 
for some $\xi,\xi'\in \Bbbk$. Then $\phi$ lifts uniquely
to an $e$-isomorphism from $A_{\Omega}$ to $A_{\Omega'}$, 
still denoted by $\phi$. Furthermore, the following hold.
\begin{enumerate}
\item[(1)]
In this case, $\phi(\Omega)=a\Omega'$ where $a=eJ(\phi)\in 
\Bbbk^{\times}$ and $\xi=a\xi'$ at the level of $\phi: A_{\Omega}\to 
A_{\Omega'}$.
\item[(2)]
If $\phi$ is a Poisson algebra isomorphism, then 
$\phi(\Omega)=a\Omega'$ where $a=J(\phi)\in \Bbbk^{\times}$ 
and $\xi=a\xi'$ at the level of $\phi: A_{\Omega}\to 
A_{\Omega'}$.
\item[(3)]
If $\phi$ is a Poisson algebra automorphism of $A_{\Omega}$, then 
$\phi(\Omega)=a\Omega$ where $a=J(\phi)\in \Bbbk^{\times}$ and 
$\xi=a\xi$ at the level of $\phi: A_{\Omega}\to 
A_{\Omega}$.
\end{enumerate}
\end{lemma}

\begin{proof}
First, by applying \Cref{xxlem7.5}(2) to both $\phi$ and $\phi^{-1}$ extended to the Poisson fraction fields, we have $\deg \Omega=\deg \Omega'$. Secondly, we claim that $\phi$ lifts to a graded algebra 
isomorphism from $A_{\Omega}\to A_{\Omega'}$. By Lemma 
\ref{xxlem7.6}(1), $\phi$ is a $\Bbbk$-linear automorphism 
of $\Bbbk x+\Bbbk y+\Bbbk z$. So $\phi$ lifts to a graded 
algebra isomorphism from $A_{\Omega}$ to $A_{\Omega'}$, still 
denoted by $\phi$. It remains to show that $\phi$ preserves
the Poisson structure (or $e$-morphism structure). This 
means that $\phi$ maps a Poisson relation such as 
$\{x,y\}=\Omega_z$ to 
\begin{equation}
\label{E7.9.1}\tag{E7.9.1}
e\{\phi(x),\phi(y)\}=\phi(\Omega_z).
\end{equation}
Note that \eqref{E7.9.1} holds in $P_{\Omega'-\xi'}$ and has degree less than $\deg \Omega$. Since each element in $(\Omega'-\xi')
A_{\Omega'}$ has degree at least $\deg \Omega$, \eqref{E7.9.1}
holds in $A_{\Omega'}$. Therefore, we have proved the claim. 

(1) By Lemma \ref{xxlem7.8}(1) (or since $\phi$ is graded), 
$\phi(\Omega)=a \Omega'$ where $a=e J(\phi)$. Since 
$\phi$ is an algebra isomorphism, it maps $\Omega-\xi$ to 
some scalar multiple of $\Omega'-\xi'$. Thus $\phi(\Omega-\xi)=
a(\Omega'-\xi')$. Consequently, $\xi=a\xi'$. The proofs of (2,3)
are similar.
\end{proof}

To summarize, the isomorphism problem (resp. embedding problem) for the class 
\[\{P_{\Omega}\mid {\text{$\Omega$ i.s. potentials of degree $\ge 5$}}\}\] 
is solved by Lemma \ref{xxlem7.9}(2) (resp. by Lemmas \ref{xxlem7.6}(2)
and \ref{xxlem7.9}(2)). The isomorphism problem (resp. embedding problem) 
for the class 
\[\{Q(P_{\Omega})\mid {\text{$\Omega$ i.s. potentials of degree $\ge 5$}}\}\] 
is solved by Corollary \ref{xxcor7.7} and Lemma \ref{xxlem7.9}(2) (resp. 
by Corollary \ref{xxcor7.7} and Lemmas \ref{xxlem7.6}(2) and 
\ref{xxlem7.9}(2)).

\section{Applications}
\label{xxsec8}
In the following subsections, we present more applications of Poisson valuations. Firstly, we will examine the automorphism problem.

\subsection{Automorphism problem}
\label{xxsec8.1}
In this paper, the automorphism problem is to compute the 
automorphism group of Poisson fields/algebras.

Throughout this subsection, we consider a general i.s. potential of degree $\geq 5$.

\begin{theorem}
\label{xxthm8.1}
Let $\Omega$ be an i.s. potential of degree $\geq 5$.
\begin{enumerate}
\item[(1)]
Let $P_{\Omega}$ be as defined in Construction 
\ref{xxcon0.6}. Then 
$$\Aut_{Poi}(Q(P_{\Omega}))=\Aut_{Poi}(P_{\Omega}).$$
Furthermore, every automorphism of $P_{\Omega}$ is Adams graded.
\item[(2)]
Let $P_{\Omega-\xi}$ be as defined in Construction 
\ref{xxcon0.6} where $\xi\in \Bbbk^{\times}$. Then 
$$\Aut_{Poi}(Q(P_{\Omega-\xi}))=\Aut_{Poi}(P_{\Omega-\xi}).$$
\end{enumerate}
\end{theorem}

\begin{proof}
(1) Assertions follow from Lemmas \ref{xxlem7.5}(1) and \ref{xxlem7.6}(2, 4).

(2) It follows from Lemmas \ref{xxlem7.5}(1) and \ref{xxlem7.6}(2).
\end{proof}

\begin{theorem}
\label{xxthm8.2}
Let $\Omega$ be an i.s. potential of degree $\geq 5$. Then 
$$\Aut_{Poi}(Q(A_{\Omega}))=\Aut_{Poi}(A_{\Omega})=
\Aut_{Poi}(P_{\Omega})$$
and each automorphism of $A_{\Omega}$ is Adams graded.
Consequently, for any $\xi\neq 0$, there is an exact sequence of automorphism groups
\[
1\to \Aut_{Poi}(Q(P_{\Omega-\xi}))\to \Aut_{Poi}(Q(A_{\Omega}))\xrightarrow{J(-)}GL_1(\kk)
\]
where $\Aut_{Poi}(Q(A_{\Omega}))$ is a finite 
subgroup of $GL_3(\Bbbk)$ of order bounded above by 
$42(\deg \Omega)(\deg \Omega-3)^2$. 
\end{theorem}

\begin{proof}
It is clear that $\Aut_{Poi}(A_{\Omega})\subseteq 
\Aut_{Poi}(Q(A_{\Omega}))$. By Theorem \ref{xxthm4.9}, we have 
$\,^1\Gamma (Q(A_{\Omega}))=A_{\Omega}$, which implies that 
$\Aut_{Poi}(Q(A_{\Omega}))\subseteq \Aut_{Poi}(A_{\Omega})$ by 
restriction. Therefore $\Aut_{Poi}(Q(A_{\Omega}))
=\Aut_{Poi}(A_{\Omega})$.

Next we prove that $\Aut_{Poi}(A_{\Omega})=\Aut_{Poi}(P_{\Omega})$ 
and every Poisson automorphism of $A_{\Omega}$ is Adams graded. For any Poisson 
automorphism $\phi$ of $A_{\Omega}$, let $\phi'$ denote the induced Poisson automorphism of $P_{\Omega}$ by \Cref{xxlem7.8}(1). By Lemma \ref{xxlem7.6}(4), $\phi'$ preserves the grading and we can uniquely lift 
$\phi'$ to an Adams graded Poisson automorphism of $A_{\Omega}$, 
denoted by $\sigma$, by \Cref{xxlem7.9}(3). Clearly $\sigma'=\phi'$. Let 
$\varphi= \phi\circ \sigma^{-1}$. Then $\varphi'=Id_{P_{\Omega}}$. 
Since $\sigma$ preserves the Adams grading, it suffices to show 
that $\varphi$ is the identity. Since $\varphi'$ is the identity, 
we have
\begin{align}
\notag
\varphi(x)&=x+\Omega q_1,\\
\notag
\varphi(y)&=y+\Omega q_2,\\
\notag
\varphi(z)&=z+\Omega q_3
\end{align}
for some $q_1,q_2,q_3\in A_{\Omega}$. Since $\Omega$ does not have 
a linear term, a computation will show that $\varphi(\Omega)=
\Omega+ \Omega \alpha(q_1,q_2,q_3)$ where $\alpha(q_1,q_2,q_3)\in 
(A_{\Omega})_{\geq 1}$. Since $\varphi$ preserves the Poisson 
center $\Bbbk[\Omega]$ of $A_{\Omega}$ \cite[Proposition 1]{UZ}, 
$\varphi(\Omega)=\Omega+\Phi$ where 
$\Phi\in \Omega \Bbbk[\Omega]$. Since $\varphi$ is an 
algebra automorphism of $\Bbbk[\Omega]$, $\Phi=0$ and
$\varphi(\Omega)=\Omega$.

Let ${\mathbb B}:=\{1, x,y,z\} \cup \{b_s\}$ be a 
fixed $\Bbbk$-linear basis of $P_{\Omega}$ consisting of 
a set of monomials $x^{s_1}y^{s_2} z^{s_3}$. For example,
if $\Omega=x^{n+3}+y^{n+3}+z^{n+3}$, then we can choose
$${\mathbb B}=\{x^{s_1}y^{s_2} z^{s_3}
\mid 
{\text{$s_i\geq 0$ for $i=1,2$ and $0\leq s_{3}
\leq n+2$}}\}.$$ 
We consider ${\mathbb B}$ as a fixed subset of monomial elements 
in $A_{\Omega}$ in a canonical way. By the induction on the 
degree of elements, every element $f$ in $A_{\Omega}$ is of the 
form 
$$f=1 f^{1}(\Omega)+x f^{x}(\Omega)+y f^{y}(\Omega)
+z f^{z}(\Omega)+\sum_{b_s} b_s f^{b_s}(\Omega)$$ 
where each $f^{\ast}(\Omega)$ is in $\Bbbk[\Omega]$. By 
re-cycling the letters, we write
$$\varphi(x)=1 f^{1}(\Omega)+x f^{x}(\Omega)
+y f^{y}(\Omega)+z f^{z}(\Omega)
+\sum_{b_s} b_s f^{b_s}(\Omega)$$
for some $f^{\ast}(\Omega)\in \Bbbk[\Omega]$.

For each $\xi\neq 0$, let $\pi_{\xi}$ be the canonical quotient 
map from $A_{\Omega} \to P_{\Omega-\xi}$. It is clear that the 
image of ${\mathbb B}$ is a $\Bbbk$-linear basis of 
$P_{\Omega-\xi}$. For simplicity, we continue to use $b_s$ etc 
for $\Bbbk$-linear basis elements in $P_{\Omega-\xi}$.

Let $\varphi'_{\xi}$ be the induced automorphism of 
$P_{\Omega-\xi}$. Then $\varphi'_{\xi}$ is a Poisson 
algebra automorphism of $P_{\Omega-\xi}$ and
$$\varphi'_{\xi}(x)=1f^{1}(\xi)+x f^{x}(\xi)
+y f^{y}(\xi)+z f^{z}(\xi)+\sum_{b_s} b_s f^{b_s}(\xi).$$
By Lemma \ref{xxlem7.6}(1), $\phi'_{\xi}$ is linear.
Thus $f^{b_s}(\xi)=0$ for all $\xi\neq 0$. Hence 
$f^{b_s}(\Omega)=0$, consequently, $\varphi(x)=1f^{1}(\Omega)
+x f^{x}(\Omega)+y f^{y}(\Omega)+z f^{z}(\Omega)$.
Since $\varphi(x)=x+\Omega q_1$, $\varphi(x)=x+ x\Omega w_{11}(\Omega)
+y\Omega w_{12}(\Omega)+z\Omega w_{13}(\Omega)$. Similarly,
one has $\varphi(y)=y+ x\Omega w_{21}(\Omega)
+y\Omega w_{22}(\Omega)+z\Omega w_{23}(\Omega)$
and $\varphi(z)=z+ x\Omega w_{31}(\Omega)
+y\Omega w_{32}(\Omega)+z\Omega w_{33}(\Omega)$
for some polynomials $w_{ij}(t)\in \kk[t]$.
Now by Lemma \ref{xxlem8.3} below, $\varphi$ is the identity
as required.

Finally, we show the exact sequence of automorphism groups and the order bound. By Lemmas \ref{xxlem7.6}(1) and \ref{xxlem7.9}(3) and Theorem \ref{xxthm8.1}(2), we have the following natural identifications
\begin{align*}
   \Aut_{Poi}(Q(P_{\Omega-\xi}))&=\Aut_{Poi}(P_{\Omega-\xi})\\
   &=\{\phi\in \Aut_{gr.alg}(\kk[x,y,z])\mid \phi(\Omega)=\Omega,\ J(\phi)=1\}\\
    &\subset \{\phi\in \Aut_{gr.alg}(\kk[x,y,z])\mid \phi(\Omega)=J(\phi)\,\Omega\}\\
    &=\Aut_{Poi}(P_{\Omega})\\
    &=\Aut_{Poi}(Q(A_{\Omega})).
\end{align*}
Hence, we can view $\Aut_{Poi}(Q(P_{\Omega-\xi}))$ as a normal subgroup of $\Aut_{Poi}(Q(A_{\Omega}))$ consisting of those Poisson automorphisms whose Jacobian determinant is trivial. So, the exact sequence follows. Now, for the order bound, after replacing $\Bbbk$ by its algebraic closure, we may assume that $\Bbbk$ is algebraically
closed. Since every Poisson automorphism of $P_{\Omega}$ is 
graded, $G:=\Aut_{Poi}(P_{\Omega})$ is a subgroup of 
$GL_3(\Bbbk)$. It remains to show that $|G|\leq 42 d (d-3)^2$ 
where $d:=\deg \Omega$. Let $B$ be the graded commutative algebra 
$P_{\Omega}$ which is $A/(\Omega)$ and let $X$ be the 
corresponding smooth curve $\Proj B$. Then the genus $g$ of $X$ 
is $\frac{1}{2}(d-1)(d-2)$ by the genus-degree formula. By the 
Hurwitz's automorphism theorem \cite{Hu}, the order of 
the group $\Aut(X)$ is bounded by $84(g-1)$. There is a natural 
group homomorphism $\pi$ from $\Aut_{gr.alg}(B)\to \Aut(X)$ 
whose kernel consists of automorphisms of the form 
$\eta_{\xi}: B\to B, b\mapsto \xi^{\deg b} b$ for some 
$\xi\in \Bbbk$. Thus, we have an exact sequence
$$1\to \ker \pi\cap G \to 
G \to \Aut(X).$$
One can check that $\ker \pi\cap G=\{\eta_{\xi}
\mid \xi^{d-3}=1\}$. Thus $|\ker \pi\cap G|= d-3$. Now
$$|G|\leq |\ker \pi\cap G| |\Aut(X)|
\leq (d-3)84 (g-1)=42d(d-3)^2.$$
\end{proof}

\begin{lemma}
\label{xxlem8.3}
Suppose $\Omega$ is an i.s. potential of degree $n\geq 3$. 
\begin{enumerate}
\item[(1)]
Let $f_1,f_2,f_3$ be elements in $\Bbbk x+\Bbbk y+\Bbbk z$
such that $f_1 \Omega_x+f_2\Omega_y+f_3\Omega_z=0$. Then
$f_i=0$ for all $i=1,2,3$.
\item[(2)]
Let $\varphi$ be a Poisson automorphism of $A_{\Omega}$ such 
that $\varphi(\Omega)=\Omega$ and 
$$\begin{aligned}
\varphi(x)&=x+ x\Omega w_{11}(\Omega)
+y\Omega w_{12}(\Omega)+z\Omega w_{13}(\Omega)\\
\varphi(y)&=y+ x\Omega w_{21}(\Omega)
+y\Omega w_{22}(\Omega)+z\Omega w_{23}(\Omega)\\
\varphi(z)&=z+ x\Omega w_{31}(\Omega)
+y\Omega w_{32}(\Omega)+z\Omega w_{33}(\Omega)
\end{aligned}
$$
for some polynomials $w_{ij}(t)$. Then 
$\varphi=Id$.
\end{enumerate}
\end{lemma}

\begin{proof}
(1) In this case 
$(f_1,f_2,f_3)$ is in the kernel of the map
$\overrightarrow{\nabla}\Omega \cdot$ in the Koszul
complex
{\begin{equation}
\notag
\begin{tabular}{ccccc}
&\,$A[-2n+2]$ & &$A[1-n]$& \vspace*{-2mm}\\
$0\to A[-3n+3]\xrightarrow{\overrightarrow{\nabla}\Omega}$
& \hspace*{-2.5mm}$\oplus A[-2n+2]$
&\hspace*{-2.5mm}$\xrightarrow{\overrightarrow{\nabla}\Omega\times }$
& \hspace*{-2.5mm}$ \oplus A[1-n]$
&\hspace*{-2.5mm}$\xrightarrow{\overrightarrow{\nabla}\Omega\cdot }A
\to A/(\Omega_x,\Omega_y,\Omega_z)\hspace*{-1mm}\to \hspace*{-.8mm}0.$\\
&\hspace*{-2.5mm}$\oplus A[-2n+2]$& & \hspace*{-2.5mm}$\oplus A[1-n]$
&
\end{tabular}
\end{equation}}
Since $\Omega$ has an isolated singularity at the origin, the 
above complex is exact by \cite[Proposition 3.5]{Pi}. Since
$f_i$ has degree 1 and $\deg \Omega\geq 3$, $(f_1,f_2,f_3)$ 
is not in the image of $\overrightarrow{\nabla}\Omega\times$
if it is nonzero. Therefore $f_i=0$ for all $i=1,2,3$

(2) We need to prove that $w_{ij}(t)=0$ for all $i,j$. 
Suppose to the contrary that some $w_{ij}\neq 0$. Then we can 
write
$$\begin{aligned}
\varphi(x)&=x+ f_1 \Omega^{s}+hdt\\
\varphi(y)&=y+ f_2 \Omega^{s}+hdt\\
\varphi(z)&=z+ f_3 \Omega^{s}+hdt
\end{aligned}
$$
where $f_1,f_2,f_3\in \Bbbk x+\Bbbk y+\Bbbk z$ are not all 
zero, $s\geq 1$, and $hdt$ stands for linear combination
of higher Adams degree terms. Through Taylor expansion,
$$\begin{aligned}
\Omega&=\phi(\Omega)\\
&=\Omega(x+f_1 \Omega^s+ hdt, y+f_2 \Omega^s+hdt,
z+f_3 \Omega^s+hdt)\\
&=\Omega+\Omega_x(f_1\Omega^s)+\Omega_y(f_2\Omega^s)
+\Omega_z (f_3 \Omega^s)+hdt
\end{aligned}
$$
which implies that $\Omega_x f_1+\Omega_y f_2
+\Omega_z f_3=0$. By part (1), $f_i=0$ for all $i=1,2,3$, is a contradiction.
\end{proof}

It is worth pointing out that the automorphism group of the Poisson field $Q(P_{\Omega})$ is closely related to the automorphism group of the projective curve $X=\Proj(A/(\Omega))$, and usually computable. To illustrate this fact, we consider a particular case where $\Omega=x^{d}+y^{d}+z^{d}=0$ is the Fermat curve for some $d\geq 5$.

Let 
$$G(0,d):=\{(a,b,c)\in (\Bbbk^{\times})^{3} \mid
ab=c^{d-1}, bc=a^{d-1}, ac=b^{d-1}\}$$
and
$$G(1,d):=\{(a,b,c)\in (\Bbbk^{\times})^{3} \mid
ab=c^{d-1}, bc=a^{d-1}, ac=b^{d-1}, abc=1\}.$$
Suppose $\kk$ is algebraically closed. Denote by $C_n$ the cyclic group of order $n$. One can check directly that 
\[G(0,d)=\{(uv,\frac{v^{d-2}}{u},v)\in (\kk^\times )^3\mid u^d=1,v^{d^2}=v^{3d}\}\cong  C_d\times C_{d(d-3)}
\]
and 
\[G(1,d)=\{(u,\frac{1}{uv},v)\in (\kk^\times )^3\mid u^d=v^d=1\}\cong C_d\times  C_{d},\] 
are finite groups.

\begin{proposition}
  \label{xxpro8.4}
Suppose $\kk$ is algebraically closed. Let 
$\Omega=x^{d}+y^{d}+z^{d}$ where $d\geq 5$. Let 
$P_{\Omega}$ and $P_{\Omega-\xi}$ be defined as in Construction 
\ref{xxcon0.6} where $\xi\in \Bbbk^\times$.
\begin{enumerate}
\item[(1)]
There is a short exact sequence of groups
$$1\to C_{d}\times C_{d(d-3)}\to \Aut_{Poi}(P_{\Omega})\to S_3\to 1.$$ 
Moreover, $d$ is even if and 
only if $\Aut_{Poi}(P_{\Omega})\cong S_3\ltimes (C_{d}\times C_{d(d-3)})$. 
\item[(2)]
If $d$ is odd, there is a short exact sequence of groups
$$1\to C_d\times C_d\to \Aut_{Poi}(P_{\Omega-\xi})\to C_3\to 1$$
and $\Aut_{Poi}(P_{\Omega-\xi})\cong C_3\ltimes (C_d\times C_d)$.
\item[(3)]
If $d$ is even, there is a short exact sequence of groups
$$1\to C_d\times C_d\to \Aut_{Poi}(P_{\Omega-\xi})\to S_3\to 1$$
and $\Aut_{Poi}(P_{\Omega-\xi})\cong S_3\ltimes (C_d\times C_d)$.
\end{enumerate}  
\end{proposition}
\begin{proof} 
(1) Let $P=P_{\Omega}$. By Lemma \ref{xxlem7.6}(4), every Poisson 
automorphism of $P$ is graded and hence $\Aut_{Poi}(P)$ is a subgroup of $GL_3(\kk)$. Let $X=\Proj(A/(\Omega))$. In the proof of Theorem \ref{xxthm8.2}, we have an exact sequence of groups 
\[
1\to \mu_{d-3} \to \Aut_{Poi}(P)\xrightarrow{\pi}\Aut(X)
\]
where $\mu_{d-3}$ is the multiplicative subgroup of $\kk^\times$ consisting of all the $(d-3)$th roots of unity. For the rest of the proof of part (1), we write $(x,y,z)$ as $(x_1,x_2,x_3)$. It is well-known (e.g.,  \cite{ODR,Tz}) that $\Aut(X)=\langle \phi_1,\phi_2,\phi_3,\phi_4\rangle\cong S_3\ltimes (C_d\times C_d)$ such that
\begin{align*}
    \phi_1([x_1:x_2:x_3])&=[x_2: x_1: x_3],\\
    \phi_2([x_1:x_2:x_3])&=[x_3: x_1: x_2],\\
    \phi_3([x_1:x_2:x_3])&=[\omega\, x_1: x_2: x_3],\\
    \phi_4([x_1:x_2:x_3])&=[x_1:\omega\, x_2: x_3]
\end{align*}
where $\omega$ is a primitive $d$th root of unity. Let $\phi$ be a Poisson algebra automorphism of $P$. Since $\pi(\phi)$ can be written in terms of $\phi_i$ for $1\le i\le 4$, it is clear that $\phi$ preserves $\Bbbk x_1\cup 
\Bbbk x_2\cup \Bbbk x_3$. Hence there are $\sigma\in 
S_3$ and $(a_1,a_2,a_3)\in (\Bbbk^{\times})^3$ such that
\begin{equation}
\label{E8.1.1}\tag{E8.4.1}
\phi(x_1)=a_1 x_{\sigma(1)}, \quad
\phi(x_2)=a_2 x_{\sigma(2)}, \quad
\phi(x_3)=a_3 x_{\sigma(3)}.
\end{equation}
Applying $\phi$ to the Poisson bracket
$\{x_i,x_j\}=d x_k^{d-1}$ where $(i,j,k)=(1,2,3)$,
or $(2,3,1)$ or $(3,1,2)$, we obtain that 
\begin{equation}
\label{E8.1.2}\tag{E8.4.2}
a_1 a_2 ={\rm{sgn}}(\sigma) a_3^{d-1},\quad
a_2 a_3 ={\rm{sgn}}(\sigma) a_1^{d-1},\quad
a_3 a_1 ={\rm{sgn}}(\sigma) a_2^{d-1}.
\end{equation}
If $\phi$ preserves $\Bbbk x_1$, $\Bbbk x_2$ and $\Bbbk x_3$ 
individually, then $\sigma$ is the identity. In this case 
\eqref{E8.1.2} is equivalent to $(a_1,a_2,a_3)\in G(0,d)$. By 
an elementary computation and the hypothesis on $\Bbbk$ for 
every $\sigma \in S_3$, \eqref{E8.1.2} has a solution 
$(a_1,a_2,a_3)$ with $d\ge 5$. Combining these facts, $G(0,d)$ 
is a normal subgroup of $\Aut_{Poi}(P)$ and 
$\Aut_{Poi}(P)/G(0,d) \cong S_3$. Therefore, the main assertion in part (1) is proved.

Suppose $d$ is even. For every $\sigma\in S_3$, let 
$\sigma'$ be the automorphism of $P$ defined by 
$\sigma'(x_i)={\rm{sgn}}(\sigma)x_{\sigma(i)}$. One can check 
by \eqref{E8.1.2} that $\sigma'$ is a Poisson algebra 
automorphism of $P$. Moreover, the subgroup of 
$\Aut_{Poi}(P)$ generated by 
$\{\sigma'\,|\, \sigma\in S_3\}$ is isomorphic to $S_3$. 
Therefore $\Aut_{Poi}(P_{\Omega})\cong S_3\ltimes G(0,d)$.

Conversely, suppose $\Aut_{Poi}(P)\cong S_3 \ltimes 
G(0,d)$. This means that $S_3$ is a subgroup of 
$\Aut_{Poi}(P)$. Let $\sigma=(12)$ 
(${\rm{sgn}}(\sigma)=-1$) and suppose
$\varphi\in \Aut_{Poi}(P)$ corresponds to $\sigma$. 
Then there are $a_1,a_2,a_3$ such that
$\varphi(x_1)=a_1 x_2$ and $\varphi(x_2)=a_2 x_1$ and 
$\varphi(x_3)=a_3 x_3$. Since $\varphi^2=Id$, we 
have $a_1a_2=1=a_3^2$. If $d$ is odd, an elementary 
computation shows that there is no $(a_1,a_2,a_3)$ such 
that both $a_1a_2=1=a_3^2$ and \eqref{E8.1.2} hold. This
yields a contradiction. Therefore, $d$ must be even.

(2) By Theorem \ref{xxthm8.2}, $\Aut_{Poi}(P_{\Omega-\xi})$ is a normal subgroup of $\Aut_{Poi}(P_{\Omega})$ consisting of those automorphisms whose Jacobian determinant is trivial. Hence, its proof is similar to the proof of part (1), except that
we need to show that there is no Poisson algebra automorphism 
$\varphi$ such that $\varphi(x)=ay$ and $\varphi(y)=bx$ and 
$\varphi(z)=cz$ for some $a,b,c\in \Bbbk^{\times}$. Suppose 
to the contrary that such an automorphism exists. Then 
we have 
$$-ab=c^{d-1},\ -bc=a^{d-1},\ -ac=b^{d-1},\ abc=-1.$$
Since $d$ is odd, the above system of equations has
no solution. Therefore, the quotient group 
$\Aut_{Poi}(P_{\Omega-1})/G(1,n)$ is isomorphic to
$C_3$ (not $S_3$). The rest is a routine verification.

(3) The proof is similar to the proof of part (2),
whence it is omitted.

\end{proof}

It is interesting to note that if $\Omega=x^d+y^d+z^d$ where $d\geq5$ is even and $\xi\neq0$, then $\Aut_{Poi}(P_{\Omega-\xi})\cong\Aut(X)$ where $X=\Proj(A/(\Omega))$. We're curious about the following question:
\begin{question}
What type of i.s. potentials $\Omega$ can satisfy $\Aut_{Poi}(P_{\Omega-\xi})\cong\Aut(X)$ with $X=\Proj(A/(\Omega))$?
\end{question}

In the remaining part of this section, we will explore other applications of Poisson valuations. Our goal is to provide essential ideas on 
how to use valuations to solve problems related to Poisson algebras. 
Hence, the presentation of each topic will be very brief. 

\subsection{Dixmier property}
\label{xxsec8.2}

\begin{definition}
\label{xxdef8.5}
Let $P$ be a Poisson algebra. We say $P$ satisfies the 
{\it Dixmier property} if every injective Poisson algebra 
morphism $f: P\to P$ is bijective.
\end{definition}

This property is related to the Dixmier conjecture \cite{Di}, which states 
that every endomorphism of the $n$-th Weyl algebra over a field 
of characteristic zero is bijective. The Dixmier Conjecture is 
stably equivalent to the Jacobian Conjecture 
\cite{BK, BCW, Ts1, Ts2, vdEKC}. By \cite{AvdE}, these two 
conjectures are also equivalent to the Poisson conjecture, which 
states that every Poisson endomorphism of the $n$th canonical 
Poisson algebra $W^{\otimes n}$ over a field of characteristic 
zero is bijective (where $W$ is the Weyl Poisson polynomial ring 
given in Example \ref{xxexa2.10}). The Poisson conjecture 
asserts that $W^{\otimes n}$ satisfies the Dixmier property. Proving the Dixmier property for certain Poisson algebras is undeniably challenging, but it is an invaluable pursuit.

Our main result of this subsection is

\begin{theorem}
\label{xxthm8.6}
Let $\Omega$ be an i.s. potential of degree $\geq 5$. 
Let $P_{\Omega}$ and $P_{\Omega-\xi}$ be defined as in 
Construction \ref{xxcon0.6}.
\begin{enumerate}
\item[(1)]
Both $P_{\Omega}$ and $Q(P_{\Omega})$ satisfy the 
Dixmier property.
\item[(2)]
Both $P_{\Omega-\xi}$ and $Q(P_{\Omega-\xi})$, for
$\xi\in \Bbbk^{\times}$, satisfy the 
Dixmier property.
\end{enumerate}
\end{theorem}

\begin{proof} We only give the proof for part (1). The proof of the part (2) is similar and omitted.

Let $P=P_{\Omega}$. By Lemma \ref{xxlem7.6}(2), $P$ has the 
Dixmier property.

Let $K=Q(P_{\Omega})$ and let $f: K\to K$ be an injective 
Poisson algebra map. By Theorem \ref{xxthm4.9}, $\,^1\Gamma(K)=P$. 
By Lemma \ref{xxlem4.2}(1), the restriction of $f$ on $P$, 
denoted by $f\mid_{P}$, is an injective Poisson algebra map 
$P:=\,^1\Gamma(K)\to \,^1\Gamma(K)=P$. By Lemma \ref{xxlem7.6}(2), 
$f\mid_{P}$ is bijective. Since $K=Q(P)$, $f$ is bijective. 
Thus, $K$ has the Dixmier property. 
\end{proof}

\begin{remark}
\label{xxrem8.7}
It is known that the Weyl Poisson field and the skew Poisson field do not possess the Dixmier property \cite{GZ}. Therefore, $Q(P_{\Omega})$ and $Q(P_{\Omega-\xi})$ are among the first
few Poisson fields that are proved to have the Dixmier property.
\end{remark}

\begin{example}
\label{xxexa8.8}
Let $\Omega=x^3+y^3+z^3+\lambda xyz$. Then $x\to \Omega x$,
$y\to \Omega y$, $z\to \Omega z$ defines an injective
Poisson algebra endomorphism of $A_{\Omega}$ that is not bijective.
So $A_{\Omega}$ does not possess the Dixmier property.
\end{example}

\subsection{Rigidity of grading}
\label{xxsec8.3}

In this subsection, we prove the following. Recall that we say a 
${\mathbb Z}$-graded algebra $A=\bigoplus_{i\in {\mathbb Z}} A_i$ 
is {\it connected graded} if $A_i=0$ for all $i>0$ and $A_0=\Bbbk$.

\begin{theorem}
\label{xxthm8.9}
Let $\Omega$ be an i.s. potential of degree $\geq 4$. Then 
$P_{\Omega}$ has a unique connected grading such that it is 
Poisson $(\deg \Omega-3)$-graded. 
\end{theorem}

\begin{proof}
Let $P$ denote $P_{\Omega}$ and let $n=\deg \Omega-3$. By 
Construction \ref{xxcon0.6}, there is an Adams grading with 
$\deg(x)=\deg(y)=\deg(z)=1$ such that $P$ is Poisson $(-n)$-graded.
To make $P$ $n$-graded, we need to set that $\deg(x)=\deg(y)=
\deg(z)=-1$ by our Definition \ref{xxdef2.1}. The associated
valuation is denoted by $\nu_{old}$ with $\nu_{old}(f)=-1$
for $f=x,y,z$. 

Suppose $P$ has a new connected grading such that $P$ is 
connected graded and Poisson $n$-graded and let $\mu$ be the 
$n$-valuation associated with this new grading as defined by
\eqref{E3.1.1}. In this case $\mu(f)<0$ for all 
$f\in P\setminus\Bbbk$. By Lemma \ref{xxlem4.8}(4), 
$\mu=\nu_{old}$ as given in the previous paragraph. As a 
consequence, $\mu(x)=\mu(y)=\mu(z)=-1$. 

If $x, y$, and $z$ are homogeneous with respect to the new 
grading, then two gradings coincide. It remains to show
that $x,y,z$ are homogenous in this new grading. Write 
$x=x_0+a$, $y=y_0+b$, and $z=z_0+c$ where $x_0,y_0,z_0$ are 
homogeneous of new degree $-1$ and $a,b,c\in \Bbbk$. Since every 
linear combination of $x,y,z$ has $\mu$-value $-1$ 
(as $\mu=\nu_{old}$), $x_0, y_0, z_0$ are linearly 
independent. By Taylor expansion, working in $P_{\Omega}$,
$$\begin{aligned}
0&=\Omega(x,y,z)\\
&=\Omega(x_0,y_0,z_0)+a\Omega_x(x_0,y_0,z_0)+b\Omega_y(x_0,y_0,z_0)
+c\Omega_z(x_0,y_0,z_0)+ hdt
\end{aligned}
$$
where $hdt$ is a linear combination of higher degree terms with respect to the new grading. Since $\Omega_x(x_0,y_0,z_0)$,
$\Omega_y(x_0,y_0,z_0)$, $\Omega_z(x_0,y_0,z_0)$ are linearly 
independent (as $\Omega$ has an isolated singularity at origin), 
$a=b=c=0$. So $x,y,z$ are homogeneous as required.
\end{proof}

\subsection{Rigidity of filtration}
\label{xxsec8.4}

In this subsection, we prove the following.

\begin{theorem}
\label{xxthm8.10}
Let $\Omega$ be an i.s. potential of degree $\geq 4$ and 
$\xi\neq 0$. Then $P_{\Omega-\xi}$ has a unique filtration 
${\mathbb F}$ such that the associated graded ring 
$\gr_{\mathbb F} (P_{\Omega-\xi})$ is a connected graded 
Poisson $(\deg \Omega-3)$-graded domain. 
\end{theorem}

\begin{proof}
Let $P$ be $P_{\Omega-\xi}$ and $n=\deg \Omega-3$. Let 
${\mathbb F}^{c}$ be the original filtration by the construction 
and let $\nu^{c}$ be the associated valuation given in Lemma 
\ref{xxlem3.2}(3). In particular, $\nu^c(f)=-1$ for $f=x,y,z$. 
Since $\gr_{\nu^c}(P)\cong P_{\Omega}$ which is Poisson
$n$-graded, $\nu^{c}$ is an $n$-valuation of $P$.

Suppose $P$ has another filtration, say ${\mathbb F}$, such 
that $\gr_{\mathbb F} P$ is a connected graded Poisson 
$n$-graded domain. Let $\nu$ be the associated $n$-valuation 
of $P$. In this case $\nu(f)<0$ for some $f\in P$. By Lemma 
\ref{xxlem4.8}(4), $\nu=\nu^{c}$. Lemma \ref{xxlem2.6}, 
${\mathbb F}={\mathbb F}^c$ as required.
\end{proof}

\subsection{Proof of statements in Table 1}
\label{xxsec8.5}
For the Weyl Poisson field $K_{Weyl}$, the valuations constructed
in Example \ref{xxexa2.10}, $\{\nu_{\xi}\mid \xi\in \Bbbk\}$, are
distinct faithful $0$-valuations (when we choose $w=0$). If $\Bbbk$
is uncountable, we have obtained uncountably many faithful 
$0$-valuations.

For the skew Poisson field $K_q$, the statement follows from 
Theorem \ref{xxthm6.2}. For $Q(P_{\Omega})$ where 
$\Omega=x^3+y^3+z^3+\lambda xyz$ and $\lambda^3\neq -3^3$,
the statement follows from Theorem \ref{xxthm3.8}(1).
For $Q(P_{\Omega-1})$ where 
$\Omega=x^3+y^3+z^3+\lambda xyz$ and $\lambda^3\neq -3^3$,
the statement follows from Theorem \ref{xxthm3.11}(1).

For $\Omega$ of degree $\geq 4$ with an isolated singularity,
the statement follows from Lemma \ref{xxlem4.8}(3,5).

\subsection*{Acknowledgments.} 
The authors are grateful to the referee for carefully reading the article and providing helpful suggestions. The authors thank Ken Goodearl and Milen Yakimov for many valuable conversations and correspondences on the subject and thank S{\' a}ndor Kov{\' a}cs for the proof of 
Lemma \ref{xxlem5.5}(1). Wang was partially supported by 
Simons collaboration grant \#688403 and Air Force Office of 
Scientific Research grant FA9550-22-1-0272. Zhang was 
partially supported by the US National Science Foundation 
(No. DMS-2001015 and DMS-2302087). Part of this research 
work was done during the first three authors' visits to the Department of Mathematics at the University of 
Washington in June 2022 and January 2023. They are 
grateful for the fourth author's invitation and wish to 
thank the University of Washington for its hospitality.

\providecommand{\bysame}{\leavevmode\hbox to3em{\hrulefill}\thinspace}
\providecommand{\MR}{\relax\ifhmode\unskip\space\fi MR }
\providecommand{\MRhref}[2]{%

\href{http://www.ams.org/mathscinet-getitem?mr=#1}{#2} }
\providecommand{\href}[2]{#2}

\end{document}